\documentclass[a4paper,10pt]{article}
\usepackage[a4paper, margin=1.2in]{geometry}
\usepackage{mathrsfs,amsfonts,latexsym,theorem,amssymb,amsmath,newlfont}
\usepackage{graphicx}
\newcommand{\tr}{^\mathsf{T}}
\usepackage{pgf}
\usepackage{fancyhdr}
\usepackage{times}
\usepackage{dsfont} 
\usepackage{hyperref}

\theoremstyle{plain}

\newtheorem{thm}{Theorem}[section]
\newtheorem{defn}[thm]{Definition}
\newtheorem{lem}[thm]{Lemma}
\newtheorem{prop}[thm]{Proposition}

\newtheorem{remark}[thm]{Remark}

\newtheorem{hypo}[thm]{Hypothesis}
\numberwithin{equation}{section}

\def\beginpf{\noindent {\bf Proof:} \quad}
\def\endpf{\rightline{$\square$}}

\def\HH{{\mathcal H}}

\def\CC{{\mathbb C}}

\def\NN^*{{\mathbb N}^*}
\def\NN{{\mathbb N}}

\def\RR{{\mathbb R}}

\def\HH{{\mathcal H}}
\def\ZZ{{\mathbb Z}}

\def\B{{\cal B}}

\def\<{\langle}
\def\>{\rangle}

\newcommand{\lag}{\langle}
\newcommand{\rag}{\rangle}
\newcommand{\A}{\mathcal{A}}
\renewcommand{\S}{\mathcal{S}}
\newcommand{\md}{\mathrm{d}}
\DeclareMathOperator{\Span}{span}

\allowdisplaybreaks[3]        

\date{\empty}
\title{Rapid Stabilization of a Linearized Bilinear 1-D Schr\"odinger Equation}
\author{Jean-Michel \textsc{Coron}\footnote{Universit\'e Pierre et Marie Curie, LJLL UMR 7598, 4 place Jussieu, 75005, Paris, France.} , Ludovick \textsc{Gagnon}\footnote{Universit\'e de Nice Sophia-Antipolis, CNRS UMR 7351, Laboratoire J.-A. Dieudonn\'e, Parc Valrose, 06108, Nice, France.} , Morgan \textsc{Morancey}\footnote{Aix Marseille Universit\'e, CNRS, Centrale Marseille, I2M UMR 7373
13453, Marseille, France.}}

\begin{document}

\maketitle

\begin{abstract}

We consider the one dimensional Schr\"odinger equation with a bilinear control and prove the rapid stabilization of the linearized equation around the ground state. The feedback law ensuring the rapid stabilization is obtained using a transformation mapping the solution to the linearized equation on the solution to an exponentially stable target linear equation. A suitable condition  is imposed on the transformation in order to cancel the non-local terms arising in the kernel system. This conditions also insures the uniqueness of the transformation.   The continuity and invertibility of the transformation follows from exact controllability of the linearized system. 



\end{abstract}

\section{Introduction}

\subsection{Main result}

Let $T>0$. Consider the Schr\"odinger equation
\begin{equation}\label{nl}
\begin{cases}
i\partial_t \psi = -\Delta \psi - u(t)\mu(x)\psi, &\quad (t,x)\in (0,T)\times (0,1),
\\
\psi(t,0)=\psi(t,1)=0, &\quad t\in (0,T).
\end{cases}
\end{equation}
In (\ref{nl}), $\psi$ is the complex-valued wave function, of $L^2$-norm 1, of a particle confined in a $1-D$ infinite square potential well. The particle is subjected to an electric field inside of the domain, where $u\in L^2((0,T);\RR)$ is the amplitude of the electric field and $\mu \in H^3((0,1);\RR)$ is the dipolar moment of the particle.

\noindent
Before stating our main result, we set some notations. Let $\A : D(\A)\subset L^2((0,1);\CC)\to L^2((0,1);\CC)$ be defined by
\begin{equation} \label{DefOperateurA}
\A\phi:= -\Delta \phi, \quad D(\A):=H^2((0,1);\CC)\cap H^1_0((0,1);\CC).
\end{equation}
The eigenvalues and eigenfunctions of $\A$ are given by
\[
\lambda_k:=(k\pi)^2,\quad \varphi_k(x):=\sqrt{2}\sin(k\pi x), \quad k\in \NN^*.
\]
The eigenstates of \eqref{nl} ($u=0$) are given by $\Phi_k(t,x):=e^{-i\lambda_k t}\varphi_k(x)$. The eigenstate $\Phi_1$ associated to the smallest eigenvalue is called the ground state.

\noindent
Define the space $H^s_{(0)}((0,1);\CC):=D(\A^{s/2})$ equipped with the inner product
\[
\<\phi,\psi\>_{H^s_{(0)}}:=\sum_{k=1}^{+\infty}\lambda_k^s \<\phi,\varphi_k \> \overline{\<\psi,\varphi_k \>} ,
\]
where $\< \cdot, \cdot \>$ is the $L^2((0,1);\CC)$-inner product. The space $H^s_{(0)}$ is endowed with the $\|.\|_{H^s_{(0)}}$-norm associated to the $H^s_{(0)}$-inner product.
We underline that the spaces used in this article can also be explicitly described by
$H^2_{(0)}((0,1);\CC)= H^2 \cap H^1_0 ((0,1);\CC)$ and
\begin{align*}
H^3_{(0)}((0,1);\CC)=& \left\{ \phi \in H^3 \cap H^1_0((0,1); \CC) \: ; \: \phi''(0)=\phi''(1)=0 \right\}, \\
H^5_{(0)}((0,1);\CC)=& \left\{ \phi \in H^5 \cap H^3_{(0)}((0,1); \CC) \: ; \: \phi^{(4)}(0)=\phi^{(4)}(1)=0 \right\}.
\end{align*}
We denote by $\S$ the radius $1$ sphere of $L^2((0,1); \CC)$.

Throughout this article, we assume that $\mu$ satisfies the following assumption.
\begin{hypo}\label{hyp}
The function $\mu$ belongs to $H^3((0,1);\RR)$ and there exists $c>0$ satisfying
\begin{equation}
\label{esthyp1}
\left|\<\mu\varphi_1,\varphi_k\>\right| \geq \dfrac{c}{k^3}, \quad \forall k\in \NN^*.
\end{equation}
\end{hypo}
\noindent
\begin{remark}\label{behavehyp}
A direct computation shows that, for every function $\mu \in H^3((0,1);\RR)$, we have
\begin{align}
\label{expressionmuphiphikresteintegral}
\<\mu \varphi_1,\varphi_k\>=\dfrac{4}{(k\pi)^3}\left((-1)^{k+1}\mu'(1)-\mu'(0)\right)-
\dfrac{\sqrt{2}}{(k\pi)^3}\int_0^1 (\mu \varphi_1)'''(x)\cos(k\pi x) dx
\end{align}
and therefore there exists $C>0$ such that
\begin{equation}
\label{boundmuphi1phik}
\left|\<\mu\varphi_1,\varphi_k\>\right| \leq \dfrac{C}{k^3}, \quad \forall k\in \NN^*.
\end{equation}
Moreover, since
\begin{equation}\label{integraletendvers0}
\lim_{k\rightarrow +\infty} \int_0^1 (\mu \varphi_1)'''(x)\cos(k\pi x) dx =0,
\end{equation}
it follows from \eqref{expressionmuphiphikresteintegral} that \eqref{esthyp1} implies that
\begin{equation}\label{deriveesdifferentes}
\mu'(1)\not = \mu'(0) \text{ and } \mu'(1)\not = -\mu'(0).
\end{equation}
\end{remark}
\noindent
As proved in~\cite{BM_Tmin}, Hypothesis \ref{hyp} is not necessary to get local exact controllability of \eqref{nl}. However (see~\cite{BeauchardLaurent})  it is a necessary and sufficient condition to get exact controllability of the following linearized equation around the ground state
\begin{equation}\label{lin}
\begin{cases}
i\partial_t \Psi = -\Delta \Psi - v(t)\mu(x)\Phi_1(t,x), &\quad (t,x)\in (0,T)\times (0,1),
\\
\Psi(t,0)=\Psi(t,1)=0, &\quad t\in (0,T),
\\
\Psi(0,x)=\Psi_0(x), &\quad x\in (0,1).
\end{cases}
\end{equation}
\begin{thm}[ \cite{BeauchardLaurent} ]\label{KC}
Let $T>0$ and assume that $\mu$ satisfies Hypothesis \ref{hyp}. Then, for every
\begin{equation} \label{DefEspaceLinearise}
\Psi_0 \in \left\{ \phi \in H^3_{(0)} \: ; \: \Re \lag \phi, \varphi_1 \rag = 0 \right\} =: \HH_0,
\qquad
\Psi_T \in \left\{ \phi \in H^3_{(0)} \: ; \: \Re \lag \phi, \Phi_1(T) \rag = 0 \right\} = : \HH_T,
\end{equation}
there exists $v \in L^2((0,T);\RR)$ such that the solution $\Psi$ of \eqref{lin} with the initial condition $\Psi(0,.)=\Psi_0$ satisfies,
$\Psi(T,.)=\Psi_T$.
\end{thm}
Condition~(\ref{DefEspaceLinearise}) means that $\Psi_0$ and $\Psi_T$ lie in the tangent vector space of $\	S$ in $\varphi_1$ and $\Phi_1(T)$, respectively. Due to the linearization of the preservation of the norm for the bilinear problem, the solution of~(\ref{lin}) satisfies
$\Re \lag \Psi(t), \Phi_1(t) \rag =0$ for every $t\geq0$.

\noindent
The main result of this paper is the construction of feedback laws leading to
rapid stabilization of the linear control system (\ref{lin}).
\begin{thm}\label{main_linearise}
Let $T>0$. Assume that $\mu$ satisfies Hypothesis~\ref{hyp}.
Then, for every $\lambda>0$, there exists $C>0$ and a feedback law $v(t)=K(\Psi(t,\cdot))$ such that
for every $\Psi_0 \in \HH_0$ the associated solution of~(\ref{lin}) satisfies
\begin{equation*}
\|\Psi(t,\cdot)\|_{H^3_{(0)}}\leq C e^{-\lambda t}\|\Psi_0\|_{H^3_{(0)}}.
\end{equation*}
\end{thm}

\medskip
For the sake of simplicity we will focus, for the rest of this article, on the rapid stabilization of the following linear Schr\"odinger equation
\begin{equation}\label{SystLin}
\begin{cases}
i\partial_t \Psi = -\Delta \Psi - v(t)\mu(x)\varphi_1(x), &\quad (t,x)\in (0,T)\times (0,1),
\\
\Psi(t,0)=\Psi(t,1)=0, &\quad t\in (0,T),
\\
\Psi(0,x)=\Psi_0(x), &\quad x\in (0,1).
\end{cases}
\end{equation}
The only difference between~(\ref{lin}) and~(\ref{SystLin}) is that the control term ``$v(t) \mu(x) \Phi_1(t,x)$'' has been replaced by ``$v(t) \mu(x) \varphi_1$''. Using again \cite[Proposition 4]{BeauchardLaurent}, we get the analogous of Theorem~\ref{KC} that is exact controllability with $L^2((0,T);\RR)$ controls of system~(\ref{SystLin}) but now in the state space $H^3_{(0)}((0,1);\CC)$.

We prove the following rapid stabilization result.
\begin{thm}\label{main}
Let $T>0$. Assume that $\mu$ satisfies Hypothesis~\ref{hyp}.
Then, for every $\lambda>0$, there exists $C>0$ and a real-valued feedback law $v(t)=K(\psi(t,\cdot))$ such that,
for every $\Psi_0 \in H^3_{(0)}((0,1);\CC)$, the associated solution of~(\ref{SystLin}) satisfies
\begin{equation*}
\|\Psi(t,\cdot)\|_{H^3_{(0)}}\leq C e^{-\lambda t}\|\Psi_0\|_{H^3_{(0)}}.
\end{equation*}
\end{thm}

\begin{remark}\label{Remarque_lambda1}
The final goal would be to achieve local rapid stabilization of the bilinear problem~(\ref{nl}) toward the ground state $\Phi_1$. To avoid dealing with a moving target, notice that
\begin{equation*}
\|\psi(t,\cdot)-\Phi_1(t,\cdot)\|_{H^3_{(0)}} = \| e^{i \lambda_1 t} \psi(t,\cdot)- \varphi_1 \|_{H^3_{(0)}}.
\end{equation*}
Thus it is simpler to look at the system satisfied by $e^{i \lambda_1 t} \psi(t,\cdot)$.
In the same spirit, we will develop the proof of Theorem~\ref{main} in this article and detail in Appendix~\ref{Annexe_lambda1} how we can modify the proof to obtain Theorem~\ref{main_linearise}. The obtained feedback law does not allow us, for now, to obtain rapid stabilization of~(\ref{nl}).
\end{remark}

\subsection{A finite dimensional example}

Let us explain the general idea of the proof of Theorem~\ref{main} in a finite dimensional setting. Let $A\in \RR^{n\times n}$ and $B\in \RR^{n}$. Consider the control system
\begin{equation}\label{finitecont}
x'(t)=Ax(t)+Bu(t)
\end{equation}
where, at time $t$,  the state is $x(t)\in \RR^n$ and the control is $u(t)\in \RR$. We assume that the system is controllable, which is equivalent to the Kalman rank condition
\begin{equation} \label{Kalman}
rank( B, AB, \dots, A^{n-1}B ) = n
\end{equation}
(see e.g.~\cite[Theorem 1.16]{CoronBook}). It is well-known that the controllability allows one to use the pole-shifting theorem~\cite[Theorem 10.1]{CoronBook}) to design a feedback law $u(t)=Kx(t)$ to obtain the exponential stability with an arbitrary exponential decay rate of \eqref{finitecont}. Let us present a different approach, more suitable for PDEs, to this result. Let $\lambda\in \RR$ and denote the identity matrix of size $n$ by $I$. Consider the target system
\begin{equation}\label{finitestable}
y'(t)=(A-\lambda I)y(t)+Bv(t)
\end{equation}
where, at time $t$, $y(t)\in \RR^n$ is the state and $v(t)\in \RR$ is the control. A straightforward computation shows that, for $v\equiv 0$, the solutions to \eqref{finitestable} satisfy
\[
\|y(t)\|\leq  e^{-(\lambda -\|A\|)t}\|y(0)\|.
\]
\noindent

Let us assume, for the moment, that we can design a transformation $(T,K) \in \RR^{n\times n} \times \RR^{1 \times n}$ such that  if $x(t)$ is the solution of~\eqref{finitecont} with
\begin{equation}
\label{def-trans-finite}
u(t) := Kx(t) + v(t),
\end{equation}
then $y(t) := Tx(t)$ is the solution of \eqref{finitestable}.
Notice that if moreover $T$ is invertible, then
\begin{align*}
\|x(t)\|=\|T^{-1}y(t)\|
&\leq \|T^{-1}\|e^{-(\lambda-\|A\|) t}\|y(0)\|
\\
&\leq \|T^{-1}\|e^{-(\lambda-\|A\|) t}\|Tx(0)\|
\leq \|T^{-1}\|\|T\|e^{-(\lambda-\|A\|) t}\|x(0)\|.
\end{align*}
Therefore, the exponential stability of \eqref{finitecont} with an arbitrary exponential decay rate is reduced to find such a transformation $(T,K)$ with $T$ invertible.
\noindent
The transformation $(T,K)$ maps \eqref{finitecont} into
\begin{equation*}
y'(t) = Tx'(t)= T \left( Ax(t) + B Kx(t) + B v(t) \right) = \left( TA +TBK \right)x(t) + TB v(t),
\end{equation*}
Hence this transformation maps \eqref{finitecont} into \eqref{finitestable}
if and only if
\begin{align}
TA+BK&=AT-\lambda T, \label{l1} \\
TB&=B.\label{l2}
\end{align}
One has the following theorem.
\begin{thm}[\cite{Coron_ICIAM15}]\label{finitedim}
There exists one and only one $(T,K)\in GL_n(\RR)\times \RR^{1\times n}$ satisfying (\ref{l1})-(\ref{l2}).
\end{thm}

The proof of Theorem \ref{finitedim} provided in \cite{Coron_ICIAM15} relies on the
phase variable canonical form (also called controller form) of \eqref{finitecont}.  We present here a different proof (in the case where the eigenvalues of $A$ are simple) more suitable to deal with the infinite dimensional setting, with the additional assumption
\begin{equation}\label{assumpinvfinie}
\lambda>0 \text{ is such that } \left((\lambda_i+\lambda)I-A\right) \text{ is invertible } \, \forall \, 1 \leq i \leq n.
\end{equation}

\beginpf
We first prove that the result holds for $(T,K)\in GL_n(\CC)\times \CC^{1\times n}$. The fact that $(T,K)$ are real-valued follows from the uniqueness of the transformations and the fact that $A$ and $B$ are real-valued.

Denote by $\{\lambda_i,e_i\}_{1\leq i \leq n}$ the eigenvalues and eigenvectors of $A$. Then, \eqref{l1}-\eqref{l2} become
\begin{align}
\left((\lambda_i+\lambda)I-A\right) Te_i&=-BKe_i,  \label{lin1} \\
TBe_i&=Be_i.\label{lin2}
\end{align}
The proof is then divided in four steps. \medskip

\textit{Step 1: Existence of a basis of the state space.} \medskip

Assumption \eqref{assumpinvfinie} implies that that there exists $n$ vectors $f_i, 1\leq i \leq n$ satisfying
\begin{equation}\label{lin1tilde}
\left((\lambda_i+\lambda)I-A\right)f_i=-B, \quad 1\leq i \leq n.
\end{equation}
We begin by proving that the set $\{ f_i\}$ forms a basis of $\CC^n$. Notice that if $K$ is known, then one recovers $Te_i$ from the relation $Te_i=f_i Ke_i$.

Suppose there exists $\{a_i\}_{ 1\leq i \leq n} \subset \CC, \, \{a_i\}_{ 1\leq i \leq n} \neq 0$ such that
\begin{equation}\label{depfin}
\sum_{i=1}^n a_i f_i=0.
\end{equation}
Applying $A$ to this equation and using \eqref{lin1tilde}, we obtain
\[
\sum_{i=1}^n a_i(\lambda_i+\lambda)f_i=-\left(\sum_{i=1}^n a_i \right)B.
\]
Applying successively $A$, we end up with
\begin{equation}\label{relationfini}
\sum_{i=1}^n a_i(\lambda_i+\lambda)^pf_i=-\sum_{k=1}^p \left(\sum_{i=1}^n a_i(\lambda_i+\lambda)^{p-k}\right)A^{k-1}B, \quad \forall p\in \NN^*.
\end{equation}
Note that, for all $j\in \NN^*$, each coefficient
\begin{equation}\label{coefffini}
\sum_{i=1}^n a_i(\lambda_i+\lambda)^{j},
\end{equation}
appears in \eqref{relationfini} for all $p\geq j+1$ in front of $A^{p-j-1}B$. We distinguish two cases. If there exists $j\in \NN^*$ such that there is a coefficient \eqref{coefffini} that is not equal to zero, then it implies that $\{A^{p-j-1}B\}_{p\geq j+1}\subset \Span\{f_i\}$. From the controllability assumption it comes that $\textrm{span}\{A^{p-j-1}B\}_{p\geq j+1}=\CC^n$. Therefore, in this case, the set of $n$ vectors $\{f_i\}$ generates the whole space and consequently a basis of $\CC^n$.

The remaining case is the situation where every coefficient \eqref{coefffini} vanishes i.e.
\begin{equation}\label{derivfini}
\sum_{i=1}^n a_i(\lambda_i+\lambda)^{j} =0, \quad \forall j\in \NN^*.
\end{equation}
In this case, consider the entire function
\[
G:z \in \CC \mapsto \sum_{i=1}^n a_ie^{(\lambda_i+\lambda)z}.
\]
From \eqref{derivfini}, we obtain for all $j\in \NN$,
\[
G^{(j)}(0)=\sum_{i=1}^n a_i(\lambda_i+\lambda)^j =0.
\]
Therefore $G\equiv 0$. Let $\cal{C}:=\textrm{Conv}\{ \, \lambda_i+\lambda \, ; \, 1\leq i \leq n \textrm{ such that } a_i \neq 0\}$ where, for a nonempty subset $R$ of $\CC$,  $\textrm{Conv}R$ is the closed convex-hull of $R$. The set $\cal{C}$ has at least one nonzero extremal, that is a point of $\cal{C}$ such that there exists at least one hyperplane that meets $\cal{C}$ only on this point. One such point must be of the form $\lambda_{k_0}+\lambda$ for $1\leq k_0\leq n$. Therefore, there exists $\theta \in [0,2\pi]$ such that
\begin{equation}\label{termedominantfinie}
\Re(e^{i \theta}(\lambda_{k}+\lambda))<\Re(e^{i \theta}(\lambda_{k_0}+\lambda)), \quad \forall \, 1\leq k \leq n, \, k\neq k_0.
\end{equation}
Let $z=se^{i\theta}$ where $s\in \RR$. We have
\begin{equation}
\label{e-s=0}
e^{-s(\lambda_{k_0}+\lambda)e^{i\theta}}G(se^{i\theta})= a_{k_0}+\sum_{i=1, i\neq k_0}^n a_ie^{s(\lambda_i-\lambda_{k_0})e^{i\theta}} \equiv 0.
\end{equation}
 From \eqref{termedominantfinie} and by letting $s \rightarrow \infty$ in \eqref{e-s=0},
we obtain that $a_{k_0}=0$. It is in contradiction with the fact that the set $\cal{C}$ contains only nonzero $a_i$. Therefore $a_i=0, 1\leq i \leq n$ so the set $\{f_i\}$ is independant and consequently a basis of $\CC^n$. The two cases were covered which implies that $\{f_i\}$ is a basis of $\CC^n$.
\begin{remark}
The latter part of the proof of the existence of a basis could have been done using the Vandermonde matrix. The proof presented here has the advantage that it may be applied in the infinite dimensional setting.
\end{remark}

\textit{Step 2: Existence of the transformation $(T,K)$.} \medskip

The transformation is obtained using \eqref{lin2}. Indeed, let
\[
B=\sum_{i=1}^n b_i e_i.
\]
Notice that, by the controllability assumption, $ b_i \neq 0, 1\leq i \leq n$. Then,
\begin{equation}\label{decompfini}
TB = B \: \Longleftrightarrow \: B = \sum_{i=1}^n b_i Te_i = \sum_{i=1}^n b_i Ke_i f_i.
\end{equation}
Since $\{f_i\}_{1\leq i\leq n}$ is a basis of $\CC^n$, there exists $\{Ke_i\}_{1\leq i\leq n} \subset \CC$ such that the last equation is verified, allowing to define $T\in \CC^{n,n}$ and $K\in \CC^{1,n}$ such that \eqref{lin1} and \eqref{lin2} hold.
\medskip

\textit{Step 3: Uniqueness of $(T,K)$.}
\medskip

To prove the uniqueness of the transformations $(T,K)$, consider $(T_1,K_1)$ and $(T_2,K_2)$ solutions of \eqref{lin1}-\eqref{lin2}. Therefore $(T_1-T_2,K_1-K_2)$ satisfies \eqref{lin1} and
\begin{equation}\label{uniquenesslin}
(T_1-T_2)B=0.
\end{equation}
Since $(T_1-T_2,K_1-K_2)$ satisfies \eqref{lin1}, we use the basis constructed previously and \eqref{uniquenesslin} to prove that $K_1=K_2$ and $T_1=T_2$. With the uniqueness of the transformation and the fact that $A$ and $B$ are real-valued, one ensures that the transformations are real-valued since $(\overline{T},\overline{K})$ is also a solution of \eqref{lin1}-\eqref{lin2}.
\medskip

\textit{Step 4: Invertibility of $T$.}
\medskip

Let $T\in \CC^{n,n}$ and $K\in \CC^{1,n}$ be such that \eqref{lin1} and \eqref{lin2} hold. We prove that $T$ is invertible by showing that $\textrm{Ker }T^*=\{0\}$. Let $x\in \textrm{Ker }T^*$. From~\eqref{l1}-\eqref{l2}, we obtain
\begin{equation*}
T^*A^*x= (A^*T^*+K^*B^*T^*+\lambda T^*)x=0.
\end{equation*}
Therefore $\textrm{Ker }T^*$ is stable by $A^*$. From~\eqref{l2} it comes that
\begin{equation*}
B^* x = B^* T^* x = 0.
\end{equation*}
Thus there exists $\tilde{x}$ eigenvector of $A^*$ in $\textrm{Ker }T^* \subset \textrm{Ker } B^*$. From the controllability assumption and the Hautus test (see for instance~\cite[Prop. 1.5.5]{TucsnakWeissBook}) it comes that $\textrm{Ker }T^*=\{0\}$.

\endpf

If the functional setting in the infinite dimensional case makes the proof more tricky, the strategy we use remains the same. Riesz basis results will be used to prove the existence of a basis of the state space and the invertibility of the transformation will be proved using the approximate controllability of the studied system.

The main technical difficulty of this paper lies in the decomposition of $B$ \eqref{decompfini} in the basis of the state space. Indeed, the control operator $B$ is admissible but not bounded in the state space. A careful analysis of the Fourier components of the control operator $B$ allow us to define a transformation $T$ which is bounded from the state space into itself but the feedback transformation won't be bounded from the state space into $\RR$. Even so, the transformation $T$ will be proved to be invertible and the closed-loop linear equation will be proved to be well-posed in the state space. It is important to note that this technical difficulty is in fact essential for the invertibility of $T$ (see Remark~\ref{Rk1_compacite}). Indeed, in our case, if $B$ were to be bounded, then the transformation $T$ would be compact and thus not continuously invertible. However, the unboundedness of $K$ from the state space into $\RR$ prevents us to prove directly the well-posedness of the closed-loop nonlinear equation.

Let us underline that the uniqueness condition $TB=B$, which was used implicitly in similar previous works, will be crucial not only to obtain the existence and uniqueness of the transformation, but also to define rigorously \eqref{l1} in the infinite dimensional setting.
\subsection{The linear Schr\"odinger equation}

As presented in the previous paragraph, the strategy to prove the rapid stabilization of the linear equation \eqref{SystLin} is inspired by the backstepping method. Recast the equations for $\Psi^1+i\Psi^2=\Psi$ as
\begin{equation}\label{slin}
\begin{cases}
\partial_t \begin{pmatrix} \Psi^1 \\ \Psi^2 \end{pmatrix} = \begin{pmatrix} 0 & -\Delta  \\ \Delta & 0 \end{pmatrix} \begin{pmatrix} \Psi^1 \\ \Psi^2 \end{pmatrix}
+ v(t) \begin{pmatrix} 0 \\ (\mu \varphi_1)(x) \end{pmatrix},  &
(t,x) \in (0,T) \times (0,1),
\\
\Psi^1(t,0) = \Psi^1(t,1) = 0, \quad \Psi^2(t,0) = \Psi^2(t,1) = 0,
& t \in  (0,T),
\\
\Psi^1(0,x) = \Psi^1_0(x), \quad \Psi^2(0,x) = \Psi^2_0(x),
& x \in (0,1).
\end{cases}
\end{equation}
where $\Psi^1_0$ and $\Psi^2_0$ are the real and imaginary part of $\Psi^0$ respectively.
From now on, all the functional spaces are real-valued, except when specified. Moreover to deal with those real and imaginary parts, we denote, for simplicity,
\begin{equation} \label{DefX^s_0}
X^s_{(0)}((0,1);\RR) := \left( H^s_{(0)}((0,1);\CC) \cap L^2((0,1);\RR) \right)^2,
\end{equation}
with the product topology.
We will use the following operators
\begin{align*}
A: D(A) & \longrightarrow X^3_{(0)}   &  B : \RR  & \longrightarrow \left(C^{\infty}  \times (H^3\cap H^1_0)((0,1); \RR)\right)\tr \\
 \begin{pmatrix} \Psi^1 \\ \Psi^2 \end{pmatrix}\, \, &\longmapsto \begin{pmatrix}  -\Delta \Psi^2 \\ \Delta \Psi^1 \end{pmatrix}\, , &
 a \, \, &\longmapsto a \begin{pmatrix} 0 \\ (\mu \varphi_1)(x) \end{pmatrix}, \\
\end{align*}
with $D(A):=X^5_{(0)}$. Based on the previous work of the first author and Q.~L\"u (\cite{CoronLu14, CoronLu15}), instead of  Volterra transformations of the second kind usually used for the backstepping approach, we seek for transformations $(T,K)$ of the form
\begin{align}
T: X^3_{(0)}& \longrightarrow X^3_{(0)} \nonumber \\
 \begin{pmatrix} \Psi^1 \\ \Psi^2 \end{pmatrix}\, \, &\longmapsto \int_0^1 \begin{pmatrix} k_{11}(x,y) & k_{12}(x,y) \\ k_{21}(x,y) & k_{22}(x,y) \end{pmatrix}
\begin{pmatrix} \Psi^1(y) \\ \Psi^2(y) \end{pmatrix} dy, \label{transformationR} \\
K: D(K) \subset X^3_{(0)}& \longrightarrow \RR \nonumber \\
 \begin{pmatrix} \Psi^1 \\ \Psi^2 \end{pmatrix} &\longmapsto
\int_0^1 \alpha^1(y) \Psi^1(y) + \alpha^2(y) \Psi^2(y) dy,\label{feedbackR}
\end{align}
such that if $(\Psi^1,\Psi^2){\tr}$ is solution of \eqref{slin} with
\begin{equation}
\label{v=Kpsi}
v(t) :=
K\begin{pmatrix} \Psi^1 (t,.) \\ \Psi^2(t,.) \end{pmatrix},
\end{equation}
then $(\xi^1(t,.),\xi^2(t,.)){\tr}:=T(\Psi^1(t,.),\Psi^2(t,.)){\tr}$
is solution to
\begin{equation}\label{slin_stable}
\begin{cases}
\partial_t \begin{pmatrix} \xi^1 \\ \xi^2 \end{pmatrix}
= \begin{pmatrix} 0 & -\Delta  \\ \Delta & 0 \end{pmatrix}
\begin{pmatrix} \xi^1 \\ \xi^2 \end{pmatrix}
- \lambda \begin{pmatrix} \xi^1 \\ \xi^2 \end{pmatrix},  &
(t,x) \in (0,T) \times (0,1),
\\
\xi^1(t,0) = \xi^1(t,1) = 0, \quad \xi^2(t,0) = \xi^2(t,1) = 0,&
t \in \times (0,T),
\\
\xi^1(0,x) = \xi^1_0(x), \quad \xi^2(0,x) = \xi^2_0(x), &
x \in \times (0,1),
\end{cases}
\end{equation}
with $(\xi^1_0 , \xi^2_0)\tr
= T (\Psi^1_0 ,\Psi^2_0)\tr$ and $T$ is invertible in the state space.
The decomposition in real and imaginary part of the solution of \eqref{SystLin} is made in order to ensure that the feedback $v(t)=K(\Psi^1(t,.),\Psi^2(t,.)){\tr}$ is real-valued. Note that
\begin{equation}
\label{decayxi}
\left\|\begin{pmatrix} \xi^1(t) \\ \xi^2(t) \end{pmatrix} \right\|_{{X^3_{(0)}}}
\leq e^{-\lambda t} \left\|\begin{pmatrix} \xi^1_0 \\ \xi^2_0 \end{pmatrix} \right\|_{{X^3_{(0)}}}, \quad t\in [0,+\infty).
\end{equation}
The kernels are defined through the equations they must satisfy for $(T,K)$ to map solutions of \eqref{slin} to solutions of \eqref{slin_stable}. This is done formally in Section~\ref{sec_heuristique} together with a more detailed presentation of this strategy.

\subsection{A brief review of previous results}

The controllability properties for the Schr\"odinger equation were mostly studied in the usual (in opposition to the bilinear model presented here) linear setting.
For the control of the linear Schr\"odinger equation with internal control (localized on a subdomain), we refer to the survey~\cite{Laurent13_survey} and the references therein. In this more classical setting we also mention~\cite{MachtyngierZuazua94, LasieckaTriggianiZhang_Schrodinger04, DehmanGerardLebeau06} concerning stabilization.

\paragraph{Exact controllability of the bilinear Schr\"odinger equation.}
~~\\
The first local controllability results on the bilinear Schr\"odinger equation appear in \cite{Beauchard05, BeauchardCoron06, BeauchardLaurent}.
These local controllability results have been extended with weaker assumptions in \cite{BM_Tmin}, in a more general setting in infinite time~\cite{NersesyanNersisyan1D} and also in the case of simultaneous controllability of a finite number of particles in \cite{Morancey_simultane, MoranceyNersesyan13}. Note that, despite the infinite speed of propagation, it was proved in \cite{2006-Coron-SCL, 2014-Beauchard-Coron-Teismann-SCL, BM_Tmin, Morancey_simultane} that a minimal amount of time is required for the controllability of some bilinear Schr\"odinger equations. More recently, local exact controllability has been established in~\cite{BeauchardLaurent15} for a Schr\"odinger-Poisson model in 2D (see also~\cite{MehatsPrivatSigalotti14} for approximate controllability) and for the analogue of~\eqref{nl} with a control depending on time and space in dimension less or equal than $3$~\cite{Puel15}.

\paragraph{Approximate controllability and stabilization of the bilinear Schr\"odinger equation.}
~~\\
The above mentioned results of exact controllability are mostly limited to the one dimensional case. In a more general setting, the available results deal with approximate controllability. Using geometric control techniques on appropriate Galerkin approximations, approximate controllability in different settings has been proved~\cite{BCC_Regular14, BoscainCaponigroSigalotti13, BCCS11, CMSB09}. For a detailed presentation, see the survey~\cite{BoscainChambrionSigalotti_review}.

However, most of these results are not suitable to prove approximate controllability in higher norms (typically $H^3_{(0)}$) and thus approximate controllability for bilinear Schr\"odinger equations has also been studied from the Lyapunov functional point of view~\cite{Mirrahimi09, BeauchardMirrahimi09, Nersesyan, Nersesyan10}. Though it enabled global controllability results, this strategy usually gives no indication on the convergence rate.

\paragraph{Rapid stabilization.}
~~\\
The strategy used in this article is inspired from backstepping techniques. Initially developed to design, in a recursive way, more effective feedback laws for globally asymptotic stable finite dimensional system for which a feedback law and a Lyapunov function are already known (we refer to~\cite{KKK_Book, Sontag_Book, CoronBook} for a comprehensive introduction in finite dimension and~\cite{CAN, KrsticLiu00} for an application of this method to partial differential equations), the backstepping approach was later used in the infinite dimensional setting to design feedback laws by mapping the system to stabilize to a target stable system. To our knowledge, this strategy was first introduced in the context of partial differential equations to design a feedback law for heat equation \cite{BK1} and, later on, for a class of parabolic PDE \cite{BK2}. The backstepping-like change of coordinates of the semi-discretized equation maps the discrete solution to stabilize to a stable solution. The corresponding continuous mapping obtained is a Volterra transformation of the second kind, seen as an infinite dimensional backstepping transformation from the triangular domain of definition of the transformation.
This backstepping strategy in infinite dimensions led to a series of work (see~\cite{KrsticSmyshlyaev_Book} for a global presentation and \cite{KrsticGuoSmyshlyaev11} for the use of the backstepping approach for the rapid stabilization  of a Schr\"odinger equation with a boundary control). The backstepping approach can be used to get stabilization in finite time as shown in \cite{2015-Coron-Nguyen-preprint, CoronVazquezKrsticBastin13}. Moreover
it is not limited to linear equations, as shown in~\cite{CerpaCoron13, CoronVazquezKrsticBastin13}. Though using Volterra transformations of the second kind provides easily invertible transformations it is also limited. Fredholm transformations, from which this paper is inspired, has already been used~\cite{CoronLu14, CoronLu15, CoronHuOlive16} for rapid stabilization using boundary feedback laws.

Abstract methods have been developed to obtain the rapid stabilization of linear partial differential equations. Among them, we cite the works \cite{Komornik}, \cite{Vest} and \cite{Urquiza}, based on the Gramian approach and the Riccati equations, which could be applied to obtain the rapid stabilization of the linearized equation \eqref{slin} as $(0, \mu \varphi_1) \in D(\mathcal{A}^*)'$. However, it seems difficult to obtain, for various physical systems, the local rapid stabilization of the nonlinear equation using those methods. For example, at least for the moment, one does not know how to deduce from
\cite{CerpaCrepeau09}, where the rapid stabilization of a linearized Korteweg-de Vries equation is obtained by using the method developed in \cite{Urquiza}, the rapid stabilization of the associated nonlinear Korteweg-de Vries equation. This is in contrast with the method used here (linear transformations to suitable target systems) applied to the same Korteweg-de Vries equation. Indeed, as shown in \cite{CoronLu14}, the feedback laws obtained by means of this method allows to get the rapid stabilization for the nonlinear Korteweg-de Vries equation. One may therefore hope that, as in \cite{CoronLu14}, our feedback law $K$ being quite explicit it might allow to obtain the rapid stabilization of the nonlinear equation. Note however that in \cite{CoronLu14} the feedback
law $K$ is continuous, which is not the case in our situation. It makes the application to the nonlinear system more complicated to study and requires suitable nonlinear modifications of the linear feedback law $K$.
\subsection{Structure of the article}

In Sec.~\ref{sec_heuristique}, we give a detailed presentation of the strategy used to construct the transformation $(T,K)$ and give a formal expression of this transformation.
In Sec.~\ref{sec_transformation} we prove that this formal transformation $T$ is well defined and is continuous in the state space $X^3_{(0)}$.
Then, we prove in Sec.~\ref{sec_inversibilite} that the previous transformation is indeed invertible in the state space. These properties of $T$ will follow using Hypothesis~\ref{hyp} i.e. exact controllability of the linearized system.
We end the proof of Theorem~\ref{main} in Sec.~\ref{sec_WellPosedness} by proving that the constructed feedback $K$ leads to a well-posed closed-loop system (i.e. the equation (\ref{slin}) with $v$ defined by \eqref{v=Kpsi}) and that $T$ actually maps the closed-loop system to the exponentially stable solutions of~(\ref{slin_stable}).
In Appendix~\ref{Annexe_SaintVenant} we study in a similar way a simplified Saint-Venant equation which exhibits the same phenomenon but on which we explicitly compute the transformation $(T,K)$.

\section{Heuristic construction of the transformations}
\label{sec_heuristique}

We recall that we look for a transformation $(T,K)$ of the form~(\ref{transformationR})-~(\ref{feedbackR}).
Let us derive the set of equations for $(T,K)$ to map solutions of \eqref{slin} to solutions of \eqref{slin_stable}. First, to ensure that the boundary conditions of~(\ref{slin_stable}) are satisfied, we assume that
\begin{equation} \label{conditionBord_noyaux}
k_{ij}(0,y) = k_{ij} (1,y) = 0, \quad \forall y \in (0,1), \quad \forall i, j \in \{1,2\}.
\end{equation}

\noindent
Using the fact that $(\Psi^1,\Psi^2)\tr \in X^3_{(0)}$, and imposing the conditions
\begin{equation} \label{conditionBordIPP_noyaux}
k_{ij}(x,0) = k_{ij} (x,1) = 0, \quad \forall x \in (0,1), \quad \forall i, j \in \{1,2\},
\end{equation}
formal computations, denoting $\Delta_x$ and $\Delta_y$ the Dirichlet Laplacian with respect to $x$ and $y$ variables respectively, lead to
\begin{align*}
&\partial_t \xi^1(t,x) + \Delta \xi^2(t,x) + \lambda \xi^1(t,x)
\\
&= \int_0^1 k_{11}(x,y) \left( -\Delta_y \Psi^2(t,y) \right) + k_{12}(x,y) \left( \Delta_y \Psi^1(t,y) + v(t) (\mu \varphi_1)(y) \right) dy
\\
&+ \int_0^1 \Delta_x k_{21}(x,y) \Psi^1(t,y) + \Delta_{x}  k_{22}(x,y) \Psi^2(t,y) dy
\\
&+ \int_0^1 \lambda k_{11}(x,y) \Psi^1(t,y) + \lambda k_{12}(x,y) \Psi^2(t,y) dy
\\
&= \int_0^1 \left(\Delta_x k_{21} + \Delta_y k_{12} + \lambda k_{11}  \right)(x,y) \Psi^1(t,y) dy
\\
&+ \int_0^1 \left(\Delta_x k_{22} - \Delta_y k_{11} + \lambda k_{12} \right)(x,y) \Psi^2(t,y) dy
\\
&+ v(t) \int_0^1 \left( \mu \varphi_1 \right)(z) k_{12}(x,z) dz.
\end{align*}
The boundary conditions~(\ref{conditionBordIPP_noyaux}) were imposed to avoid boundary terms in the integrations by parts. Using the expression \eqref{feedbackR} of the feedback leads to
\begin{align*}
&\partial_t \xi^1(t,x) + \Delta \xi^2(t,x) + \lambda \xi^1(t,x)
\\
&= \int_0^1 \left[ \left(\Delta_y k_{12} + \Delta_x k_{21} + \lambda k_{11}\right)(x,y) + \alpha^1(y) \left( \int_0^1 \left( \mu \varphi_1 \right)(z) k_{12}(x,z) dz \right)  \right] \Psi^1(t,y) dy
\\
&- \int_0^1 \left[ \left(\Delta_y k_{11} - \Delta_x k_{22} - \lambda k_{12}\right)(x,y) - \alpha^2(y) \left( \int_0^1 \left( \mu \varphi_1 \right)(z) k_{12}(x,z) dz \right)  \right] \Psi^2(t,y) dy.
\tag{\theequation} \addtocounter{equation}{1}
\label{equationnoyau1_sansTB=B}
\end{align*}

\noindent
In the same way one gets
\begin{align*}
&\partial_t \xi^2(t,x) - \Delta \xi^1(t,x) + \lambda \xi^2(t,x)
\\
&= \int_0^1 \left[ \left(\Delta_y k_{22} - \Delta_x k_{11} + \lambda k_{21}\right)(x,y) + \alpha^1(y) \left( \int_0^1 \left( \mu \varphi_1 \right)(z) k_{22}(x,z) dz \right)  \right] \Psi^1(t,y) dy
\\
&- \int_0^1 \left[ \left(\Delta_y k_{21} + \Delta_x k_{12} - \lambda k_{22}\right)(x,y) - \alpha^2(y) \left( \int_0^1 \left( \mu \varphi_1 \right)(z) k_{22}(x,z) dz \right)  \right] \Psi^2(t,y) dy.
\tag{\theequation} \addtocounter{equation}{1}
\label{equationnoyau2_sansTB=B}
\end{align*}

\noindent
If we want $(\xi^1, \xi^2)$ to satisfy~(\ref{slin_stable}) then we need to find the functions $k_{ij}$ and $\alpha^j$ satisfying
\begin{equation} \label{SystemeNoyau_sansTB=B}
\left\{
\begin{aligned}
&\left(\Delta_y k_{11} - \Delta_x k_{22} - \lambda k_{12}\right)(x,y) - \alpha^2(y) \left( \int_0^1 \left( \mu \varphi_1 \right)(z) k_{12}(x,z) dz \right) = 0, \; (x,y)\in (0,1)^2,
\\
&\left(\Delta_y k_{12} + \Delta_x k_{21} + \lambda k_{11}\right)(x,y) + \alpha^1(y) \left( \int_0^1 \left( \mu \varphi_1 \right)(z) k_{12}(x,z) dz \right) = 0,\; (x,y)\in (0,1)^2,
\\
&\left(\Delta_y k_{21} + \Delta_x k_{12} - \lambda k_{22}\right)(x,y) - \alpha^2(y) \left( \int_0^1 \left( \mu \varphi_1 \right)(z) k_{22}(x,z) dz \right) = 0, \; (x,y)\in (0,1)^2,
\\
&\left(\Delta_y k_{22} - \Delta_x k_{11} + \lambda k_{21}\right)(x,y) + \alpha^1(y) \left( \int_0^1 \left( \mu \varphi_1 \right)(z) k_{22}(x,z) dz \right) = 0, \; (x,y)\in (0,1)^2,
\\
&k_{ij}(x,0) = k_{ij}(x,1) = 0, \quad x\in (0,1),
\\
&k_{ij}(0,y) = k_{ij}(1,y) = 0, \quad y\in (0,1).
\end{aligned}
\right.
\end{equation}

\paragraph{A fundamental extra condition.}
One could try to solve rightaway \eqref{SystemeNoyau_sansTB=B} and prove the invertibility of the transformation $T$ but the non-local terms yield a tedious task. To overcome this difficulty, one introduces, as in the finite dimensional setting, what will be referred throughout this article as the $TB=B$ condition,
\begin{equation*}
T \begin{pmatrix} 0 \\ (\mu \varphi_1)(x) \end{pmatrix}
= \int_0^1 \begin{pmatrix} k_{12}(x,y) (\mu \varphi_1)(y) \\ k_{22}(x,y) (\mu \varphi_1)(y) \end{pmatrix} dy=\begin{pmatrix} 0 \\ (\mu \varphi_1)(x) \end{pmatrix}.
\end{equation*}
Plugging this into~(\ref{SystemeNoyau_sansTB=B}) we obtain that we now seek for a solution to
\begin{equation} \label{SystemeNoyau}
\begin{cases}
\left(\Delta_y k_{11} - \Delta_x k_{22} - \lambda k_{12}\right)(x,y) = 0, & (x,y)\in (0,1)^2,
\\
\left(\Delta_y k_{12} + \Delta_x k_{21} + \lambda k_{11}\right)(x,y) = 0, & (x,y)\in (0,1)^2,
\\
\left(\Delta_y k_{21} + \Delta_x k_{12} - \lambda k_{22}\right)(x,y) - \alpha^2(y) (\mu \varphi_1)(x) = 0, & (x,y)\in (0,1)^2,
\\
\left(\Delta_y k_{22} - \Delta_x k_{11} + \lambda k_{21}\right)(x,y) + \alpha^1(y) (\mu \varphi_1)(x) = 0, & (x,y)\in (0,1)^2,
\\
k_{ij}(x,0) = k_{ij}(x,1) = 0, & x\in (0,1),
\\
k_{ij}(0,y) = k_{ij}(1,y) = 0, & y\in (0,1),
\end{cases}
\end{equation}
together to  the $TB=B$ condition
\begin{equation} \label{ConditionTB=B}
\begin{cases}
\displaystyle \int_0^1 k_{12}(x,y) (\mu \varphi_1)(y) dy =0, & x\in (0,1),
\\
\displaystyle  \int_0^1 k_{22}(x,y) (\mu \varphi_1)(y) dy =(\mu \varphi_1)(x), & x\in (0,1).
\end{cases}
\end{equation}

\begin{remark}

In \cite{CoronLu14, CoronLu15}, the authors were dealing with a boundary control for the Korteweg-de Vries equation and for the Kuramoto-Sivashinsky equation. In their case, the $TB=B$ condition writes
\[k_y(x,L)=0, \quad x\in (0,L),
\]
for the former and
\[k_{yy}(x,L)=0, \quad x\in (0,L),
\]
for the latter, where $k$ is the kernel of the Fredholm transformation in each case. Contrary to our framework, this boundary condition appeared naturally from the integration by parts performed in order to obtain the equation on the kernel. It was not seen as a particular boundary condition although, a careful analysis of their work shows that the $TB=B$ condition was used to prove the uniqueness of the transformation $(T,K)$. The common ground between their work and this article is the additional regularity that the kernel needs in order to satisfy the $TB=B$ condition. Notice that in what we present in this article, the relation between the kernels $k_{ij}$ and $\alpha^j$ is more intricate and considerably modify the analysis.
\end{remark}

\paragraph{Formal decomposition}
The global strategy to construct a solution of~(\ref{SystemeNoyau})-(\ref{ConditionTB=B}) is the following. First assume that $\alpha^1$ and $\alpha^2$ are known. This enables us to compute the kernels $k_{ij}$ satisfying~(\ref{SystemeNoyau}) as functions of $\alpha^1$ and $\alpha^2$. Then we prove that we find $\alpha^1$ and $\alpha^2$ such that~(\ref{ConditionTB=B}) is satisfied.

\noindent
We decompose the functions in the following form
\begin{equation} \label{decompositionNoyaux}
k_{ij}(x,y) = \sum_{n=1}^{+\infty} f_n^{ij}(x) \varphi_n(y), \qquad
\alpha^j(y) = \sum_{n=1}^{+\infty} \alpha_n^j \varphi_n(y).
\end{equation}
This leads to
\begin{equation} \label{Systeme_fn}
\left\{
\begin{aligned}
&f_n^{11}{''}(x) + \lambda_n f_n^{22}(x) - \lambda f_n^{21}(x) = \alpha^1_n (\mu \varphi_1)(x),
\\
&f_n^{12}{''}(x) - \lambda_n f_n^{21}(x) - \lambda f_n^{22}(x) = \alpha^2_n (\mu \varphi_1)(x),
\\
&f_n^{21}{''}(x) - \lambda_n f_n^{12}(x) + \lambda f_n^{11}(x) = 0,
\\
&f_n^{22}{''}(x) + \lambda_n f_n^{11}(x) + \lambda f_n^{12}(x) = 0,
\\
&f_n^{ij}(0) = f_n^{ij}(1) = 0.
\end{aligned}
\right.
\end{equation}
The $TB=B$ condition (\ref{ConditionTB=B}) becomes
\begin{equation} \label{ContionTB=B_fn}
\left\{
\begin{aligned}
&\sum_{n=1}^{+\infty} \< \mu \varphi_1, \varphi_n \> f_n^{12}(x) = 0,
\\
&\sum_{n=1}^{+\infty} \< \mu \varphi_1, \varphi_n \> f_n^{22}(x) = (\mu \varphi_1)(x).
\end{aligned}
\right.
\end{equation}

\noindent
As mentioned, we begin by assuming that the feedback law is known.
We consider two sequences $(\beta^1_n)_{n \in \NN^*}$ and $(\beta^2_n)_{n \in \NN^*}$ to be precised later on.

Let $(g_n^{11}, g_n^{12}, g_n^{21}, g_n^{22}){\tr}$ be the solution of system~(\ref{Systeme_fn}) with right-hand side $\left(\beta^1_n (\mu \varphi_1), 0, 0, 0\right){\tr}$ and let $(h_n^{11}, h_n^{12}, h_n^{21}, h_n^{22}){\tr}$ be the solution of system~(\ref{Systeme_fn}) with right-hand side $\left(0, \beta^2_n (\mu \varphi_1), 0, 0\right){\tr}$.
System~(\ref{Systeme_fn}) being linear, it comes that
\begin{equation} \label{Decompositionfn_gn_hn}
f_n^{ij} = \frac{\alpha^1_n}{\beta^1_n} g_n^{ij} + \frac{\alpha^2_n}{\beta^2_n} h_n^{ij}.
\end{equation}

\noindent
Decomposing $g_n^{ij}$ in the $L^2$-orthonormal basis $(\varphi_k)_{k\in \NN^*}$, if we denote by $A_{nk}$ the following matrix
\begin{equation} \label{Def_Ank}
A_{nk} =
\begin{pmatrix}
-\lambda_k & 0 & -\lambda & \lambda_n  \\
0 & -\lambda_k & -\lambda_n & -\lambda \\
\lambda & -\lambda_n & -\lambda_k & 0  \\
\lambda_n & \lambda & 0 & -\lambda_k
\end{pmatrix},
\end{equation}
system~(\ref{Systeme_fn}) leads to
\begin{equation*}
A_{nk} \begin{pmatrix} \<g_n^{11}, \varphi_k\> \\ \<g_n^{12}, \varphi_k\> \\ \<g_n^{21}, \varphi_k\> \\ \<g_n^{22}, \varphi_k\> \end{pmatrix}
= \begin{pmatrix} \beta^1_n \< \mu \varphi_1, \varphi_k \> \\ 0 \\ 0 \\ 0 \end{pmatrix}.
\end{equation*}
Then
\begin{equation} \label{Formule_gn}
\left\{
\begin{aligned}
&g_n^{11} (x) = \sum_{k=1}^{+\infty} \frac{\lambda_k \left( \lambda_n^2 - \lambda^2 - \lambda_k^2 \right)}{\delta_{nk}(\lambda)} \beta^1_n \< \mu \varphi_1, \varphi_k \> \varphi_k(x),
\\
&g_n^{12} (x)= \sum_{k=1}^{+\infty} \frac{2 \lambda \lambda_k \lambda_n}{\delta_{nk}(\lambda)} \beta^1_n \< \mu \varphi_1, \varphi_k \> \varphi_k(x),
\\
&g_n^{21} (x) = \sum_{k=1}^{+\infty} \frac{-\lambda \left( \lambda^2 + \lambda_k^2 + \lambda_n^2\right)}{\delta_{nk}(\lambda)} \beta^1_n \< \mu \varphi_1, \varphi_k \> \varphi_k(x),
\\
&g_n^{22} (x) = \sum_{k=1}^{+\infty} \frac{\lambda_n \left( \lambda^2 - \lambda_k^2 + \lambda_n^2\right)}{\delta_{nk}(\lambda)} \beta^1_n \< \mu \varphi_1, \varphi_k \> \varphi_k(x),
\end{aligned}
\right.
\end{equation}
where
\begin{equation}
\label{defdeltank}
\delta_{nk}(\lambda) = \det(A_{nk}) = \left( \lambda^2 + (\lambda_k - \lambda_n)^2\right) \left( \lambda^2 + (\lambda_k + \lambda_n)^2 \right).
\end{equation}
The same computations can be carried out for $h_n^{ij}$, leading to the following relations
\begin{equation} \label{Relation_gn_hn}
\begin{aligned}
&h_n^{11} = -\frac{\beta^2_n}{\beta^1_n} g_n^{12}, \qquad
h_n^{12} = \frac{\beta^2_n}{\beta^1_n} g_n^{11},
\\
&h_n^{21} = -\frac{\beta^2_n}{\beta^1_n} g_n^{22}, \qquad
h_n^{22} = \frac{\beta^2_n}{\beta^1_n} g_n^{21}.
\end{aligned}
\end{equation}
For the sake of simplicity, we denote by $c_{nk}^{ij}$ and $d_{nk}^{ij}$ the coefficients such that
\begin{equation} \label{DecompositionFourier_gn_hn}
g_n^{ij}(x) = \sum_{k=1}^{+\infty} c_{nk}^{ij} \beta^1_n \< \mu \varphi_1, \varphi_k \> \varphi_k(x),
\quad
h_n^{ij}(x) = \sum_{k=1}^{+\infty} d_{nk}^{ij} \beta^2_n \< \mu \varphi_1, \varphi_k \> \varphi_k(x).
\end{equation}

\paragraph{Summary of the construction.}
Finally, using the definition of the transformation~(\ref{transformationR}), the decompositions~(\ref{decompositionNoyaux}), (\ref{Decompositionfn_gn_hn}) and the relations~(\ref{Relation_gn_hn}) it comes that the transformation $T$ we are looking for is defined by
\begin{equation} \label{Def_T}
\begin{aligned}
T \begin{pmatrix} \Psi^1 \\ \Psi^2 \end{pmatrix}
&= \begin{pmatrix}
\begin{displaystyle}
\sum_{n=1}^{+\infty} \left( \frac{\alpha^1_n}{\beta^1_n} g_n^{11} + \frac{\alpha^2_n}{\beta^2_n} h_n^{11} \right) \<\Psi^1, \varphi_n\>
+
\sum_{n=1}^{+\infty} \left( \frac{\alpha^1_n}{\beta^1_n} g_n^{12} + \frac{\alpha^2_n}{\beta^2_n} h_n^{12} \right) \<\Psi^2, \varphi_n\>
\end{displaystyle}
\\
\begin{displaystyle}
\sum_{n=1}^{+\infty} \left( \frac{\alpha^1_n}{\beta^1_n} g_n^{21} + \frac{\alpha^2_n}{\beta^2_n} h_n^{21} \right) \<\Psi^1, \varphi_n\>
+
\sum_{n=1}^{+\infty} \left( \frac{\alpha^1_n}{\beta^1_n} g_n^{22} + \frac{\alpha^2_n}{\beta^2_n} h_n^{22} \right) \<\Psi^2, \varphi_n\>
\end{displaystyle}
\end{pmatrix}
\\
&=
\sum_{n=1}^{+\infty}  \left[ \frac{\alpha^1_n \<\Psi^2,\varphi_n\>}{\beta^1_n} - \frac{\alpha^2_n
\<\Psi^1,\varphi_n\>}{\beta^1_n} \right]
\begin{pmatrix} g_n^{12} \\ g_n^{22} \end{pmatrix}
+ \left[ \frac{\alpha^1_n \<\Psi^1,\varphi_n\>}{\beta^2_n} + \frac{\alpha^2_n \<\Psi^2,\varphi_n\>}{\beta^2_n} \right]
\begin{pmatrix} h_n^{12} \\ h_n^{22} \end{pmatrix}.
\end{aligned}
\end{equation}
In the same spirit, the $TB=B$ condition ~(\ref{ContionTB=B_fn}) becomes
\begin{equation} \label{Def_TB=B}
\begin{pmatrix} 0 \\ \mu \varphi_1 \end{pmatrix}
= \sum_{n=1}^{+\infty} \frac{\alpha^1_n \<\mu \varphi_1, \varphi_n\>}{\beta^1_n} \begin{pmatrix} g_n^{12} \\ g_n^{22} \end{pmatrix} + \frac{\alpha^2_n \<\mu \varphi_1, \varphi_n\>}{\beta^2_n} \begin{pmatrix} h_n^{12} \\ h_n^{22} \end{pmatrix}.
\end{equation}
This ends the heuristic of the construction of the transformation and the feedback law. Indeed, in the following section we prove that for a suitable choice of the sequences $(\beta^1_n)_{n \in \NN^*}$ and $(\beta^2_n)_{n \in \NN^*}$ then
\begin{equation} \label{Def_B}
\mathcal{B} := \left \{
\begin{pmatrix} g_n^{12} \\ g_n^{22} \end{pmatrix} , \:
\begin{pmatrix} h_n^{12} \\ h_n^{22} \end{pmatrix} \: ; \: n \in \NN^*
\right\},
\end{equation}
is a Riesz basis of $X^2_{(0)}$ (see Proposition~\ref{prop_baseRiesz}). Then from (\ref{Def_TB=B}) we get the feedback laws from the expansion of $(
0, \mu \varphi_1)\tr$
in the basis $\mathcal{B}$.
Finally, we study the behaviour of the coefficients $\alpha_n^1$ and $\alpha_n^2$ as $n$ goes to infinity to prove that the transformation $T$ given by (\ref{Def_T}) is indeed continuous from $X^3_{(0)}$ to $X^3_{(0)}$.

\begin{remark}
From~(\ref{Def_TB=B}) it already appears that the behaviour of the coefficients $\< \mu \varphi_1, \varphi_n \>$, and thus Hypothesis~\ref{hyp}, will play a crucial role.
\end{remark}

\section{Definition and properties of the transformation}
\label{sec_transformation}

In this section, we make rigorous the heuristic developed in the previous section. In subsection~\ref{subsec_BaseRiesz}, we prove that for a suitable choice of $\beta^1_n$ and $\beta^2_n$ then  $\mathcal{B}$ defined in (\ref{Def_B}) is a Riesz basis of $X^2_{(0)}$ where the functions $g_n^{ij}$ and $h_n^{ij}$ are defined by~(\ref{Formule_gn}) and (\ref{Relation_gn_hn}). This enables us to define the feedback law and the transformation $T$ using the relations~(\ref{Def_TB=B}) and (\ref{Def_T}).
However, this does not give enough regularity to prove that $T : X^3_{(0)} \to X^3_{(0)}$. We prove the extra regularity we need on the feedback laws in subsection~\ref{subsec_ExtraRegularity}. This leads, in subsection~\ref{subsec_RegulariteT}, to the expected regularity for the transformation $T$.

\subsection{Riesz basis property} \label{subsec_BaseRiesz}

Let us recall some results on Riesz basis.
\begin{defn}
Let $H$ be an Hilbert space and $\{g_n\}_{n\in \NN^*}\subset H$. We say that $\{g_n\}_{n\in \NN^*}$ is $\omega$-independent if
\[
\sum_{n\in \NN^*} a_n g_n =0, \text{ with } \{a_n\}_{n\in \NN^*}\subset \RR
\quad \Longrightarrow \quad
a_n=0, \forall n\in \NN^*.
\]
\end{defn}

\begin{thm} \cite[Theorem 15]{YoungBook} \label{thm_independent}
Let $H$ be a separable Hilbert space and let $\{e_n\}_{n\in \NN^*}$ be an orthonormal basis for $H$. If $\{g_n\}_{n\in \NN^*}$ is an $\omega$-independent sequence quadratically close to $\{e_n\}_{n\in \NN^*}$, that is
\[
\sum_{n\in \NN^*}\|e_n-g_n\|^2_{H} < +\infty,
\]
then $\{g_n\}_{n\in \NN^*}$ is a Riesz basis for $H$.
\end{thm}

\begin{thm} \label{thm_dense}
Let $H$ be a separable Hilbert space and let $\{e_n\}_{n\in \NN^*}$ be an orthonormal basis for $H$. If $\{g_n\}_{n\in \NN^*}$ is dense in $H$ and is quadratically close to $\{e_n\}_{n\in \NN^*}$, i.e.
\[
\sum_{n\in \NN^*}\|e_n-g_n\|^2_{H} < +\infty,
\]
then $\{g_n\}_{n\in \NN^*}$ is a Riesz basis for $H$.
\end{thm}
Let us provide a proof of Theorem \ref{thm_dense}, stated as a remark in \cite[Remark 2.1, p. 318]{Gohberg}.

\beginpf

Let us prove Theorem \ref{thm_dense} by contradiction. Suppose that the $\{g_n\}_{n\in \NN^*}$ is dense in $H$ and is quadratically close to $\{e_n\}_{n\in \NN^*}$ but that $\{g_n\}_{n\in \NN^*}$ is not a Riesz basis. Therefore, by Theorem \ref{thm_independent}, there must exist a non-zero sequence $\{a_n\}_{n\in \NN^*}\subset \RR$ such that
\[
\sum_{n\in \NN^*} a_n g_n =0.
\]
Since $\{g_n\}_{n\in \NN^*}$ is quadratically close to $\{e_n\}_{n\in \NN^*}$, there exists $N\in \NN^*$ such that
\[
\sum_{n\geq N+1}\|e_n-g_n\|^2_{H} < 1.
\]
Therefore, from \cite[Theorem 13]{YoungBook}, the $\{e_n\}_{1\leq n\leq N} \cup \{g_n\}_{n\geq N+1}$ is a Riesz basis of $H$. This implies that there exists $k\leq N$ such that $a_k \neq 0$. Hence, with the density assumption, we have
\[
H=\overline{\Span\{g_n \, | \, n\in \NN^*\}}=\overline{\Span  \{g_n \, | \,n\in \NN^* \setminus \{k\}\}}.
\]
Thus,
\[
\textrm{codim}(\Span  \{g_n \, | \, n\geq N+1 \})\leq N-1.
\]
However
\[
H/\overline{\Span  \{g_n \, | \, n\geq N+1 \}},
\] is isomorphic to $\Span  \{e_n \, | \, n\leq  N \}$, which is of dimension $N$, leading to a contradiction.

\endpf

We will use the previous criteria to prove the following proposition.
\begin{prop} \label{prop_baseRiesz}
Let $\beta_n^i$ being defined as~\eqref{Def_beta_jn}. Then,
\begin{equation*}
\B:=\left\{
\begin{pmatrix} g_n^{12} \\ g_n^{22} \end{pmatrix} ,
\begin{pmatrix} h_n^{12} \\ h_n^{22} \end{pmatrix}
 \: ; \: n \in \NN^* \right\}
\end{equation*}
is a Riesz basis of $X^2_{(0)}$ and
\begin{equation*}
\tilde{\B}:=\left\{
\begin{pmatrix} g_n^{12}/\lambda_n^{1/2} \\ g_n^{22}/\lambda_n^{1/2} \end{pmatrix} ,
\begin{pmatrix} h_n^{12}/\lambda_n^{1/2} \\ h_n^{22}/\lambda_n^{1/2} \end{pmatrix}
 \: ; \: n \in \NN^* \right\}
\end{equation*}
is a Riesz basis of $X^3_{(0)}$.
\end{prop}

To apply the previous criterion for the Riesz basis, we prove that $\B$ (resp. $\tilde{\B}$) is quadratically close to the orthonormal basis of $X^2_{(0)}$ (resp. $X^3_{(0)}$) given by
\[
\left\{ \begin{pmatrix} \varphi_n/\lambda_n^{s/2} \\ 0 \end{pmatrix},
\begin{pmatrix} 0 \\ \varphi_n/\lambda_n^{s/2} \end{pmatrix} \: ; \: n \in \NN^* \right\},
\quad \text{with } s=2 \: \text{ (resp. $s=3$)}.
\]
Thus, we choose $\beta^1_n$ and $\beta^2_n$ such that
$g_n^{12}$ and $h_n^{22}$ are close to $\varphi_n/\lambda_n$ i.e. $\< g_n^{12}, \varphi_n \>= \<h_n^{22}, \varphi_n \>= 1/\lambda_n$, that is,
\begin{equation} \label{Def_beta_jn}
\begin{aligned}
\beta_n^1&:=\dfrac{\lambda(\lambda^2+4\lambda^2_n)}{2\lambda^3_n \< \mu \varphi_1,\varphi_n \>},
\\
\beta_n^2&:=-\dfrac{\lambda(\lambda^2+ 4\lambda^2_n)}{(\lambda^2+2\lambda^2_n)\lambda_n \<\mu \varphi_1,\varphi_n \>}.
\end{aligned}
\end{equation}
Notice that, from Hypothesis \ref{hyp} and Remark \ref{behavehyp}, there exists $C>0$ such that
\begin{equation} \label{comportement_beta_jn}
\frac{n}{C} \leq \left| \beta^1_n \right| \leq Cn,
\quad
\frac{n}{C} \leq \left| \beta^2_n \right| \leq Cn,
\quad
\forall n \in \NN^*.
\end{equation}

During the proof of Proposition~\ref{prop_baseRiesz}, we will use the following lemma. Its proof is purely technical and postponed to Appendix~\ref{Annexe_Riesz}.
\begin{lem} \label{Lemme_quadratiquement_proche}
With the above definition, for $s=2$ and $s=3$, one gets
\begin{gather*}
\sum_{n\in \NN^*} \left\|
\left(
\begin{array}{c}
\varphi_n/\lambda_n^{s/2}\\
0
\end{array}
\right)
-
\left(
\begin{array}{c}
g^{12}_n/\lambda_n^{(s-2)/2} \\
g^{22}_n/\lambda_n^{(s-2)/2}
\end{array}
\right)
\right\|^2_{{X^{s}_{(0)}}} <+\infty,
\\
\sum_{n\in \NN^*} \left\|
\left(
\begin{array}{c}
0 \\
\varphi_n/\lambda_n^{s/2}
\end{array}
\right)
-
\left(
\begin{array}{c}
h^{12}_n/\lambda_n^{(s-2)/2}\\
h^{22}_n/\lambda_n^{(s-2)/2}
\end{array}
\right)
\right\|^2_{{X^s_{(0)}}} <+\infty.
\end{gather*}
\end{lem}
We are now ready to prove Proposition~\ref{prop_baseRiesz}.

\beginpf
Let $s=2$ or $s=3$. In view of  Theorem~\ref{thm_independent} and Lemma~\ref{Lemme_quadratiquement_proche}, assume that there exists a sequence $\{a_n, b_n\}_{n\in \NN^*}$ such that
\begin{equation}\label{dependant}
\sum_{n\in \NN^*} a_n \left(
\begin{array}{c}
g^{12}_n \\
g^{22}_n
\end{array}
\right)+b_n
\left(
\begin{array}{c}
h^{12}_n \\
h^{22}_n
\end{array}
\right)
=
\left(
\begin{array}{c}
0 \\
0
\end{array}
\right) \text{ in } X^s_{(0)}.
\end{equation}

\medskip
\emph{First step: } we apply negative powers of the Laplacian to characterize elements of $S := \overline{\Span \B}$.

Recall that $\A$ is defined in~\eqref{DefOperateurA}. Let us write \eqref{Systeme_fn} for $g_n^{ij}$ as
\begin{align*}
\A g_n=\underbrace{\left[
\lambda_n \left( \begin{array}{cccc}
0 & 0 & 0 & 1 \\
0 & 0 & -1 & 0 \\
0 & -1 & 0 & 0 \\
1 & 0 & 0 & 0 \\
\end{array}
\right)
+\lambda \left( \begin{array}{cccc}
0 & 0 & -1 & 0 \\
0 & 0 & 0 & -1 \\
1 & 0 & 0 & 0 \\
0 & 1 & 0 & 0 \\
\end{array}
\right)
\right]}_{=:J_n} g_n  - \underbrace{
\left(
\begin{array}{c}
\beta_n^1 \mu \varphi_1 \\
0  \\
0\\
0
\end{array}
\right)}_{=:F_n},
\end{align*}
where $g_n=(g_n^{11},g_n^{12},g_n^{21},g_n^{22}){\tr}$ and where $\A g_n=(\A g_n^{11},\A g_n^{12},\A g_n^{21},\A g_n^{22}){\tr}$.

Since $\det(J_n)=(\lambda^2+\lambda^2_n)^2 \neq 0$, we have
\[
\A^{-1} g_n=J_n^{-1}g_n+\A^{-1}J_n^{-1}F_n,
\]
where
\[
J_n^{-1}=\dfrac{1}{\lambda^2+\lambda_n^2}\left( \begin{array}{cccc}
0 & 0 & \lambda & \lambda_n \\
0 & 0 & -\lambda_n & \lambda \\
-\lambda & -\lambda_n & 0 & 0 \\
\lambda_n & -\lambda & 0 & 0 \\
\end{array}\right).
\]
Therefore, using \eqref{Relation_gn_hn}, we obtain
\begin{equation} \label{RelationA-1}
\begin{aligned}
\A^{-1} g_n^{12}&=\dfrac{1}{\lambda_n^2+\lambda^2}\left[\lambda g_n^{22}-\lambda_n \dfrac{\beta_n^1}{\beta_n^2}h_n^{22}\right], \\
\A^{-1} g_n^{22}&=\dfrac{1}{\lambda_n^2+\lambda^2}\left[\lambda_n \dfrac{\beta_n^1}{\beta_n^2}h_n^{12} - \lambda g_n^{12}+\lambda_n\beta_n^1 \A^{-1} (\mu \varphi_1)\right], \\
\A^{-1} h_n^{12}&=\dfrac{1}{\lambda_n^2+\lambda^2}\left[\lambda_n \dfrac{\beta_n^2}{\beta_n^1} g_n^{22}+\lambda h_n^{22}\right], \\
\A^{-1} h_n^{22}&=\dfrac{-1}{\lambda_n^2+\lambda^2}\left[\lambda_n \dfrac{\beta_n^2}{\beta_n^1} g_n^{12}+\lambda h_n^{12}+\lambda \beta_n^2 \A^{-1} (\mu \varphi_1)\right].
\end{aligned}
\end{equation}
Thus, applying $\A^{-1}$ to~\eqref{dependant} leads to
\begin{equation} \label{A-1independent}
\begin{pmatrix} c_0 \A^{-1} (\mu \varphi_1) \\ 0 \end{pmatrix} =
\sum_{n \in \NN^*}  \frac{\lambda a_n + \lambda_n \frac{\beta_n^2}{\beta_n^1} b_n}{\lambda_n^2+\lambda^2} \begin{pmatrix} g_n^{12} \\ g_n^{22} \end{pmatrix} +
\frac{-\lambda_n \frac{\beta_n^1}{\beta_n^2}a_n + \lambda b_n}{\lambda_n^2 + \lambda^2}
\begin{pmatrix} h_n^{12} \\ h_n^{22} \end{pmatrix},
\end{equation}
where
\begin{equation} \label{Def_c0}
c_0 := \sum_{n \in \NN^*} \frac{\lambda_n \beta_n^1 a_n - \lambda \beta_n^2 b_n}{\lambda_n^2 + \lambda^2}.
\end{equation}
Applying $\A^{-1}$ to~\eqref{A-1independent} we obtain
\begin{equation} \label{A-2independent}
\begin{pmatrix} c_1 \A^{-1} (\mu \varphi_1) \\ c_0 \A^{-2} (\mu \varphi_1) \end{pmatrix}
= \sum_{n \in \NN^*} \frac{(\lambda^2-\lambda_n^2)a_n + 2\lambda \lambda_n \frac{\beta_n^2}{\beta_n^1}b_n}{(\lambda_n^2 + \lambda^2)^2} \begin{pmatrix} g_n^{12} \\ g_n^{22} \end{pmatrix}
+ \frac{-2 \lambda \lambda_n \frac{\beta_n^1}{\beta_n^2}a_n + (\lambda^2-\lambda_n^2) b_n}{(\lambda_n^2 + \lambda^2)^2} \begin{pmatrix} h_n^{12} \\ h_n^{22} \end{pmatrix},
\end{equation}
where
\begin{equation} \label{Def_c1}
c_1 := \sum_{n \in \NN^*} \frac{2 \lambda \lambda_n \beta_n^1 a_n + (\lambda_n^2-\lambda^2) \beta_n^2 b_n}{(\lambda_n^2+\lambda^2)^2}.
\end{equation}

In order to iterate \eqref{A-1independent} and \eqref{A-2independent}, notice that applying $\A^{-1}$ to the relations~\eqref{RelationA-1} leads to
\begin{equation} \label{RelationA-2}
\begin{aligned}
&\A^{-2} g_n^{12}=\dfrac{1}{(\lambda_n^2+\lambda^2)^2}\Big[(\lambda_n^2-\lambda^2) g_n^{12}+2\lambda\lambda_n \frac{\beta_n^1}{\beta_n^2}h_n^{12} +2\lambda \lambda_n \beta_n^1 \A^{-1} (\mu \varphi_1)\Big], \\
&\A^{-2} g_n^{22}=\dfrac{1}{(\lambda_n^2+\lambda^2)^2}\Big[(\lambda_n^2-\lambda^2) g_n^{22}+2\lambda \lambda_n \frac{\beta_n^1}{\beta_n^2}h_n^{22}\Big]  +\dfrac{\lambda_n \beta_n^1 }{\lambda_n^2+\lambda^2}\A^{-2} (\mu \varphi_1) ,   \\
&\A^{-2} h_n^{12}=\dfrac{1}{(\lambda_n^2+\lambda^2)^2}\Big[(\lambda_n^2-\lambda^2) h_n^{12}-2\lambda \lambda_n \frac{\beta_n^2}{\beta_n^1}g_n^{12}  +(\lambda_n^2-\lambda^2)\beta_n^2 \A^{-1} (\mu \varphi_1)\Big],    \\
&\A^{-2} h_n^{22}=\dfrac{1}{(\lambda_n^2+\lambda^2)^2}\Big[(\lambda_n^2-\lambda^2) h_n^{22}-2\lambda \lambda_n \frac{\beta_n^2}{\beta_n^1}g_n^{22}\Big] -\dfrac{\lambda \beta_n^2}{\lambda_n^2+\lambda^2} \A^{-2} (\mu \varphi_1).
\end{aligned}
\end{equation}
Applying successively $\A^{-2}$ to \eqref{A-1independent} we obtain, for every $p \in \mathbb{N}$ (with the convention that the sum from $1$ to $0$ is equal to 0),
\begin{equation} \label{A-impair/pair_independent}
\begin{aligned}
& \begin{pmatrix}
c_0 \A^{-(2p+1)} (\mu \varphi_1) + \sum_{j=1}^p \left( \sum_{n \in \NN^*} \frac{2 \lambda \lambda_n \beta_n^1 k_n^j + (\lambda_n^2-\lambda^2) \beta_n^2 l_n^j}{(\lambda_n^2+\lambda^2)^{2j+1}} \right) \A^{-2(p-j)-1}(\mu \varphi_1)
\\
\sum_{j=1}^p \left( \sum_{n \in \NN^*} \frac{\lambda_n \beta_n^1 k_n^j - \lambda \beta_n^2 l_n^j}{(\lambda_n^2+\lambda^2)^{2j}} \right) \A^{-2(p-j+1)}(\mu \varphi_1)
\end{pmatrix}
\\
&= \sum_{n \in \NN^*} \frac{1}{(\lambda_n^2+\lambda^2)^{2p+1}} \left[ k_n^{p+1} \begin{pmatrix} g_n^{12} \\ g_n^{22} \end{pmatrix} + l_n^{p+1} \begin{pmatrix} h_n^{12} \\ h_n^{22} \end{pmatrix} \right],
\end{aligned}
\end{equation}
where the coefficients $(k_n^j, l_n^j)$ are defined by
\begin{equation} \label{Recurrence_k_n_l_n}
\begin{pmatrix} k_n^{j+1} \\ l_n^{j+1} \end{pmatrix}
= \begin{pmatrix} (\lambda^2-\lambda_n^2) & 2 \lambda \lambda_n \frac{\beta_n^2}{\beta_n^1} \\
-2 \lambda \lambda_n \frac{\beta_n^1}{\beta_n^2} & (\lambda^2- \lambda_n^2) \end{pmatrix}
\begin{pmatrix} k_n^j \\ l_n^j \end{pmatrix},
\quad \text{with }
\begin{pmatrix} k_n^1 \\ l_n^1 \end{pmatrix} := \begin{pmatrix} \lambda a_n + \lambda_n \frac{\beta_n^2}{\beta_n^1}b_n \\ -\lambda_n \frac{\beta_n^1}{\beta_n^2} a_n +\lambda b_n \end{pmatrix}.
\end{equation}
Applying successively $\A^{-2}$ to \eqref{A-2independent} we obtain, for every $p \in \NN^*$,
\begin{equation} \label{A-impair/pair2_independent}
\begin{aligned}
& \begin{pmatrix}
c_1 \A^{-(2p+1)} (\mu \varphi_1) + \sum_{j=1}^p \left( \sum_{n \in \NN^*} \frac{2 \lambda \lambda_n \beta_n^1 \tilde{k}_n^j + (\lambda_n^2-\lambda^2) \beta_n^2 \tilde{l}_n^j}{(\lambda_n^2+\lambda^2)^{2(j+1)}} \right) \A^{-2(p-j)-1}(\mu \varphi_1)
\\
c_0 \A^{-2(p+1)} \sum_{j=1}^p \left( \sum_{n \in \NN^*} \frac{\lambda_n \beta_n^1 \tilde{k}_n^j - \lambda \beta_n^2 \tilde{l}_n^j}{(\lambda_n^2+\lambda^2)^{2j+1}} \right) \A^{-2(p-j+1)}(\mu \varphi_1)
\end{pmatrix}
\\
&= \sum_{n \in \NN^*} \frac{1}{(\lambda_n^2+\lambda^2)^{2(p+1)}} \left[ \tilde{k}_n^{p+1} \begin{pmatrix} g_n^{12} \\ g_n^{22} \end{pmatrix} + \tilde{l}_n^{p+1} \begin{pmatrix} h_n^{12} \\ h_n^{22} \end{pmatrix} \right],
\end{aligned}
\end{equation}
where the coefficients $(\tilde{k}_n^j, \tilde{l}_n^j)$ are defined by
\begin{equation} \label{Recurrence_ktilde_n_ltilde_n}
\begin{pmatrix} \tilde{k}_n^{j+1} \\ \tilde{l}_n^{j+1} \end{pmatrix}
= \begin{pmatrix} (\lambda^2-\lambda_n^2) & 2 \lambda \lambda_n \frac{\beta_n^2}{\beta_n^1} \\
-2 \lambda \lambda_n \frac{\beta_n^1}{\beta_n^2} & (\lambda^2- \lambda_n^2) \end{pmatrix}
\begin{pmatrix} \tilde{k}_n^j \\ \tilde{l}_n^j \end{pmatrix},
\quad \text{with }
\begin{pmatrix} \tilde{k}_n^0 \\ \tilde{l}_n^0 \end{pmatrix} = \begin{pmatrix} a_n \\ b_n \end{pmatrix}.
\end{equation}

Assuming that $c_0 \neq 0$, equality~\eqref{A-1independent} implies that
\[
\begin{pmatrix} \A^{-1} (\mu \varphi_1) \\ 0 \end{pmatrix} \in \overline{\Span \B}.
\]
Then using~\eqref{A-2independent} we obtain
\[
\begin{pmatrix} 0 \\ c_0 \A^{-2} (\mu \varphi_1) \end{pmatrix}
+ \begin{pmatrix} c_1 \A^{-1} (\mu \varphi_1) \\ 0 \end{pmatrix} \in \overline{\Span \B}
\quad \Longrightarrow \quad \begin{pmatrix} 0 \\ \A^{-2} (\mu \varphi_1) \end{pmatrix} \in \overline{\Span \B}.
\]
Iterating with the relations~\eqref{A-impair/pair_independent} and~\eqref{A-impair/pair2_independent} it comes that, for every $p \in \NN^*$,
\begin{equation} \label{ElementsDuSpan}
\begin{pmatrix} \A^{-(2p-1)} (\mu \varphi_1) \\ 0 \end{pmatrix} \in \overline{\Span \B}
\quad \text{ and }
\begin{pmatrix} 0 \\ \A^{-2p}(\mu \varphi_1) \end{pmatrix}  \in \overline{\Span \B}.
\end{equation}
Notice that if $c_0 = 0$ and $c_1 \neq 0$, the same argument can be repeated to obtain~\eqref{ElementsDuSpan}. Actually, one gets~\eqref{ElementsDuSpan} as soon as there exists a non-zero coefficient in the left-hand side of~\eqref{A-impair/pair_independent} or in the left-hand side of~\eqref{A-impair/pair2_independent}.
Thus it comes that either~\eqref{ElementsDuSpan} holds or, for every $j \in \NN^*$,
\begin{equation} \label{TousLesCoeffsNuls}
\begin{aligned}
0 &= \sum_{n \in \NN^*} \frac{2 \lambda \lambda_n \beta_n^1 k_n^j + (\lambda_n^2-\lambda^2) \beta_n^2 l_n^j}{(\lambda_n^2+\lambda^2)^{2j+1}}, \\
0&= \sum_{n \in \NN^*} \frac{\lambda_n \beta_n^1 k_n^j - \lambda \beta_n^2 l_n^j}{(\lambda_n^2+\lambda^2)^{2j}}, \\
0&= \sum_{n \in \NN^*} \frac{2 \lambda \lambda_n \beta_n^1 \tilde{k}_n^j + (\lambda_n^2-\lambda^2) \beta_n^2 \tilde{l}_n^j}{(\lambda_n^2+\lambda^2)^{2(j+1)}}, \\
0&= \sum_{n \in \NN^*} \frac{\lambda_n \beta_n^1 \tilde{k}_n^j - \lambda \beta_n^2 \tilde{l}_n^j}{(\lambda_n^2+\lambda^2)^{2j+1}}.
\end{aligned}
\end{equation}

\medskip
\emph{Second step: } we prove that if~\eqref{ElementsDuSpan} holds, then $\B$ is a Riesz basis.

Let $(d^1, d^2)\tr \in X^s_{(0)}$ such that
\[
\left\< \begin{pmatrix} d^1 \\ d^2 \end{pmatrix} , \begin{pmatrix} \phi^1 \\ \phi^2 \end{pmatrix} \right\>=0, \quad  \forall \begin{pmatrix} \phi^1 \\ \phi^2 \end{pmatrix} \in \overline{\Span \B}.
\]
Using~\eqref{ElementsDuSpan}, we obtain that, for every $p \in \NN^*$,
\begin{align*}
0 &= \left\< \begin{pmatrix} d^1 \\ d^2 \end{pmatrix} , \begin{pmatrix} \A^{-(2p-1)}(\mu \varphi_1)	 \\ 0 \end{pmatrix} \right\>_{X^s_{(0)}}
\\
&= \< d^1, \A^{-(2p-1)}(\mu \varphi_1) \>_{H^s_{(0)}}
\\
&= \sum_{k \in \NN^*} \lambda_k^s \<d^1, \varphi_k\>  \< \varphi_k, \A^{-(2p-1)}(\mu \varphi_1) \>
\\
&= \sum_{k \in \NN^*} \lambda_k^{s+1} \frac{\< \mu \varphi_1, \varphi_k \>  \<d^1, \varphi_k \>}{\lambda_k^{2p}}.
\end{align*}
Let $d_k^1 := \< d^1, \varphi_k \>$ and define
\begin{equation*}
G : z \in \CC \mapsto \sum_{k \in \NN^* } \lambda_k^{s-1} \< \mu \varphi_1, \varphi_k \> d_k^1 e^{z/ \lambda_k^2} \in \CC.
\end{equation*}
From uniform convergence on compact sets, it comes that $G$ is an entire function and the previous relation imply that, for every $p \in \NN^*$,
\[
G^{(p-1)}(0) = \sum_{k \in \NN^*} \lambda_k^{s+1} \frac{\< \mu \varphi_1, \varphi_k \>  d_k^1}{\lambda_k^{2p}} =0.
\]
Thus, $G \equiv 0$.
If $d^1 \neq 0$, let
\[n_0 := \min \{ n \in \NN^* ; d_n^1 \neq 0 \}.
\]
It comes that
\begin{equation}
\label{sim=0-z}
0 = \lambda_{n_0}^{s-1} \< \mu \varphi_1, \varphi_{n_0} \> d_{n_0}^1 + \sum_{k > n_0 } \lambda_k^{s-1} \< \mu \varphi_1, \varphi_k \> d_k^1 \exp \left( z \left( \frac{1}{\lambda_k^2} - \frac{1}{\lambda_{n_0}^2} \right) \right) .
\end{equation}
As $\< \mu \varphi_1, \varphi_{n_0} \> \neq 0$, considering $z$ real and letting $z$ go to $+\infty$
 in \eqref{sim=0-z} imply $d_{n_0}^1 =0$. Therefore, from the definition of $n_0$, $d^1=0$.

Using the exact same strategy for
\[
\begin{pmatrix} 0 \\ \A^{-2p}(\mu \varphi_1) \end{pmatrix}  \in \overline{\Span \B},
\quad \forall p \in \NN^*,
\]
implies $d^2 = 0$. Thus, it comes that $d^1 = d^2 = 0$ i.e. $\overline{\Span \B} = X^s_{(0)}$. From Theorem~\ref{thm_dense}, we obtain that $\B$ is a Riesz basis of $X^2_{(0)}$ (resp. $\tilde{\B}$ is a Riesz basis of $X^3_{(0)}$).

\medskip
\emph{Third step: } we prove that in the remaining case~\eqref{TousLesCoeffsNuls}, one has $a_n = b_n=0$ for any $n \in \NN^*$.

Let us define
\begin{equation} \label{Def_Gtilde}
\tilde{G} : z \in \CC \mapsto \sum_{n \in \NN^*} \frac{1}{(\lambda_n^2+\lambda^2)^2}
\left\< \begin{pmatrix} 2 \lambda \lambda_n \beta_n^1 \\ (\lambda_n^2 - \lambda^2) \beta_n^2  \end{pmatrix} , \exp \left( z M \right) \begin{pmatrix} a_n \\ b_n \end{pmatrix} \right\>,
\end{equation}
with
\begin{equation} \label{DefM}
M := \frac{1}{(\lambda_n^2 + \lambda^2)^2}
\begin{pmatrix} (\lambda^2-\lambda_n^2) & 2 \lambda \lambda_n \frac{\beta_n^2}{\beta_n^1} \\
-2 \lambda \lambda_n \frac{\beta_n^1}{\beta_n^2} & (\lambda^2- \lambda_n^2) \end{pmatrix}.
\end{equation}
Notice that the matrix appearing in this definition is the one used in the recurrence relations~\eqref{Recurrence_k_n_l_n} and~\eqref{Recurrence_ktilde_n_ltilde_n}.
This matrix is diagonalized as follows
\[
M = \frac{1}{(\lambda_n^2 + \lambda^2)^2} \begin{pmatrix} -i \frac{\beta_n^2}{\beta_n^1} & i\frac{\beta_n^2}{\beta_n^1} \\ & \\ 1 & 1 \end{pmatrix}
\begin{pmatrix} (\lambda^2-\lambda_n^2)+ 2i \lambda \lambda_n & 0 \\ 0 & (\lambda^2-\lambda_n^2)- 2i \lambda \lambda_n \end{pmatrix}
\begin{pmatrix} \frac{i}{2} \frac{\beta_n^1}{\beta_n^2} & \frac{1}{2}  \\ & \\  \frac{-i}{2} \frac{\beta_n^1}{\beta_n^2}  &  \frac{1}{2}  \end{pmatrix}.
\]
From this diagonalization we deduce that
\begin{align}\label{Gtilde}
\tilde{G}(z) = \sum_{n \in \NN^*}&
\exp \left( \frac{(\lambda+i\lambda_n)^2}{(\lambda^2+\lambda_n^2)^2} z \right)\Bigg[-\frac{i\beta_n^1 a_n}{2}
\frac{(\lambda-i\lambda_n)^2}{(\lambda_n^2+\lambda^2)^2}-\frac{\beta_n^2 b_n}{2}\frac{(\lambda+i\lambda_n)^2}{(\lambda_n^2+\lambda^2)^2}\Bigg] \nonumber \\
+&\exp \left( \frac{(\lambda-i\lambda_n)^2}{(\lambda^2+\lambda_n^2)^2} z \right)\Bigg[\frac{i\beta_n^1 a_n}{2}
\frac{(\lambda+i\lambda_n)^2}{(\lambda_n^2+\lambda^2)^2}+\frac{\beta_n^2 b_n}{2}\frac{(\lambda-i\lambda_n)^2}{(\lambda_n^2+\lambda^2)^2}\Bigg] .
\end{align}

Again, $\tilde{G}$ is an entire function from the uniform convergence on compact sets. The recurrence relation~\eqref{Recurrence_ktilde_n_ltilde_n} implies that, for every $p \in \mathbb{N}$,
\begin{align*}
\tilde{G}^{(p)}(0)
&= \sum_{n \in \NN^*} \frac{1}{(\lambda_n^2+\lambda^2)^2}
\left\< \begin{pmatrix} 2 \lambda \lambda_n \beta_n^1 \\ (\lambda_n^2 - \lambda^2) \beta_n^2  \end{pmatrix} , M^p \begin{pmatrix} a_n \\ b_n \end{pmatrix} \right\>
\\
&= \sum_{n \in \NN^*} \frac{2 \lambda\lambda_n \beta_n^1 \tilde{k}_n^p + (\lambda_n^2 -\lambda^2) \beta_n^2 \tilde{l}_n^p}{(\lambda_n^2 +\lambda^2)^{2(p+1)}}
\\
&=0.
\end{align*}
Thus we obtain
\begin{equation}\label{tildeG=0}
\tilde{G} \equiv 0.
\end{equation}
 We claim that
\begin{equation}\label{freq}
\frac{(\lambda - i\lambda_n)^2}{(\lambda_n^2+\lambda^2)^2}\neq \frac{(\lambda + i\lambda_m)^2}{(\lambda_m^2+\lambda^2)^2}, \quad (n,m)\in (\NN^*)^2.
\end{equation}
Indeed, let $(n,m)\in (\NN^*)^2$ be such that
\begin{equation}\label{freq-eq}
\frac{(\lambda - i\lambda_n)^2}{(\lambda_n^2+\lambda^2)^2}\neq \frac{(\lambda + i\lambda_m)^2}{(\lambda_m^2+\lambda^2)^2}.
\end{equation}
Taking the modulus of both sides of \eqref{freq-eq}, one gets
\[
\frac{1}{\lambda_n^2+\lambda^2}=\frac{1}{\lambda_m^2+\lambda^2},
\]
which implies that $\lambda_n=\lambda_m$ and therefore $n=m$. Hence
\begin{equation}
\label{casnn}
\frac{(\lambda - i\lambda_n)^2}{(\lambda_n^2+\lambda^2)^2}= \frac{(\lambda + i\lambda_n)^2}{(\lambda_n^2+\lambda^2)^2},
\end{equation}
Taking the imaginary part of \eqref{casnn} yields a contradiction and proves that
\eqref{freq} holds.

 For simplicity, we rewrite $G$ as
\[
G(z)=\sum_{n\in \NN^*} C_ne^{\mu_n z},
\]
with $\mu_n$ equal to $(\lambda - i\lambda_k)^2/(\lambda_k^2+\lambda^2)^2$ or $(\lambda + i\lambda_k)^2/(\lambda_k^2+\lambda^2)^2$ for some $k\in \NN^*$ and $C_n$ is the corresponding coefficient in \eqref{Gtilde}. Notice that $\mu_n$ are all different, $\mu_n \rightarrow 0$ as $n \rightarrow \infty$ and $C_n \in \ell^2(\NN^*;\CC)$.

We repeat the same argument as in the finite dimensional case. Let
\[
\cal{C}:=\textrm{Conv}\{\mu_n \, | \, n\in \NN^* \textrm{ such that } C_n \neq 0\}.
\]
Consider a nonzero extremal point of $\cal{C}$, which is therefore of the form $\mu_{n_0}$ for some $n_0\in \NN^*$. Hence, there exists $\theta \in [0,2\pi]$ such that
\begin{equation}\label{ordrefreq}
\Re\left(e^{i\theta}\mu_n\right) < \Re\left(e^{i\theta}\mu_{n_0}\right), \quad \forall n \in  \NN^*\setminus \{n_0\}
\end{equation}
By letting $z=\rho e^{i\theta}$ with $\rho \in \RR$ and using \eqref{tildeG=0}, we obtain
\begin{equation}
\label{Grhotheta=0}
0=G\left(\rho e^{i\theta}\right) e^{-\rho e^{i\theta}\mu_{n_0}}= C_{n_0}+\sum_{n \in \NN^*\setminus\{n_0\}}
  C_n e^{\rho\left(e^{i\theta}\mu_n-e^{i\theta}\mu_{n_0}\right)}.
\end{equation}
We then let $\rho \rightarrow +\infty$ in \eqref{Grhotheta=0} to obtain, using \eqref{ordrefreq}, that $C_{n_0}=0$ which is a contradiction with the construction $\cal{C}$. It implies that $C_n=0, \forall n\in \NN$. The expression of the $C_n$ implies for all $n\in \NN^*$
\begin{align*}
-\frac{i\beta_n^1 a_n}{2}
\frac{(\lambda-i\lambda_n)^2}{(\lambda_n^2+\lambda^2)^2}-\frac{\beta_n^2 b_n}{2}\frac{(\lambda+i\lambda_n)^2}{(\lambda_n^2+\lambda^2)^2}=0,\\
\frac{i\beta_n^1 a_n}{2}
\frac{(\lambda+i\lambda_n)^2}{(\lambda_n^2+\lambda^2)^2}+\frac{\beta_n^2 b_n}{2}\frac{(\lambda-i\lambda_n)^2}{(\lambda_n^2+\lambda^2)^2}=0.
\end{align*}
One then easily concludes that $a_n=b_n=0$ by using \eqref{freq-eq} and  the fact that $\beta_n^1\neq 0$ and $\beta_n^2\neq 0$ for all $n\in \NN^*$. Theorem~\ref{thm_independent} thus implies the Riesz basis property.

\endpf

\subsection{Definition and regularity of the feedback law} \label{subsec_ExtraRegularity}

So far, we have obtained from \eqref{Systeme_fn} the expression of the kernels $k^{ij}$ with respect to the Fourier coefficients $\alpha_n^i$ of the feedback
\begin{equation} \label{DefK}
K \begin{pmatrix} \Psi^1 \\ \Psi^2 \end{pmatrix} := \sum_{n=1}^{+\infty} \alpha^1_n \lag \Psi^1, \varphi_n \rag + \alpha^2_n \lag \Psi^2, \varphi_n \rag.
\end{equation}
The regularity of the kernels and, consequently, the regularity of $T$, will be directly related to the decay rate of those coefficients as $n\rightarrow +\infty$. It remains to use the $TB=B$ condition~(\ref{Def_TB=B}) to construct $K$ and $T$.

As $\mu \varphi_1 \in H^2_{(0)}$ and $\B$ is a Riesz basis of $X^2_{(0)}$ (see Proposition~\ref{prop_baseRiesz}), it comes that there exist sequences $(a^1_n)_{n \in \NN^*}$, $(a^2_n)_{n \in \NN^*} \in \ell^2(\NN^*, \RR)$ such that
\begin{equation} \label{decomposition_muPhi1}
\begin{pmatrix}
0 \\ \mu \varphi_1
\end{pmatrix}
= \sum_{n=1}^{+\infty}
a^1_n \begin{pmatrix} g_n^{12} \\ g_n^{22} \end{pmatrix}
+ a^2_n \begin{pmatrix} h_n^{12} \\ h_n^{22} \end{pmatrix}, \quad \text{in } X^2_{(0)}.
\end{equation}
Getting back to the $TB=B$ condition~(\ref{Def_TB=B}), we define
\begin{equation} \label{def_alpha1n_alpha2n}
\alpha^1_n := \frac{\beta^1_n}{\< \mu \varphi_1,\varphi_n \>} a^1_n,
\qquad
\alpha^2_n := \frac{\beta^2_n}{\< \mu \varphi_1,\varphi_n \>} a^2_n.
\end{equation}
Then, following the heuristic of Sec.~\ref{sec_heuristique}, the transformation $T$ is defined by~(\ref{Def_T}) and the feedback law $K$ is defined by~(\ref{DefK}).

Unfortunately, the regularity of the coefficients is only $(a^j_n)_{n \in \NN^*} \in \ell^2(\NN^*,\RR)$ for the $X^2_{(0)}$ Riesz basis and will not be sufficient to prove that $T$ is continuous in $X^3_{(0)}$.

\begin{remark} \label{Rk1_compacite}
Recall that we assumed controllability of the linearized system i.e. Hypothesis~\ref{hyp} . From~\cite[Remark 2]{BeauchardLaurent}, it follows that $\mu'(0) \pm \mu'(1) \neq 0$ and then that $(\mu \varphi_1) \not\in H^3_{(0)}$. Thus, $(0, \mu \varphi_1)\tr$ cannot be decomposed in the Riesz basis of $X^3_{(0)}$. This would have led to more regularity for the sequences $(\alpha_n^j)_n$.

It is fortunate since, as it will be noticed in Remark~\ref{rmq_inversibilite_regularite_coeff}, if $(\mu \varphi_1) \in H^3_{(0)}$, the obtained transformation $T$ would have been compact in $X^3_{(0)}$ and thus not invertible in $X^3_{(0)}$.
\end{remark}

Performing a suitable decomposition of the function $\mu \varphi_1$, we prove that the coefficients of the feedback law satisfy the following regularity.
\begin{prop} \label{prop_regulariteFeedback}
We define the sequence $(h_n)_{n \in \NN^*}$ by
\[
h_n := \frac{4}{n^3 \pi^2} \left( (-1)^{n+1} \mu'(1) - \mu'(0) \right).
\]
Then the sequences $(\alpha^1_n)_{n \in \NN^*}$ and $(\alpha^2_n)_{n \in \NN^*}$ defined in~(\ref{def_alpha1n_alpha2n}) satisfy
\begin{equation*}
\left(\frac{\alpha^1_n}{n^3} \right) \in \ell^2(\NN^*,\RR), \quad
\left(\frac{1}{n^3} \left( \alpha^2_n - \lambda_n \beta^2_n \frac{h_n}{\<\mu \varphi_1, \varphi_n\>} \right) \right) \in \ell^2(\NN^*,\RR).
\end{equation*}
\end{prop}

\begin{remark} \label{rmq_alphajn_borne}
As a corollary, it comes that, for every $j \in \{1,2\}$,
\begin{equation*}
\left( \frac{\alpha_n^j}{n^3} \right)_{n \in \NN^*} \in \ell^\infty(\NN^*,\RR).
\end{equation*}
\end{remark}

\beginpf

\emph{First step: splitting of the problem.}
We start using the ideas developed in~\cite{Puel13} to extend the regularization result \cite[Lemma 1]{BeauchardLaurent} to higher dimensions.
Let
\begin{equation*}
h : x \in (0,1) \mapsto -\frac{\pi \sqrt{2}}{3}  \left( (\mu'(0)+\mu'(1)) x^3 - 3 \mu'(0) x^2 + (2 \mu'(0) - \mu'(1))x \right).
\end{equation*}
It comes that
\begin{equation}\label{defg}
g := \mu \varphi_1 - h \in H^3_{(0)}.
\end{equation}
 The Fourier coefficients of $h$ are given by
\begin{equation} \label{coeff_fourier_h}
\< h, \varphi_k \> = h_k = \frac{4}{k^3 \pi^2} \left( (-1)^{k+1} \mu'(1) - \mu'(0) \right), \quad	 \forall k \in \NN^*,
\end{equation}
which is the leading term in the asymptotic expansion of $\lag \mu \varphi_1, \varphi_k\rag$ given in \cite[Remark 2]{BeauchardLaurent}. Let us split the left-hand side of \eqref{decomposition_muPhi1} in two parts
\[
\begin{pmatrix}
0 \\ \mu\varphi_1
\end{pmatrix}=
\begin{pmatrix}
0 \\ g
\end{pmatrix}+
\begin{pmatrix}
0 \\ h
\end{pmatrix}.
\]

As $g \in H^3_{(0)}$, using the Riesz basis of $X^3_{(0)}$ we get the existence of sequences $(\delta^1_n)_{n \in \NN^*}$ and $(\delta^2_n)_{n \in \NN^*}$ such that
\begin{equation} \label{decomposition_g}
\begin{pmatrix}
0 \\ g
\end{pmatrix}
= \sum_{n=1}^{+\infty} \<\mu\varphi_1, \varphi_n \> \left[
\frac{\lambda_n^{1/2} \delta^1_n}{\beta^1_n} \begin{pmatrix} g_n^{12}/ \lambda_n^{1/2} \\ g_n^{22}/ \lambda_n^{1/2} \end{pmatrix} +
\frac{\lambda_n^{1/2} \delta^2_n}{\beta^2_n} \begin{pmatrix} h_n^{12}/ \lambda_n^{1/2} \\ h_n^{22}/ \lambda_n^{1/2} \end{pmatrix}
\right].
\end{equation}
The coefficients of the decomposition in a Riesz basis being $\ell^2$ sequences, Hypothesis~\ref{hyp} and the behaviour of coefficients $\beta^j_n$ in~(\ref{comportement_beta_jn}) imply that, for every $j \in \{1, 2\}$,
\begin{equation} \label{comportement_aj}
\left( \frac{\delta_n^j}{n^3} \right)_{n \in \NN^*} \in \ell^2(\NN^*, \RR).
\end{equation}

\medskip
\emph{Second step: decomposition of $h$.}
Using the Riesz basis $\B$ of $X^2_{(0)}$, we get coefficients $(\rho^1_n)_{n \in \NN^*}$ and $(\rho^2_n)_{n \in \NN^*}$ such that
\begin{equation}\label{decomposition_h}
\begin{pmatrix}
0 \\ h
\end{pmatrix}
= \sum_{n=1}^{+\infty} \<\mu\varphi_1, \varphi_n \> \left[
\frac{\rho^1_n}{\beta^1_n} \begin{pmatrix} g_n^{12} \\ g_n^{22} \end{pmatrix} +
\frac{\rho^2_n}{\beta^2_n} \begin{pmatrix} h_n^{12} \\ h_n^{22} \end{pmatrix}
\right].
\end{equation}
Recall that the basis $\B$ is obtained as a perturbation of the $L^2$-orthonormal basis. To highlight this, we define the sequences $(\gamma_n^j)_{n \in \NN^*}$, $j\in\{1,2\}$, such that
\begin{equation} \label{Def_gamma_nj}
\begin{aligned}
\rho^1_n &= \lambda_n \beta^1_n \gamma^1_n,
\\
\rho^2_n& = \lambda_n \beta^2_n \frac{\<h, \varphi_n\>}{\<\mu \varphi_1, \varphi_n\>} + \lambda_n \beta^2_n \gamma^2_n.
\end{aligned}
\end{equation}

From (\ref{DecompositionFourier_gn_hn}), (\ref{decomposition_h}), and (\ref{Def_gamma_nj}),  we get,  for every $k \in \NN^*$,
\begin{align*}
\begin{pmatrix}
0 \\ \<h,\varphi_k\>
\end{pmatrix}
&= \sum_{n=1}^{+\infty} \<\mu \varphi_1, \varphi_n\> \<\mu \varphi_1, \varphi_k\>
\left[ \rho^1_n \begin{pmatrix} c_{nk}^{12} \\ c_{nk}^{22} \end{pmatrix}
+ \rho^2_n \begin{pmatrix} d_{nk}^{12} \\ d_{nk}^{22} \end{pmatrix} \right]
\\
&= \sum_{n=1}^{+\infty}  \<\mu \varphi_1, \varphi_n\> \<\mu \varphi_1, \varphi_k\>
\left( \frac{\< h,\varphi_n \>}{\< \mu \varphi_1, \varphi_n \>}  \lambda_n \beta^2_n \right)
\begin{pmatrix} d_{nk}^{12} \\ d_{nk}^{22} \end{pmatrix}
\\
&+ \sum_{n=1}^{+\infty}  \<\mu \varphi_1, \varphi_n\> \<\mu \varphi_1, \varphi_k\>
\left[ \lambda_n \beta^1_n \gamma^1_n \begin{pmatrix} c_{nk}^{12} \\ c_{nk}^{22} \end{pmatrix}
+ \lambda_n \beta^2_n \gamma^2_n \begin{pmatrix} d_{nk}^{12} \\ d_{nk}^{22} \end{pmatrix} \right].
\end{align*}
which, using \eqref{Formule_gn}, \eqref{defdeltank}, \eqref{Relation_gn_hn}, and (\ref{Def_beta_jn}), is equivalent to
\begin{align*}
&-\begin{pmatrix}
\sum_{n=1}^{+\infty}   \<\mu \varphi_1, \varphi_k\> \<h, \varphi_n\> \lambda_n \beta^2_n d_{nk}^{12}
\\
\sum_{n=1, n \neq k}^{+\infty}  \<\mu \varphi_1, \varphi_k\> \<h, \varphi_n\> \lambda_n \beta^2_n d_{nk}^{22}
\end{pmatrix}
\\
&= \sum_{n=1}^{+\infty}  \<\mu \varphi_1, \varphi_n\> \<\mu \varphi_1, \varphi_k\>
\left[ \lambda_n \beta^1_n \gamma^1_n \begin{pmatrix} c_{nk}^{12} \\ c_{nk}^{22} \end{pmatrix}
+ \lambda_n \beta^2_n \gamma^2_n \begin{pmatrix} d_{nk}^{12} \\ d_{nk}^{22} \end{pmatrix} \right],
\end{align*}
for every $k \in \NN^*$.
Hence, if we define
\begin{equation}
\label{deftildeh}
\begin{pmatrix}
\tilde{h}_1 \\ \tilde{h}_2
\end{pmatrix}
:= -\sum_{k=1}^{+\infty}
\begin{pmatrix}
\left(\sum_{n=1}^{+\infty} \<h, \varphi_n\> \<\mu \varphi_1, \varphi_k\> \lambda_n \beta^2_n d_{nk}^{12} \right) \varphi_k
\\
\left(\sum_{n=1, n \neq k}^{+\infty} \<h, \varphi_n\> \<\mu \varphi_1, \varphi_k\> \lambda_n \beta^2_n d_{nk}^{22} \right) \varphi_k
\end{pmatrix},
\end{equation}
we get, using also \eqref{DecompositionFourier_gn_hn},
\begin{align*}
\begin{pmatrix}
\tilde{h}_1 \\ \tilde{h}_2
\end{pmatrix}
&= \sum_{n=1}^{+\infty} \sum_{k=1}^{+\infty} \<\mu \varphi_1, \varphi_n\>
\left[ \lambda_n  \gamma^1_n \begin{pmatrix} \<\mu \varphi_1, \varphi_k\> \beta^1_n c_{nk}^{12} \\
\<\mu \varphi_1, \varphi_k\> \beta^1_n c_{nk}^{22} \end{pmatrix}
+ \lambda_n \gamma^2_n \begin{pmatrix} \<\mu \varphi_1, \varphi_k\> \beta^2_n d_{nk}^{12} \\
\<\mu \varphi_1, \varphi_k\> \beta^2_n d_{nk}^{22} \end{pmatrix} \right]\varphi_k
\\
&= \sum_{n=1}^{+\infty} \lambda_n^{3/2} \<\mu \varphi_1, \varphi_n\>
\left[ \gamma^1_n \begin{pmatrix} g_n^{12}/\lambda_n^{1/2} \\ g_n^{22}/\lambda_n^{1/2} \end{pmatrix}
+ \gamma^2_n \begin{pmatrix} h_n^{12}/\lambda_n^{1/2} \\ h_n^{22}/\lambda_n^{1/2} \end{pmatrix} \right].
\tag{\theequation} \addtocounter{equation}{1}
\label{ExtraRegularity_fonctionauxiliaire}
\end{align*}
 Finally, if $(\tilde{h}_1,\tilde{h}_2)\tr \in X^3_{(0)}$ (which will be proved in the next two steps), it comes that, for every $j \in \{1, 2\}$,
\begin{equation*}
\left( \lambda_n^{3/2} \< \mu \varphi_1, \varphi_n\>  \gamma_n^j \right)_n \in \ell^2(\NN^*,\RR),
\end{equation*}
which, thanks to Hypothesis~\ref{hyp}, implies that
\begin{equation}\label{estgamma}
(\gamma_n^j)_n \in \ell^2(\NN^*,\RR).
\end{equation}
 Using \eqref{decomposition_muPhi1},  \eqref{def_alpha1n_alpha2n}, \eqref{defg}, \eqref{decomposition_g}, \eqref{decomposition_h}, and \eqref{Def_gamma_nj}, we obtain
\begin{align*}
\frac{\alpha^1_n }{n^3} &= \frac{\rho^1_n}{n^3} + \frac{\delta^1_n}{n^3} = \frac{\lambda_n \beta^1_n}{n^3} \gamma^1_n + \frac{\delta^1_n}{n^3},
\\
\frac{\alpha^2_n }{n^3} &= \frac{\rho^2_n}{n^3} + \frac{\delta^2_n}{n^3} = \frac{\lambda_n \beta^2_n}{n^3} \gamma^2_n + \frac{\lambda_n \beta^2_n}{n^3} \frac{\< h , \varphi_n \>}{\< \mu \varphi_1, \varphi_n \>} + \frac{\delta^2_n}{n^3},
\end{align*}
which, with \eqref{esthyp1}, \eqref{comportement_beta_jn}, \eqref{coeff_fourier_h}, \eqref{comportement_aj}, and
\eqref{estgamma}  will end the proof of Proposition~\ref{prop_regulariteFeedback}.

\medskip
\emph{Third step: $H^3_{(0)}$ regularity of $\tilde{h}_1$.}
We start by proving that $\tilde{h}_1 \in H^3_{(0)}$ i.e.
\begin{equation} \label{H3regularity_tildeh1}
\sum_{k=1}^{+\infty} \left|k^3 \< \tilde{h}_1,\varphi_k \> \right|^2 = \sum_{k=1}^{+\infty} \left| k^3 \sum_{n=1}^{+\infty} \< \mu \varphi_1, \varphi_k \> \<h, \varphi_n\> \lambda_n \beta^2_n d_{nk}^{12} \right|^2 < +\infty.
\end{equation}
From (\ref{comportement_beta_jn}) and (\ref{coeff_fourier_h}) it comes that there exists $C >0$ such that
\begin{equation*}
|\<h ,\varphi_n\> \lambda_n \beta^2_n| \leq C, \quad \forall n \in \NN^*.
\end{equation*}
Using \eqref{boundmuphi1phik} together with \eqref{Formule_gn}, \eqref{defdeltank}, \eqref{Relation_gn_hn} and \eqref{DecompositionFourier_gn_hn},  it suffices to prove that
\begin{equation} \label{H3regularity_tildeh1bis}
M_{12}:=\sum_{k=1}^{+\infty} \left| \sum_{n=1}^{+\infty} d_{nk}^{12} \right|^2 =
\sum_{k=1}^{+\infty} \left| \sum_{n=1}^{+\infty} \frac{\lambda_k(\lambda^2+\lambda_k^2-\lambda_n^2)}{\left( \lambda^2 + (\lambda_k - \lambda_n)^2\right) \left( \lambda^2 + (\lambda_k + \lambda_n)^2 \right)} \right|^2 < +\infty.
\end{equation}
Notice that
\begin{equation*}
M_{12} \leq 2 \sum_{k=1}^{+\infty} \left[ \left( \frac{\lambda^2 \lambda_k}{\lambda^2(\lambda^2+4\lambda_k^2)} \right)^2
+ \left( \sum_{\underset{n \neq k}{n=1}}^{+\infty} d_{nk}^{12}
\right)^2 \right].
\end{equation*}
As
\begin{equation*}
\sum_{k=1}^{+\infty} \left( \frac{\lambda^2 \lambda_k}{\lambda^2(\lambda^2+4\lambda_k^2)} \right)^2  < +\infty,
\end{equation*}
we only deal with the second term.
Straightforward computations lead to
\begin{equation} \label{regulariteH30_1ere_decomposition}
\sum_{\underset{n \neq k}{n=1}}^{+\infty} d_{nk}^{12} =
\frac{1}{2} \sum_{\underset{n \neq k}{n=1}}^{+\infty} \left( \frac{\lambda_k-\lambda_n}{\lambda^2 + (\lambda_k - \lambda_n)^2} + \frac{\lambda_k+\lambda_n}{\lambda^2 + (\lambda_k + \lambda_n)^2} \right).
\end{equation}
The second term of this sum is easily dealt with. As $\lambda^2 + (\lambda_k + \lambda_n)^2 \geq (\lambda_k + \lambda_n)^2$ it comes that
\begin{align*}
\sum_{\underset{n \neq k}{n=1}}^{+\infty} \frac{\lambda_k+\lambda_n}{\lambda^2 + (\lambda_k + \lambda_n)^2}
&\leq \sum_{\underset{n \neq k}{n=1}}^{+\infty} \frac{1}{\lambda_k + \lambda_n}
\\
&\leq \sum_{n=1}^{+\infty} \int_{n-1}^n \frac{1}{\lambda_k + \pi^2 x^2} \md x
\\
&= \frac{1}{2k \pi}.
\end{align*}
Thus,
\begin{equation} \label{regulariteH30_1erterme}
\sum_{k=1}^{+\infty} \left(\sum_{\underset{n \neq k}{n=1}}^{+\infty} \frac{\lambda_k+\lambda_n}{\lambda^2 + (\lambda_k + \lambda_n)^2} \right)^2 < + \infty.
\end{equation}
We now turn back to the first term of the right-hand side of~(\ref{regulariteH30_1ere_decomposition}). We have
\begin{equation*}
\left| \sum_{\underset{n \neq k}{n=1}}^{+\infty} \frac{\lambda_k-\lambda_n}{\lambda^2+(\lambda_k-\lambda_n)^2} \right|
\leq \sum_{\underset{n \neq k}{n=1}}^{+\infty} \frac{1}{|\lambda_k-\lambda_n|}
= \sum_{n=1}^{k-1} \frac{1}{\lambda_k-\lambda_n} + \sum_{n=k+1}^{+\infty} \frac{1}{\lambda_n-\lambda_k}.
\end{equation*}
Straightforward computations lead to
\begin{align*}
\sum_{n=1}^{k-1} \frac{1}{\lambda_k-\lambda_n}
&= \frac{1}{2k \pi^2} \sum_{n=1}^{k-1} \left( \frac{1}{k-n} + \frac{1}{k+n} \right)
\\
&= \frac{1}{2 k \pi^2} \left( \sum_{j=1}^{k-1} \frac{1}{j} + \sum_{j=k+1}^{2k-1} \frac{1}{j} \right)
\\
&= \frac{1}{2k \pi^2} \left( \sum_{j=1}^{2k-1} \frac{1}{j} - \frac{1}{k} \right).
\end{align*}
In the same spirit, for $N\geq k+1$,
\begin{equation*}
0\geq \sum_{n=k+1}^{N} \frac{1}{\lambda_k-\lambda_n}
= \frac{1}{2k \pi^2} \sum_{j=2k+1}^{N+k} \frac{1}{j} - \frac{1}{2k \pi^2} \sum_{j=1}^{N-k} \frac{1}{j}
\end{equation*}
and thus,
\begin{equation*}
0\leq \sum_{n=k+1}^{+\infty} \frac{1}{\lambda_n -\lambda_k}
 \leq \frac{1}{2k \pi^2} \sum_{j=1}^{2k} \frac{1}{j} - \frac{1}{2k \pi^2} \sum_{j=n-k+1}^{N+k} \frac{1}{j}
  \leq \frac{1}{2k \pi^2} \sum_{j=1}^{2k} \frac{1}{j}.
\end{equation*}
 Finally,
\begin{equation} \label{H3regularite_estimeeIntermediaire}
\left| \sum_{\underset{n \neq k}{n=1}}^{+\infty}  \frac{\lambda_k-\lambda_n}{\lambda^2 + (\lambda_k - \lambda_n)^2} \right|
\leq \sum_{\underset{n \neq k}{n=1}}^{+\infty}  \frac{1}{| \lambda_k - \lambda_n |}
\leq \frac{1}{2k \pi^2} \left( 2 \sum_{j=1}^{2k-1} \frac{1}{j} - \frac{1}{2k} \right),
\end{equation}
which belongs to $\ell^2(\NN^*)$ as $\sum_{j=1}^{2k-1} \frac{1}{j} \underset{k \to +\infty}{\sim} \ln(2k-1)$. From~(\ref{regulariteH30_1ere_decomposition}), this proves that $\tilde{h}_1 \in H^3_{(0)}$.

\medskip
\emph{Fourth step: $H^3_{(0)}$ regularity of $\tilde{h}_2$.} We end the proof of Proposition~\ref{prop_regulariteFeedback} by proving that $\tilde{h}_2 \in H^3_{(0)}$ i.e., by \eqref{deftildeh},
\begin{equation} \label{H3regularity_tildeh2}
\sum_{k=1}^{+\infty} \left|k^3 \<\tilde{h}_2,\varphi_k\> \right|^2 = \sum_{k=1}^{+\infty} \left| k^3 \sum_{\underset{n \neq k}{n=1}}^{+\infty} \<\mu \varphi_1, \varphi_k\> \<h, \varphi_n\> \lambda_n \beta^2_n d_{nk}^{22} \right|^2 < +\infty.
\end{equation}
From (\ref{comportement_beta_jn}) and (\ref{coeff_fourier_h}) it comes that there exists $C>0$ such that
\begin{equation*}
|\< h ,\varphi_n\> \lambda_n \beta^2_n| \leq C, \quad \forall n \in \NN^*.
\end{equation*}
Using \eqref{boundmuphi1phik} it suffices to prove that
\begin{equation} \label{H3regularity_tildeh2bis}
\sum_{k=1}^{+\infty} \left| \sum_{\underset{n \neq k}{n=1}}^{+\infty} d_{nk}^{22} \right|^2 =
\sum_{k=1}^{+\infty} \left| \sum_{\underset{n \neq k}{n=1}}^{+\infty} \frac{\lambda(\lambda^2+\lambda_k^2+\lambda_n^2)}{\left( \lambda^2 + (\lambda_k - \lambda_n)^2\right) \left( \lambda^2 + (\lambda_k + \lambda_n)^2 \right)} \right|^2 < +\infty.
\end{equation}
Notice that
\begin{align*}
\frac{\lambda(\lambda^2+\lambda_k^2+\lambda_n^2)}{\left( \lambda^2 + (\lambda_k - \lambda_n)^2\right) \left( \lambda^2 + (\lambda_k + \lambda_n)^2 \right)}
&= \frac{\lambda}{2} \left( \frac{1}{\lambda^2 + (\lambda_k + \lambda_n)^2} + \frac{1}{\lambda^2 + (\lambda_k - \lambda_n)^2} \right)
\\
&\leq \frac{\lambda}{2} \left( \frac{\lambda_k + \lambda_n}{\lambda^2 + (\lambda_k + \lambda_n)^2} + \frac{|\lambda_k-\lambda_n|}{\lambda^2 + (\lambda_k - \lambda_n)^2} \right)
\\
&\leq \frac{\lambda}{2} \left( \frac{\lambda_k + \lambda_n}{\lambda^2 + (\lambda_k + \lambda_n)^2} + \frac{1}{|\lambda_k-\lambda_n|} \right).
\end{align*}
From~(\ref{regulariteH30_1erterme}) and~(\ref{H3regularite_estimeeIntermediaire}), it comes that
\begin{equation*}
\sum_{k=1}^{+\infty} \left| \sum_{\underset{n \neq k}{n=1}}^{+\infty} \frac{\lambda(\lambda^2+\lambda_k^2+\lambda_n^2)}{\left( \lambda^2 + (\lambda_k - \lambda_n)^2\right) \left( \lambda^2 + (\lambda_k + \lambda_n)^2 \right)} \right|^2 < + \infty.
\end{equation*}
This ends the proof of Proposition~\ref{prop_regulariteFeedback}.

\endpf

\subsection{Domain of definition and continuity of the transformation} \label{subsec_RegulariteT}

The regularity of the coefficients obtained in the previous section is sufficient to define a continuous operator $T$ in the state space.

\begin{prop} \label{Prop:regularite_T}
The transformation $T$ defined on $X^3_{(0)}$ by~(\ref{Def_T}) and~(\ref{def_alpha1n_alpha2n}) is linear continuous in $X^3_{(0)}$.
\end{prop}

\beginpf
Let $(\Psi^1, \Psi^2)\tr \in X^3_{(0)}$.
From~(\ref{Def_T}), using the Riesz basis property of Proposition~\ref{prop_baseRiesz}, it comes that $T \begin{pmatrix} \Psi^1 \\ \Psi^2 \end{pmatrix}$ belongs to $X^3_{(0)}$ if the following sequences,
\begin{gather*}
\left( \frac{\lambda_n^{1/2} \alpha^1_n \<\Psi^2,\varphi_n\>}{\beta^1_n} - \frac{\lambda_n^{1/2} \alpha^2_n \<\Psi^1,\varphi_n\>}{\beta^1_n} \right)_{n \in \NN^*},
\\
\left( \frac{\lambda_n^{1/2} \alpha^1_n \<\Psi^1,\varphi_n\>}{\beta^2_n} + \frac{\lambda_n^{1/2} \alpha^2_n \<\Psi^2,\varphi_n\>}{\beta^2_n} \right)_{n \in \NN^*},
\end{gather*}
belong to $\ell^2(\NN^*,\RR)$.
 From~(\ref{comportement_beta_jn}) and Remark~\ref{rmq_alphajn_borne}, it comes that, for every $i,j,k \in \{1,2\}$,
\begin{align*}
\sum_{n=1}^{+\infty} \left| \frac{\lambda_n^{1/2} \alpha^j_n \<\Psi^i, \varphi_n\> }{\beta^k_n} \right|^2
&\leq C \sum_{n=1}^{+\infty} \left| \frac{\alpha^j_n}{n^3} n^3 \<\Psi^i, \varphi_n\> \right|^2
\\
&\leq C \sum_{n=1}^{+\infty} \left| n^3 \<\Psi^i, \varphi_n\> \right|^2
\\
&= C \| \Psi^i \|_{H^3_{(0)}}^2.
\end{align*}
Finally, we obtain
\begin{equation*}
\left\| T \begin{pmatrix} \Psi^1 \\ \Psi^2 \end{pmatrix} \right\|_{X^3_{(0)}} \leq
C \left\| \begin{pmatrix} \Psi^1 \\ \Psi^2 \end{pmatrix} \right\|_{X^3_{(0)}}.
\end{equation*}

\endpf

\section{Invertibility of the transformation}
\label{sec_inversibilite}

This section aims at proving the invertibility of $T$. As a first step, we prove in subsection~\ref{subsec_Fredholm} that $T$ is a Fredholm operator. In subsection~\ref{subsec_Fredholm}, we prove that the analogous of~(\ref{l1}) in finite dimension holds on a certain functional space. This will be used in subsection~(\ref{subsec_Inversible}) to obtain the invertibility of $T$.

\subsection{Fredholm form} \label{subsec_Fredholm}

The goal of this subsection is the proof of the following result.
\begin{prop} \label{prop_Fredholm}
There exists $\widetilde{T} : X^3_{(0)} \to X^3_{(0)}$ invertible such that $T-\widetilde{T}$ is a compact operator.
\end{prop}

The proof of this proposition rely on the study of the feedback law done in Proposition~\ref{prop_regulariteFeedback}.

\beginpf
Let $\widetilde{T}$ be the transformation defined by~(\ref{Def_T}) where the coefficients $\alpha^1_n$ and $\alpha^2_n$ are respectively replaced by
\begin{equation*}
\tilde{\alpha}^1_n = 0, \qquad
\tilde{\alpha}^2_n =  \frac{\<h, \varphi_n\>}{\<\mu \varphi_1, \varphi_n\>} \lambda_n \beta^2_n.
\end{equation*}
From Proposition~\ref{prop_regulariteFeedback} (recall that $h_n=\<h,\varphi_n\>$ is given in~\eqref{coeff_fourier_h}), it follows that defining
$\hat{\alpha}^j_n = \alpha^j_n- \tilde{\alpha}^j_n$, we get that $( \hat{\alpha}^j_n/n^3 )_n \in \ell^2(\NN^*,\RR)$.

The computations done in the proof of Proposition~\ref{Prop:regularite_T} show that $\widetilde{T}$ is a linear continuous operator of $X^3_{(0)}$. We prove that $\widetilde{T}$ is invertible.
For any $(\Psi^1, \Psi^2)\tr \in X^3_{(0)}$,
\begin{equation*}
\widetilde{T} \begin{pmatrix} \Psi^1 \\ \Psi^2 \end{pmatrix}
= \sum_{n=1}^{+\infty} - \frac{\tilde{\alpha}^2_n \lambda_n^{1/2}}{\beta^1_n} \< \Psi^1, \varphi_n \>
\begin{pmatrix} g_n^{12}/\lambda_n^{1/2} \\ g_n^{22}/\lambda_n^{1/2} \end{pmatrix}
+ \frac{\tilde{\alpha}^2_n \lambda_n^{1/2}}{\beta^2_n}  \< \Psi^2, \varphi_n\>
\begin{pmatrix} h_n^{12}/\lambda_n^{1/2} \\ h_n^{22}/\lambda_n^{1/2} \end{pmatrix}.
\end{equation*}
 From~\eqref{boundmuphi1phik}, \eqref{deriveesdifferentes}, \eqref{comportement_beta_jn}, and  \eqref{coeff_fourier_h}, we have that, for every $j\in\{1,2\}$,  $(\beta^j_n\lambda_n^{3/2}/(\tilde{\alpha}^2_n\lambda_n^{1/2}))_n \in \ell^\infty(\NN^*,\RR)$. Then, the Riesz basis property of Proposition~\ref{prop_baseRiesz} ends the proof of the invertibility of $\widetilde{T}$.

Finally, we prove that $T-\widetilde{T}$ is compact using the Hilbert-Schmidt criterion, i.e., we prove that
\begin{equation} \label{critere_HS}
\sum_{n=1}^{+\infty} \left\| (T-\widetilde{T}) \begin{pmatrix} \varphi_n/\lambda_n^{3/2} \\ 0 \end{pmatrix} \right\|_{X^3_{(0)}}^2 + \left\| (T-\widetilde{T}) \begin{pmatrix} 0 \\ \varphi_n/\lambda_n^{3/2} \end{pmatrix} \right\|_{X^3_{(0)}}^2 <+ \infty.
\end{equation}
From~(\ref{Def_T}) and the definition of $\widetilde{T}$ it comes that
\begin{align*}
\left(T - \widetilde{T} \right) \begin{pmatrix} \Psi^1 \\ \Psi^2 \end{pmatrix}
=&\sum_{n=1}^{+\infty}  \left[ \frac{\hat{\alpha}^1_n \<\Psi^2,\varphi_n\>}{\beta^1_n} - \frac{\hat{\alpha}^2_n \<\Psi^1,\varphi_n\>}{\beta^1_n} \right]
\begin{pmatrix} g_n^{12} \\ g_n^{22} \end{pmatrix}
\\
+ &\sum_{n=1}^{+\infty} \left[ \frac{\hat{\alpha}^1_n \<\Psi^1,\varphi_n\>}{\beta^2_n} + \frac{\hat{\alpha}^2_n \<\Psi^2,\varphi_n\>}{\beta^2_n} \right]
\begin{pmatrix} h_n^{12} \\ h_n^{22} \end{pmatrix}.
\end{align*}
Thus,
\begin{align*}
&\left\| (T-\widetilde{T}) \begin{pmatrix} \varphi_n/\lambda_n^{3/2} \\ 0 \end{pmatrix} \right\|_{X^3_{(0)}}^2
\\
&= \left\| - \frac{\hat{\alpha}^2_n }{\lambda_n^{3/2} \beta^1_n}
\begin{pmatrix} g_n^{12} \\ g_n^{22} \end{pmatrix}
+ \frac{\hat{\alpha}^1_n }{\lambda_n^{3/2} \beta^2_n} \begin{pmatrix} h_n^{12} \\ h_n^{22} \end{pmatrix} \right\|_{X^3_{(0)}}^2
\\
&=\left\| - \frac{\hat{\alpha}^2_n }{\lambda_n \beta^1_n}
\begin{pmatrix} g_n^{12}/\lambda_n^{1/2} \\ g_n^{22}/\lambda_n^{1/2} \end{pmatrix}
+ \frac{\hat{\alpha}^1_n }{\lambda_n \beta^2_n}  \begin{pmatrix} h_n^{12}/\lambda_n^{1/2} \\ h_n^{22}/\lambda_n^{1/2} \end{pmatrix} \right\|_{X^3_{(0)}}^2.
\end{align*}
The first term is estimated in the following way
\begin{align*}
&\left\|
 \frac{\hat{\alpha}^2_n }{\lambda_n \beta^1_n}
\begin{pmatrix} g_n^{12}/\lambda_n^{1/2} \\ g_n^{22}/\lambda_n^{1/2} \end{pmatrix}
\right\|_{X^3_{(0)}}^2
\\
&\leq
2\left\|
 \frac{\hat{\alpha}^2_n }{\lambda_n \beta^1_n}
\begin{pmatrix} g_n^{12}/\lambda_n^{1/2} - \varphi_n/\lambda_n^{3/2} \\ g_n^{22}/\lambda_n^{1/2} \end{pmatrix}
\right\|_{X^3_{(0)}}^2
+
2\left\|
 \frac{\hat{\alpha}^2_n }{\lambda_n \beta^1_n}
\begin{pmatrix}\varphi_n/\lambda_n^{3/2} \\ 0 \end{pmatrix}
\right\|_{X^3_{(0)}}^2
\\
&\leq
2 \left(\frac{\hat{\alpha}^2_n }{\lambda_n \beta^1_n} \right)^2
\left\|
\begin{pmatrix} g_n^{12}/\lambda_n^{1/2} - \varphi_n/\lambda_n^{3/2} \\ g_n^{22}/\lambda_n^{1/2} \end{pmatrix}
\right\|_{X^3_{(0)}}^2
+
2 \left(\frac{\hat{\alpha}^2_n }{\lambda_n \beta^1_n} \right)^2.
\end{align*}
We proved in Lemma~\ref{Lemme_quadratiquement_proche} that
\begin{equation*}
\sum_{n=1}^{+\infty}
\left\|
\begin{pmatrix} g_n^{12}/\lambda_n^{1/2} - \varphi_n/\lambda_n^{3/2} \\ g_n^{22}/\lambda_n^{1/2} \end{pmatrix}
\right\|_{X^3_{(0)}}^2 < + \infty.
\end{equation*}
As $(\hat{\alpha}^2_n/n^3)_n \in \ell^2(\NN^*,\RR)$, this gives, using \eqref{comportement_beta_jn}, that
\begin{equation*}
\sum_{n=1}^{+\infty}
\left\| \frac{\hat{\alpha}^2_n }{\lambda_n \beta^1_n}
\begin{pmatrix} g_n^{12}/\lambda_n^{1/2} \\ g_n^{22}/\lambda_n^{1/2} \end{pmatrix} \right\|_{X^3_{(0)}}^2 < +\infty.
\end{equation*}
The term
\begin{equation*}
\sum_{n=1}^{+\infty}
\left\|  \frac{\hat{\alpha}^1_n }{\lambda_n \beta^2_n}
\begin{pmatrix} h_n^{12}/\lambda_n^{1/2} \\ h_n^{22}/\lambda_n^{1/2} \end{pmatrix} \right\|_{X^3_{(0)}}^2,
\end{equation*}
is dealt with in the exact same way, ending the proof of Proposition~\ref{prop_Fredholm}.

\endpf

\begin{remark} \label{rmq_inversibilite_regularite_coeff}
Notice that the key point in proving the invertibility of $\widetilde{T}$ is that for any $n \in \NN^*$, $\<h, \varphi_n\> \neq 0$.

\noindent
We also underline that the compactness of $T -\widetilde{T}$ comes from the fact that $( \hat{\alpha}^j_n/n^3)_{n \in \NN^*} \in \ell^2(\NN^*,\RR)$. If we had $\mu \varphi_1 \in H^3_{(0)}$ (this is not possible due to Hypothesis~\ref{hyp}, see \cite[Remark 2]{BeauchardLaurent}), then from~(\ref{Def_TB=B}) we would have obtained that $( \alpha^j_n/n^3)_{n \in \NN^*} \in \ell^2(\NN^*,\RR)$. This would have led to the compactness of the transformation $T$ preventing its invertibility.
\end{remark}

\subsection{Operator equality} \label{subsec_DomainDef}

We prove that the formal operator equality
\begin{equation} \label{EgaliteOperateurs}
T(A+BK) = AT - \lambda T
\end{equation}
holds true on an appropriate functional space. Recall that $K$ is defined by~(\ref{DefK}).

\begin{remark} \label{Remarque_DomaineDef_K}
Notice that due to the regularity of the coefficients $\alpha^j_n$ obtained in Proposition~\ref{prop_regulariteFeedback}, the operator $K$ is not defined on $X^3_{(0)}$. Otherwise, taking any $(0, \psi)\tr \in X^3_{(0)}$ with $\psi \notin H^s_{(0)}$ for any $s>3$ would imply
\[
\left| K\begin{pmatrix}
0 \\ \psi
\end{pmatrix} \right| = \left| \sum_{n\in \NN^*} \alpha_n^2 \< \psi,\varphi_n\> \right|
= \left| \sum_{n\in \NN^*} \frac{\alpha_n^2}{n^3} (n^3 \< \psi,\varphi_n\>) \right| < +\infty,
\]
and therefore $(\alpha_n^2/n^3 )_{n\in \NN^*} \in l^2(\NN^*,\RR)$ which is in contradiction with Proposition \ref{prop_regulariteFeedback}.
\end{remark}
Due to the previous remark, the functional setting for~(\ref{EgaliteOperateurs}) to hold needs to be specified. Let us first deal with $K$. Proposition~\ref{prop_regulariteFeedback} implies that, for every
$(\Psi^1,\Psi^2)\tr \in H^3_{(0)}\times H^4_{(0)}$,
\begin{align*}
 \sum_{n\in \NN^*} |\alpha_n^1 \< \Psi^1, \varphi_n \> + \alpha_n^2 \<& \Psi^2, \varphi_n \> |
\leq \left(\sum_{n\in \NN^*} \left|\dfrac{\alpha_n^1}{n^3}\right|^2 \right)^{1/2}\|\Psi^1\|_{H^3_{(0)}}
+ C \sum_{n\in \NN^*}\lambda_n^{3/2} \left| \< \Psi^2, \varphi_n \> \right| \\
\leq& \left(\sum_{n\in \NN^*} \left|\dfrac{\alpha_n^1}{n^3}\right|^2 \right)^{1/2}\|\Psi^1\|_{H^3_{(0)}}
+ C \left(\sum_{n\in \NN^*} \dfrac{1}{\lambda_n}\right)^{1/2} \left(\sum_{n\in \NN^*}\lambda_n^4\left|\< \Psi^2, \varphi_n \> \right|^2\right)^{1/2} \\
\leq& C \|\Psi^1\|_{H^3_{(0)}} + C \|\Psi^2\|_{H^4_{(0)}}.
\end{align*}
This shows that $K$ is well defined and continuous on $H^3_{(0)}\times H^4_{(0)}$.

We now turn to $A+BK$. Recall that we expect
\begin{equation}
\label{valueA+BK}
(A+BK)\begin{pmatrix} \Psi^1 \\ \Psi^2 \end{pmatrix}
=\begin{pmatrix} 0 & -\Delta \\ \Delta & 0 \end{pmatrix}
\begin{pmatrix} \Psi^1 \\ \Psi^2 \end{pmatrix}
+ K\begin{pmatrix} \Psi^1 \\ \Psi^2 \end{pmatrix}
\begin{pmatrix} 0 \\ \mu\varphi_1 \end{pmatrix}.
\end{equation}
We define
\begin{align*}
D(A+BK):=\left\{ \begin{pmatrix} \Psi^1 \\ \Psi^2 \end{pmatrix}
\in (H^5\cap H^3_{(0)})\times H^5_{(0)}
\: ; \: \Delta\Psi^1+K\begin{pmatrix} \Psi^1 \\ \Psi^2 \end{pmatrix}\mu \varphi_1 \in H^3_{(0)} \right\}.
\end{align*}
and then define $A+BK$ on $ D(A+BK)$ by \eqref{valueA+BK}.
Note that $(\Psi^1 , \Psi^2 )\tr \in D(A+BK)$ if and only if $(\Psi^1 ,  \Psi^2)  \in (H^5\cap H^3_{(0)})\times H^5_{(0)}$ and
\begin{equation}\label{bord_domaine}
\Delta^2\Psi^1(x)+2 K\begin{pmatrix}
\Psi^1 \\ \Psi^2
\end{pmatrix}\mu'(x) \varphi_1'(x)=0, \quad x \in \{0,1\}.
\end{equation}

We now prove the density of $D(A+BK)$.
\begin{lem}\label{Lemme:domaine_dense}
The domain $D(A+BK)$ is dense in $X^3_{(0)}$.
\end{lem}

\beginpf
Let us prove that $D(A+BK)^\perp=\{0\}$ in $X^3_{(0)}$. Let $(\Psi^1,\Psi^2 )\tr\in D(A+BK)^\perp$.

\emph{First step: we prove that $\Psi^1=0$.}

Let $k \in \NN^*$. Consider $\phi^1:=\varphi_k$. From the asymptotic behaviour of $\alpha_n^2$ (see Proposition~\ref{prop_regulariteFeedback}), there exists $N\in \NN^*$ such that $\forall n \geq N$, $\alpha_n^2\neq 0$. Therefore, setting $\phi^2=-\alpha_k^1\varphi_n/\alpha^2_n\in H^5_{(0)}$ implies that
\[
 K\begin{pmatrix}
\phi^1 \\ \phi^2
\end{pmatrix} =  \alpha_k^1+\alpha_n^2\left(-\frac{\alpha_k^1}{\alpha^2_n}\right)=0
\]
which means that \eqref{bord_domaine} is satisfied and $( \phi^1,\phi^2 )\tr \in D(A+BK)$. Thus,
\begin{align*}
0 =
\left<\begin{pmatrix} \Psi^1 \\ \Psi^2 \end{pmatrix}, \begin{pmatrix} \varphi_k \\ -\frac{\alpha_k^1}{\alpha^2_n}\varphi_n \end{pmatrix}\right\>_{X^3_{(0)}}
=&\left<\Psi^1 , \varphi_k \right\>_{H^3_{(0)}}-\frac{\alpha_k^1}{\alpha^2_n}\left<\Psi^2 ,\varphi_n \right\>_{H^3_{(0)}}\\
=&\left<\Psi^1 , \varphi_k \right\>_{H^3_{(0)}}-\frac{\lambda_n^{3/2}\alpha_k^1}{\alpha^2_n}\left(\lambda_n^{3/2}\left<\Psi^2 ,\varphi_n \right\>\right) \\
\underset{n \to + \infty}{\longrightarrow} &\left<\Psi^1 , \varphi_k \right\>_{H^3_{(0)}}.
\end{align*}
Indeed, since $\Psi^2\in H^3_{(0)}$, it implies that $\left( n^3 \lag \Psi^2 ,\varphi_n \rag \right)_{n\geq N}\in \ell^2(\NN^*;\RR )$ and Remark \ref{behavehyp}, \eqref{comportement_beta_jn} and \eqref{Def_gamma_nj} imply that $\left(\lambda_n^{3/2}\alpha_k^1/\alpha^2_n \right)_n \in \ell^\infty(\NN^*\setminus \{1,...,N-1\};\RR)$. We conclude that
\[
\left<\Psi^1 , \varphi_k \right\>_{H^3_{(0)}}=0, \quad \forall k\in \NN^*,
\]
and thus $\Psi^1 =0$.

\medskip
\emph{Second step : we prove that $\Psi^2=0$.}

For all $k\in \NN^*$, let $\phi^2=\varphi_k$. If there exists $n\in \NN^*$ such that $\alpha_n^1\neq 0$, then setting $\phi^1=-\alpha_k^2\varphi_n/\alpha_n^1$ implies that $K\begin{pmatrix} \phi^1 \\ \phi^2 \end{pmatrix}=0$ and therefore $( \phi^1, \phi^2)\tr \in D(A+BK)$. Thus,
\[
0=\left<\begin{pmatrix} \Psi^1 \\ \Psi^2 \end{pmatrix}, \begin{pmatrix} -\frac{\alpha_k^2}{\alpha_n^1}\varphi_n \\ \varphi_k \end{pmatrix}\right\>_{X^3_{(0)}}
=\left<\Psi^2 , \varphi_k \right\>_{H^3_{(0)}}.
\]
Otherwise, for every $\phi^1\in H^3_{(0)}$,
\[
K\begin{pmatrix} \phi^1 \\ \phi^2 \end{pmatrix}=\alpha_k^2.
\]
Then we consider $\phi^1$ solution to
\begin{equation*}
\begin{cases}
-\Delta \phi^1 =\alpha_k^2(\mu \varphi_1) \\
\phi^1(0)=\phi^1(1)=0.
\end{cases}
\end{equation*}
Since $\mu \varphi_1\in H^3\cap H^1_0$, then $\phi^1 \in H^5\cap H^3_{(0)}$ and, since \eqref{bord_domaine} is satisfied, $(\phi^1, \phi^2)\tr \in D(A+BK)$. Moreover, as $\Psi^1=0$,
\[
0 = \left<\begin{pmatrix} \Psi^1 \\ \Psi^2 \end{pmatrix}, \begin{pmatrix} \phi^1 \\ \phi^2 \end{pmatrix}\right\>_{X^3_{(0)}}
=\left<\Psi^2 , \varphi_k \right\>_{H^3_{(0)}}.
\]
This proves that $( \Psi^1, \Psi^2)\tr=0$.

\endpf

Let us now turn our attention to the kernel system. More precisely, we prove the following.

\begin{prop}\label{Prop:EgaliteOperateur}
For $( \Psi^1,\Psi^2)\tr \in D(A+BK)$,
\[
T(A+BK)\begin{pmatrix} \Psi^1 \\ \Psi^2 \end{pmatrix}=(A-\lambda I)T\begin{pmatrix} \Psi^1 \\ \Psi^2 \end{pmatrix}, \quad
\text{in }X^3_{(0)}.\]
\end{prop}

\beginpf
Let $( \Psi^1,\Psi^2)\tr \in D(A+BK)$. Then
\[
(A+BK)\begin{pmatrix} \Psi^1 \\ \Psi^2 \end{pmatrix}=\begin{pmatrix} -\Delta \Psi^2 \\ \Delta \Psi^1+K\begin{pmatrix} \Psi^1 \\ \Psi^2 \end{pmatrix}\mu \varphi_1 \end{pmatrix}.
\]
Therefore
\begin{align*}
T(A+BK)\begin{pmatrix} \Psi^1 \\ \Psi^2 \end{pmatrix}
=&\sum_{n\in \NN^*} \left[\dfrac{\alpha_n^1}{\beta_n^1}\left\<\Delta \Psi^1+K\begin{pmatrix} \Psi^1 \\ \Psi^2 \end{pmatrix}\mu \varphi_1, \varphi_n \right\>+\dfrac{\alpha_n^2}{\beta_n^1}\left\<\Delta \Psi^2,\varphi_n\right\> \right]\begin{pmatrix} g_n^{12} \\ g_n^{22} \end{pmatrix}
\\
&+\left[-\dfrac{\alpha_n^1}{\beta_n^2}\left\<\Delta \Psi^2,\varphi_n\right\>+\dfrac{\alpha_n^2}{\beta_n^2}\left\<\Delta \Psi^1+K\begin{pmatrix} \Psi^1 \\ \Psi^2 \end{pmatrix}\mu \varphi_1, \varphi_n \right\> \right]\begin{pmatrix} h_n^{12} \\ h_n^{22} \end{pmatrix}
\\
=&K\begin{pmatrix} \Psi^1 \\ \Psi^2 \end{pmatrix} \sum_{n\in \NN^*}\left[ \dfrac{\alpha_n^1}{\beta_n^1} \begin{pmatrix} g_n^{12} \\ g_n^{22} \end{pmatrix}  + \dfrac{\alpha_n^2}{\beta_n^2} \begin{pmatrix} h_n^{12} \\ h_n^{22} \end{pmatrix}  \right] \left\<\mu \varphi_1,\varphi_n\right\>
\\
&-\sum_{n\in\NN^*}\lambda_n\left[\dfrac{\alpha_n^1}{\beta_n^1} \left\<\Psi^1,\varphi_n\right\>+\dfrac{\alpha_n^2}{\beta_n^1} \left\<\Psi^2,\varphi_n\right\>\right]\begin{pmatrix} g_n^{12} \\ g_n^{22} \end{pmatrix}
\\
&  +\lambda_n\left[\dfrac{\alpha_n^1}{\beta_n^2} \left\<\Psi^2,\varphi_n\right\>-\dfrac{\alpha_n^2}{\beta_n^2} \left\<\Psi^1,\varphi_n\right\>\right]\begin{pmatrix} h_n^{12} \\ h_n^{22} \end{pmatrix}
\\
=&\left( \sum_{n \in \NN^*} \alpha^1_n \lag \Psi^1,\varphi_n \rag + \alpha^2_n \lag \Psi^2, \varphi_n \rag \right)  \begin{pmatrix} 0 \\ \mu \varphi_1 \end{pmatrix}
\\
&+\sum_{n\in\NN^*}-\lambda_n\left[\dfrac{\alpha_n^1}{\beta_n^1} \left\<\Psi^1,\varphi_n\right\>+\dfrac{\alpha_n^2}{\beta_n^1} \left\<\Psi^2,\varphi_n\right\>\right]\begin{pmatrix} g_n^{12} \\ g_n^{22} \end{pmatrix}
\\
&  +\lambda_n\left[\dfrac{\alpha_n^1}{\beta_n^2} \left\<\Psi^2,\varphi_n\right\>-\dfrac{\alpha_n^2}{\beta_n^2} \left\<\Psi^1,\varphi_n\right\>\right]\begin{pmatrix} h_n^{12} \\ h_n^{22} \end{pmatrix}, \text{ in }X^2_{(0)},
\end{align*}
using the expression of $T$ in \eqref{Def_T}, the $TB=B$ condition~(\ref{Def_TB=B}) and the expression of $K$ given in (\ref{DefK}).
Moreover, by the definition of $g_n^{ij}, h_n^{ij}$ given in \eqref{Formule_gn} and by the relations \eqref{Relation_gn_hn}, we have, on one hand, using \eqref{Def_T} once more,
\begin{align*}
\left\<(A-\lambda I)T \begin{pmatrix} \Psi^1 \\ \Psi^2 \end{pmatrix}, \begin{pmatrix} \varphi_k \\ 0 \end{pmatrix}\right\>_{(L^2)^2}
=& \left\< T \begin{pmatrix} \Psi^1 \\ \Psi^2 \end{pmatrix},  A^* \begin{pmatrix} \varphi_k \\ 0 \end{pmatrix} \right\>_{(L^2)^2} - \lambda \left\< T \begin{pmatrix} \Psi^1 \\ \Psi^2 \end{pmatrix}, \begin{pmatrix} \varphi_k \\ 0 \end{pmatrix} \right\>_{(L^2)^2}
\\
=& \left\<T \begin{pmatrix} \Psi^1 \\ \Psi^2 \end{pmatrix}, \begin{pmatrix} -\lambda \varphi_k  \\ \lambda_k \varphi_k \end{pmatrix}\right\>_{(L^2)^2}
\\
=& \sum_{n\in \NN^*} \left[\frac{\alpha_n^1}{\beta_n^1}\<\Psi^2,\varphi_n\>- \frac{\alpha_n^2}{\beta_n^1}\<\Psi^1,\varphi_n\>\right] \left(\lambda_k\<g_n^{22},\varphi_k\>-\lambda\<g_n^{12},\varphi_k\>\right)
\\
&+ \left[\frac{\alpha_n^1}{\beta_n^2}\<\Psi^1,\varphi_n\>+ \frac{\alpha_n^2}{\beta_n^2}\<\Psi^2,\varphi_n\>\right] \left(\lambda_k\<h_n^{22},\varphi_k\>-\lambda\<h_n^{12},\varphi_k\>\right)
\\
=& \sum_{n\in \NN^*} \left[\frac{\alpha_n^1}{\beta_n^1}\<\Psi^2,\varphi_n\>- \frac{\alpha_n^2}{\beta_n^1}\<\Psi^1,\varphi_n\>\right] \lambda_n\<g_n^{11},\varphi_k\>
\\
&+ \left[\frac{\alpha_n^1}{\beta_n^2}\<\Psi^1,\varphi_n\>+ \frac{\alpha_n^2}{\beta_n^2}\<\Psi^2,\varphi_n\>\right] \lambda_n \<h_n^{11},\varphi_k\>
\\
=& \sum_{n\in \NN^*}  \left[\frac{\alpha_n^1}{\beta_n^2}\<\Psi^2,\varphi_n\>- \frac{\alpha_n^2}{\beta_n^2}\<\Psi^1,\varphi_n\>\right] \lambda_n\<h_n^{12},\varphi_k\>
\\
&-\left[\frac{\alpha_n^1}{\beta_n^2}\<\Psi^1,\varphi_n\>+ \frac{\alpha_n^2}{\beta_n^2}\<\Psi^2,\varphi_n\>\right] \lambda_n \<g_n^{12},\varphi_k\>
\\
=& \left\<T(A+BK)\begin{pmatrix} \Psi^1 \\ \Psi^2 \end{pmatrix} , \begin{pmatrix} \varphi_k \\ 0 \end{pmatrix} \right\>_{(L^2)^2}.
\end{align*}
On the other hand, using again~\eqref{Formule_gn}, \eqref{Relation_gn_hn} and \eqref{Def_T}
\begin{align*}
\left\<(A-\lambda I)T \begin{pmatrix} \Psi^1 \\ \Psi^2 \end{pmatrix}, \begin{pmatrix} 0 \\ \varphi_k  \end{pmatrix}\right\>_{(L^2)^2}=& \left\<T \begin{pmatrix} \Psi^1 \\ \Psi^2 \end{pmatrix}, \begin{pmatrix} -\lambda_k \varphi_k  \\ -\lambda \varphi_k \end{pmatrix}\right\>_{(L^2)^2} \\
=& \sum_{n\in \NN^*}  \left[\frac{\alpha_n^1}{\beta_n^1}\<\Psi^2,\varphi_n\>- \frac{\alpha_n^2}{\beta_n^1}\<\Psi^1,\varphi_n\>\right]
\left(-\lambda \<g_n^{22},\varphi_k\>-\lambda_k \< g_n^{12},\varphi_k\>\right) \\
&+ \left[\frac{\alpha_n^1}{\beta_n^2}\<\Psi^1,\varphi_n\>+ \frac{\alpha_n^2}{\beta_n^2}\<\Psi^2,\varphi_n\>\right] \left(  -\lambda\< h_n^{22},\varphi_k\>-\lambda_k\<h_n^{12},\varphi_k\>\right) \\
=& \sum_{n\in \NN^*}  \left[\frac{\alpha_n^1}{\beta_n^2}\<\Psi^2,\varphi_n\>- \frac{\alpha_n^2}{\beta_n^2}\<\Psi^1,\varphi_n\>\right]
\lambda_n\<h_n^{22},\varphi_k\>\\
&- \left[\frac{\alpha_n^1}{\beta_n^1}\<\Psi^1,\varphi_n\>+ \frac{\alpha_n^2}{\beta_n^1}\<\Psi^2,\varphi_n\>\right]  \lambda_n\<g_n^{22},\varphi_k\>\\
&+\sum_{n\in \NN^*}  \left[\alpha_n^1\<\Psi^1,\varphi_n\>+\alpha_n^2 \lag\Psi^2,\varphi_n\rag\right]
\<\mu \varphi_1,\varphi_k\>\\
=&\left\<T(A+BK) \begin{pmatrix} \Psi^1 \\ \Psi^2 \end{pmatrix}, \begin{pmatrix} 0 \\ \varphi_k  \end{pmatrix}\right\>_{(L^2)^2}.
\end{align*}
Indeed, \eqref{Formule_gn} and \eqref{Relation_gn_hn} imply
\begin{align*}
-\lambda_k\<h_n^{12},\varphi_k\>-\lambda\<h_n^{22},\varphi_k\>=&\lambda_n\<h_n^{21},\varphi_k\>+\beta_n^2\<\mu \varphi_1,\varphi_k\>\\
=&-\lambda_n\dfrac{\beta_n^2}{\beta_n^1}\<g_n^{22},\varphi_k\>+\beta_n^2\<\mu \varphi_1,\varphi_k\>.
\end{align*}
Consequently, by the definition of $D(A+BK)$ and by the continuity of $T$ from $X^3_{(0)}$ into itself, $T(A+BK)=(A-\lambda I)T$ holds in $X^3_{(0)}$. Notice that all the previous infinite sums are converging due to the regularity assumptions on the functions of $D(A+BK)$.

\endpf

We conclude this section by noting that, following Proposition \ref{Prop:EgaliteOperateur}, if $(\Psi^1 \\ \Psi^2)\tr \in D(A+BK)$, then $T\begin{pmatrix} \Psi^1 \\ \Psi^2 \end{pmatrix}\in X^5_{(0)}$. Indeed,
\[
AT\begin{pmatrix} \Psi^1 \\ \Psi^2 \end{pmatrix}=T(A+BK)\begin{pmatrix} \Psi^1 \\ \Psi^2 \end{pmatrix}-\lambda T\begin{pmatrix} \Psi^1 \\ \Psi^2 \end{pmatrix},
\]
where the right-hand side of the last equation is in $X^3_{(0)}$.

\subsection{Invertibility} \label{subsec_Inversible}

Let us now turn to the invertibility of $T$ and prove the following result.

\begin{prop} \label{Prop:InversibiliteT}
The operator $T$ is invertible from $X^3_{(0)}$ into itself.
\end{prop}

\beginpf
In the previous section, we have proven that $\widetilde{T}$ is invertible and $T-\widetilde{T}$ is a compact operator. Consequently, the index of $T$ is equal to zero. Thus, it is sufficient to prove that $\textrm{Ker } T^*=\{0\}$ to show the invertibility of $T$.

Let us rewrite the operator equality~\eqref{EgaliteOperateurs} under the form
\begin{equation}
\label{eqcrochetavecrho}
T(A+BK+\lambda I+ \rho I)=AT + \rho T=(A + \rho I)T,
\end{equation}
where $\rho \in \mathbb{C}$ will be chosen later.
Assume for the moment that $\rho\in \mathbb{C}$ is such that
\begin{gather}\label{A+BK+rhoIinvertible}
(A+ \rho I+BK+\lambda I) \text{ is an invertible operator from $D(A+BK)$ to $X^3_{(0)}$},
\\
\label{A+rhoIinvertible}
(A+ \rho I)\text{ is an invertible operator from $D(A)$ to $X^3_{(0)}$}.
\end{gather}
Note that here, the vector spaces $D(A+BK), D(A)$ and $X^3_{(0)}$ are complexified. Then, from \eqref{eqcrochetavecrho}, one has
\begin{equation}
\label{corchet-et-inverse}
(A+\rho I)^{-1}T=T(A+BK+\lambda I+ \rho I)^{-1}.
\end{equation}
Consider $(\chi^1, \chi^2)\tr \in \textrm{Ker } T^*$. Then, for all $( \phi^1, \phi^2)\tr \in X^3_{(0)}$, we have
\begin{equation*}
\begin{array}{rcl}
0&=&\<(A+\rho I)^{-1}T\begin{pmatrix} \phi^1 \\ \phi^2 \end{pmatrix}
 -T(A+BK+\lambda I+ \rho I)^{-1}
 \begin{pmatrix} \phi^1 \\ \phi^2 \end{pmatrix}  ,
 \begin{pmatrix} \chi^1 \\ \chi^2 \end{pmatrix} \>_{X^3_{(0)}}
 \\
&=&\<\begin{pmatrix} \phi^1 \\ \phi^2 \end{pmatrix} , T^* ((A+\rho I)^{-1})^* \begin{pmatrix} \chi^1 \\ \chi^2 \end{pmatrix} \>_{X^3_{(0)}}  -\<(A+BK+\lambda I+ \rho I)^{-1}\begin{pmatrix} \phi^1 \\ \phi^2 \end{pmatrix}  , T^* \begin{pmatrix} \chi^1 \\ \chi^2 \end{pmatrix} \>_{X^3_{(0)}}
 \\
 &=&\<\begin{pmatrix} \phi^1 \\ \phi^2 \end{pmatrix} , T^* ((A+\rho I)^{-1})^* \begin{pmatrix} \chi^1 \\ \chi^2 \end{pmatrix} \>_{X^3_{(0)}}.
\end{array}
\end{equation*}
Thus the space $ \textrm{Ker } T^*$, which is of finite dimension, is stable by $((A+\rho I)^*)^{-1}$. Hence, if
\begin{equation}\label{kerT*not0}
\textrm{Ker } T^*\not = \{0\},
\end{equation}
 $((A+\rho I)^*)^{-1}$ has an eigenfunction in $ \textrm{Ker } T^*$. This eigenfunction is also an eigenfunction of $(A^*)^{-1}$.
Thus, if \eqref{kerT*not0} holds, there exists $\nu \in \CC$ and $( \chi^1, \chi^2)\tr \in \textrm{Ker } T^*\setminus\{0\}$ such that
\begin{equation}\label{vp_adjoint}
(A^*)^{-1} \begin{pmatrix} \chi^1 \\ \chi^2 \end{pmatrix} = \begin{pmatrix} -\Delta^{-1} \chi^2 \\ \Delta^{-1} \chi^1 \end{pmatrix} = \nu \begin{pmatrix} \chi^1 \\ \chi^2 \end{pmatrix}.
\end{equation}
Hence, for every $j \in \NN^*$,
\begin{equation}
\label{calculexpresionchi}
\nu^2 \<\chi^1,\varphi_j\>=\nu\<-\Delta^{-1} \chi^2,\varphi_j\>
=\dfrac{1}{\lambda_j}\nu\<\chi^2,\varphi_j\>
=\lambda_j\<\Delta^{-1} \chi^1,\varphi_j\>
=-\dfrac{1}{\lambda_j^2}\<\chi^1,\varphi_j\>.
\end{equation}
Therefore,
\begin{equation}
\label{nu2+}
(\nu^2+\dfrac{1}{\lambda_j^2})\<\chi^1,\varphi_j\>=0.
\end{equation}
Note that $\chi^1= 0$ together with \eqref{calculexpresionchi} implies $\chi^2=0$. Hence, since $(\chi^1, \chi^2)\tr \not = 0$, $\chi^1\not = 0$, which with \eqref{nu2+} implies that  there exists one and only one $k\in \NN^*$ such that
\[
\nu=\pm i \dfrac{1}{\lambda_k}, \quad \chi^1=c_k\varphi_k, \quad c_k\in \CC \setminus \{0\} .
\]
Furthermore, from \eqref{vp_adjoint}, we obtain $\chi^2=\mp i c_k\varphi_k$.
Finally, we have, by the $TB=B$ condition~(\ref{Def_TB=B}),
\begin{align*}
\pm i \overline{c_k}\lambda_k^3\<\mu \varphi_1,\varphi_k\>=&\<\mu \varphi_1,\chi^2\>_{H^1_0,H^5_{(0)}} \\
=&\<\begin{pmatrix} 0 \\  \mu \varphi_1 \end{pmatrix},\begin{pmatrix} \chi^1 \\  \chi^2 \end{pmatrix} \>_{(H^1_{0})^2,(H^5_{(0)})^2} \\
=&\< T\begin{pmatrix} 0 \\  \mu \varphi_1 \end{pmatrix},\begin{pmatrix} \chi^1 \\  \chi^2 \end{pmatrix} \>_{(H^1_{0})^2,(H^5_{(0)})^2} \\
=&\<\begin{pmatrix} 0 \\  \mu \varphi_1 \end{pmatrix},T^*\begin{pmatrix} \chi^1 \\  \chi^2 \end{pmatrix} \>_{(H^1_0)^2,(H^5_{(0)})^2}\\
=&0.
\end{align*}
Since $\<\mu \varphi_1,\varphi_k\>\neq 0$, we conclude that $c_k$ must be zero, which implies that \eqref{kerT*not0} does not hold and therefore $\textrm{Ker }T^*=\{0\}$.

It remains to prove the existence of $\rho\in \mathbb{C}$ such that \eqref{A+BK+rhoIinvertible} and \eqref{A+rhoIinvertible} hold. Let $\kappa := \rho+\lambda$. Applying $A^{-1}$ to $A+BK+\kappa I$ yields the operator
\begin{equation}
\label{A-1A+BK+rho I}
 I+A^{-1}BK+\kappa A^{-1} :  D(A+BK) \rightarrow D(A),
\end{equation}
where $A^{-1}B=(\Delta^{-1} (\mu \varphi_1), 0 )\tr $. Let us prove that
the set of $\kappa \in \mathbb{C}$ such that $I+A^{-1}BK+\kappa A^{-1}$ is invertible from $D(A+BK)$ to $D(A)$ is non-empty.

First, if $K(A^{-1}B)\neq -1$, then the operator $I+A^{-1}BK: D(A+BK) \rightarrow D(A)$ is invertible and the proof is over. Indeed, to solve
\begin{equation}\label{I+ABKinv}
(I+A^{-1}BK)\psi=f,
\end{equation}
for any $f\in D(A)$, one applies $K$ to \eqref{I+ABKinv} ($K(\psi)$, $K(A^{-1}B)$ and $K(f)$ are well-defined in this case) leading to
\[
K(\psi)(1+K(A^{-1}B))=K(f).
\]
Since $K(A^{-1}B)\neq -1$, we use the expression of $K(\psi)$ in \eqref{I+ABKinv} to obtain
\[
\psi=f-\dfrac{A^{-1}BK(f)}{1+K(A^{-1}B)}.
\]
Suppose then that $K(A^{-1}B)=-1$. It corresponds to the case where $A^{-1}B \in D(A+BK)$. Notice that $0$ is an eigenvalue of $I+A^{-1}BK$ of algebraic multiplicity $1$. Then, from \cite{Nagy}, there exists an open set $\Omega \subset \CC $ of $0\in \CC$ such that there exist an holomorphic function $\kappa \in \Omega \mapsto \lambda(\kappa) \in \CC$ and an holomorphic function $\kappa \in \Omega \mapsto x(\kappa) \in D(A+BK)$ such that
\begin{align}
x(0)=A^{-1}B=\begin{pmatrix}
 \Delta^{-1} (\mu \varphi_1) \\ 0 \end{pmatrix}, \qquad \nonumber \\
(I+A^{-1}BK+\kappa A^{-1})x(\kappa)=\lambda(\kappa) x(\kappa). \label{vpinv}
\end{align}
If $\lambda(\kappa) \neq 0$ in a small neighborhood of $0$, then $I+A^{-1}BK+\kappa A^{-1}$ is invertible for $\kappa$ close to $0$ and the proof is over. Suppose then that $\lambda(\kappa) = 0$ in a small neighborhood of $0$. In this case, consider the power serie expansion of $x$ around $0$
\[
x=\begin{pmatrix} \Delta^{-1} (\mu \varphi_1) \\ 0 \end{pmatrix}
+\sum_{k=1}^\infty \kappa^k \begin{pmatrix} x^1_k \\ x^2_k \end{pmatrix}.
\]
Notice that since $x\in D(A+BK)$ and $A^{-1}B\in D(A+BK)$, we obtain that $(x^1_k, x^2_k)\tr \in D(A)$.
At the zeroth order, \eqref{vpinv} writes
\[
\begin{pmatrix} \Delta^{-1} (\mu \varphi_1) \\ 0 \end{pmatrix}+\begin{pmatrix} \Delta^{-1} (\mu \varphi_1) \\ 0 \end{pmatrix}K\begin{pmatrix} \Delta^{-1} (\mu \varphi_1) \\ 0 \end{pmatrix}=\begin{pmatrix} 0\\ 0 \end{pmatrix},
\]
from the hypothesis $K(A^{-1}B)=-1$. At the higher order, we have
\begin{equation}\label{vpordre}
\begin{pmatrix} x^1_k \\ x^2_k \end{pmatrix}+\begin{pmatrix} \Delta^{-1} (\mu \varphi_1) \\ 0 \end{pmatrix} K \begin{pmatrix} x^1_k \\ x^2_k \end{pmatrix} + \kappa A^{-1} \begin{pmatrix} x^1_{k-1} \\ x^2_{k-1} \end{pmatrix}=\begin{pmatrix} 0 \\ 0 \end{pmatrix},
\end{equation}
where $(x^1_{0}, x^2_{0})\tr := (\Delta^{-1} (\mu \varphi_1) , 0)\tr $. Taking $K$ of \eqref{vpordre} yields
\[
K\left(A^{-1}\begin{pmatrix} x^1_k \\ x^2_k \end{pmatrix}\right)=0, \quad \forall k\geq 0.
\]
By successively taking $A^{-1}$ and $K$ of \eqref{vpordre}, we obtain
\begin{align}\label{Kexp}
K\left(A^{-n}\begin{pmatrix} x^1_k \\ x^2_k \end{pmatrix}\right)=0, \quad \forall k\geq 0, \quad \forall n\geq 1.
\end{align}
Therefore from \eqref{vpordre} and \eqref{Kexp}
\begin{align}
K\left(A^{-1-2n} \begin{pmatrix} 0 \\  (\mu \varphi_1)\end{pmatrix} \right)=& \sum_{j\in \NN}\alpha_j^1\< (-1)^n \Delta^{-1-2n}(\mu \varphi_1),\varphi_n\> \nonumber \\
=& \sum_{j\in \NN} (-1)^n \alpha_j^1\dfrac{ \< \mu \varphi_1,\varphi_n\> }{\lambda_j^{1+2n}} \nonumber \\
=&0. \label{K1}
\end{align}
Consider the entire function
\[
H(z):=\sum_{j\in \NN}  \alpha_j^1\dfrac{ \< \mu \varphi_1,\varphi_n\>e^{-z/\lambda_j^2 }}{\lambda_j^{3}}.
\]
From \eqref{K1}, we obtain that $H^{(p)}(0)=0$ and therefore $H\equiv 0$. By letting $z \rightarrow -\infty$ and by Hypothesis \ref{hyp}, we deduce that
\[
\alpha_j^1 =0, \quad \forall \, j \geq 1.
\]
In the same fashion,
\begin{align}
K\left(A^{-2n} \begin{pmatrix} 0 \\  (\mu \varphi_1)\end{pmatrix} \right)=& \sum_{j\in \NN}\alpha_j^2\< (-1)^n \Delta^{-2n}(\mu \varphi_1),\varphi_n\> \nonumber \\
=& \sum_{j\in \NN} (-1)^n \alpha_j^1\dfrac{ \< \mu \varphi_1,\varphi_n\> }{\lambda_j^{2n}} \nonumber \\
=&0. \label{K2}
\end{align}
Consider the entire function
\[
\hat{H}(z):=\sum_{j\in \NN}  \alpha_j^2\dfrac{ \< \mu \varphi_1,\varphi_n\>e^{-z/\lambda_j^2 }}{\lambda_j^{2}}.
\]
From \eqref{K2}, we obtain that $\hat{H}^{(p)}(0)=0$ and therefore $H\equiv 0$. By letting $z \rightarrow -\infty$ and by Hypothesis \ref{hyp}, we deduce that
\[
\alpha_j^2 =0, \quad \forall \, j \geq 1.
\]
From Proposition \ref{prop_regulariteFeedback}, we know that $\alpha_j^2 \neq 0, \, \forall j \geq 1$. Hence a contradiction either with $K(A^{-1}B)=-1$, which implies the invertibility of $I+A^{-1}BK+\kappa A^{-1}$ for all $\kappa \in \CC$, or with the fact that $\lambda(\kappa)=0$ in a small neighborhood of $0$, which implies that $I+A^{-1}BK+\kappa A^{-1}$ is invertible in a small neighborhood of $0$. Since $(A+\rho I)$ has discrete eigenvalues, it is possible in those two cases to choose $\rho$ such that \eqref{A+BK+rhoIinvertible}-\eqref{A+rhoIinvertible} are satisfied.

\section{Well-posedness of the closed-loop linear system and rapid stabilization}
\label{sec_WellPosedness}

This section ends the proof of Theorem~\ref{main}. Due to Remark~\ref{Remarque_DomaineDef_K}, the feedback $K$ is not well defined for functions in $X^3_{(0)}$. In subsection~\ref{subsec_wellposed}, we give a meaning to the solution $( \Psi^1, \Psi^2 )\tr $ of the closed-loop system

\begin{equation}\label{slin-closed-loop}
\begin{cases}
\partial_t \begin{pmatrix} \Psi^1 \\ \Psi^2 \end{pmatrix} = \begin{pmatrix} 0 & -\Delta  \\ \Delta & 0 \end{pmatrix} \begin{pmatrix} \Psi^1 \\ \Psi^2 \end{pmatrix}
+ K \begin{pmatrix} \Psi^1 \\ \Psi^2 \end{pmatrix} \begin{pmatrix} 0 \\ (\mu \varphi_1)(x) \end{pmatrix},  &
(t,x) \in (0,T) \times (0,1),
\\
\Psi^1(t,0) = \Psi^1(t,1) = 0, \quad \Psi^2(t,0) = \Psi^2(t,1) = 0,
& t \in  (0,T),
\\
\Psi^1(0,x) = \Psi^1_0(x), \quad \Psi^2(0,x) = \Psi^2_0(x),
& x \in (0,1)
\end{cases}
\end{equation}
by proving that $A+BK$ generates a $C^0-$semigroup. Finally we conclude to the exponential stability using the operator equality of Proposition~\ref{Prop:EgaliteOperateur} and the invertibility of the transformation $T$.

\subsection{Well-posedness of the closed-loop linear system} \label{subsec_wellposed}

Let us first show,
\begin{prop} \label{Prop:C0semigroup}
The operator $(A+BK)$ defined on $D(A+BK)$ generates a $C^0-$semigroup on $X^3_{(0)}$. Thus, there exists a unique solution $C([0,T];X^3_{(0)})$ of \eqref{slin} with $v(t)=K\begin{pmatrix}
\psi^1(t,.) \\ \psi^2(t,.)
\end{pmatrix}$.
\end{prop}

\beginpf
We prove that $(A+BK)$ is the infinitesimal generator of a $C^0-$semigroup on $X^3_{(0)}$.

\emph{First step.} The density of $D(A+BK)$ in $X^3_{(0)}$ was proven in Lemma \ref{Lemme:domaine_dense}.

\emph{Second step.} Let us prove that $(A+BK)$ is closed. Let $(\psi^1_n,\psi^2_n)\tr \in D(A+BK)$ such that
\begin{align*}
\begin{pmatrix}
\psi^1_n \\ \psi^2_n
\end{pmatrix} \underset{n \to +\infty}{\longrightarrow}&
\begin{pmatrix}
\psi^1 \\ \psi^2
\end{pmatrix}, \textrm{ in } X^3_{(0)}, \\
(A+BK)\begin{pmatrix}\psi^1_n \\ \psi^2_n
\end{pmatrix} \underset{n \to +\infty}{\longrightarrow}&
\begin{pmatrix}
\phi^1 \\ \phi^2
\end{pmatrix}, \textrm{ in }  X^3_{(0)}.
\end{align*}
We have
\[
(A+BK)\begin{pmatrix}\psi^1_n \\ \psi^2_n
\end{pmatrix}=\begin{pmatrix} -\Delta \psi^2_n \\ \Delta \psi^1_n+K\begin{pmatrix}\psi^1_n \\ \psi^2_n
\end{pmatrix}\mu \varphi_1
\end{pmatrix}.
\]
Hence
\begin{align*}
-\Delta \psi^2_n \underset{n \to +\infty}{\longrightarrow}& -\Delta \psi^2, \textrm{ in } H^1_0 \\
-\Delta \psi^2_n \underset{n \to +\infty}{\longrightarrow}& \phi^1, \textrm{ in } H^3_{(0)},
\end{align*}
and, consequently, $-\Delta \psi^2=\phi^1$, $\psi^2\in H^5_{(0)}$ and $\psi_n^2 \underset{n \to +\infty}{\longrightarrow} \psi^2$ in $H^5_{(0)}$. Therefore,
$
K\begin{pmatrix}\psi^1 \\ \psi^2 \end{pmatrix}
$
is well-defined and
\begin{align*}
&\left| K\begin{pmatrix}\psi^1_n \\ \psi^2_n \end{pmatrix} - K\begin{pmatrix}\psi^1 \\ \psi^2 \end{pmatrix}\right| = \left| \sum_{k\in \NN^*} \alpha_k^1\<\psi^1_n-\psi^1,\varphi_k\>+\alpha_k^2\<\psi^2_n-\psi^2,\varphi_k\>\right| \\
\leq & \left( \sum_{k \in \NN^*} \left( \frac{\alpha^1_k}{k^3} \right)^2 \right)^{1/2} \|\psi^1_n-\psi^1\|_{H^3_{(0)}}+ \left( \sum_{k \in \NN^*} \left( \frac{\alpha^2_k}{k^4} \right)^2 \right)^{1/2} \|\psi^2_n-\psi^2\|_{H^4_{(0)}} \underset{n \to +\infty}{\longrightarrow} 0.
\end{align*}
We obtain
\begin{align*}
\Delta \psi_n^1+K\begin{pmatrix}\psi^1_n \\ \psi^2_n \end{pmatrix} \mu \varphi_1 \underset{n \to +\infty}{\longrightarrow} & \phi^2 \quad \textrm{ in } H^3_{(0)}, \\
\Delta \psi_n^1+K\begin{pmatrix}\psi^1_n \\ \psi^2_n \end{pmatrix} \mu \varphi_1 \underset{n \to +\infty}{\longrightarrow} & \Delta \psi^1+K\begin{pmatrix}\psi^1 \\ \psi^2 \end{pmatrix} \mu \varphi_1 \quad\textrm{ in }  H^1_0.
\end{align*}
We conclude that $\Delta \psi^1+K\begin{pmatrix}\psi^1 \\ \psi^2 \end{pmatrix} \mu \varphi_1 = \phi^2$ and therefore $ (\psi^1, \psi^2)\tr \in D(A+BK)$.

\emph{Third step.} Let us now prove the dissipativity of $(A+BK)$. Since $T$ is invertible from $X^3_{(0)}$ into itself, we define the norm $\|\cdot\|_T:=\|T \cdot\|_{X^3_{(0)}}$, which is equivalent to the $X^3_{(0)}$ norm. We denote $\< \cdot, \cdot\>_T$ the associated inner product.

Consider $(\psi^1, \psi^2)\tr \in D(A+BK)$. From Proposition \ref{Prop:EgaliteOperateur}, we have, since $T \begin{pmatrix} \psi^1 \\ \psi^2 \end{pmatrix}\in X^5_{(0)}$,
\begin{align*}
\<(A+BK)\begin{pmatrix} \psi^1 \\ \psi^2 \end{pmatrix},\begin{pmatrix} \psi^1 \\ \psi^2 \end{pmatrix} \>_T=& \<T(A+BK)\begin{pmatrix} \psi^1 \\ \psi^2 \end{pmatrix},T\begin{pmatrix} \psi^1 \\ \psi^2 \end{pmatrix} \>_{X^3_{(0)}} \\
=& \<AT\begin{pmatrix} \psi^1 \\ \psi^2 \end{pmatrix},T\begin{pmatrix} \psi^1 \\ \psi^2 \end{pmatrix} \>_{X^3_{(0)}} -\lambda \left\|\begin{pmatrix} \psi^1 \\ \psi^2 \end{pmatrix} \right\|^2_T \\
=& -\lambda \left\|\begin{pmatrix} \psi^1 \\ \psi^2 \end{pmatrix} \right\|^2_T  \leq 0.
\end{align*}

Consider now the dissipativity of $(A+BK)^*$. First, let $(\psi^1, \psi^2 )\tr \in D(A+BK)$ and $(\phi^1, \phi^2)\tr \in X^3_{(0)}$. We have
\begin{align*}
\<(A+BK)\begin{pmatrix} \psi^1 \\ \psi^2 \end{pmatrix},\begin{pmatrix} \phi^1 \\ \phi^2 \end{pmatrix} \>_T=&\<T(A+BK)\begin{pmatrix} \psi^1 \\ \psi^2 \end{pmatrix},T\begin{pmatrix} \phi^1 \\ \phi^2 \end{pmatrix} \>_{X^3_{(0)}} \\
=&\<\begin{pmatrix} \psi^1 \\ \psi^2 \end{pmatrix},(T^{-1}A^*T-\lambda I)\begin{pmatrix} \phi^1 \\ \phi^2 \end{pmatrix} \>_{T}. \\
\end{align*}
which implies that $D((A+BK)^*)=\Big\{ (\phi^1, \phi^2)\tr \in X^3_{(0)} \, \Big|  \, T\begin{pmatrix} \phi^1 \\ \phi^2 \end{pmatrix}  \in X^5_{(0)} \Big\}$. Then, for $(\phi^1, \phi^2)\tr \in D((A+BK)^*)$, we have
\begin{align*}
\<(A+BK)^*\begin{pmatrix} \phi^1 \\ \phi^2 \end{pmatrix}, \begin{pmatrix} \phi^1 \\ \phi^2 \end{pmatrix}\>_T=&\<(T^{-1}A^*T-\lambda I)\begin{pmatrix} \phi^1 \\ \phi^2 \end{pmatrix},\begin{pmatrix} \phi^1 \\ \phi^2 \end{pmatrix}  \>_{T} \\
=&\<A^*T\begin{pmatrix} \phi^1 \\ \phi^2 \end{pmatrix}, T \begin{pmatrix} \phi^1 \\ \phi^2 \end{pmatrix}  \>_{X^3_{(0)}}-\lambda \left\|T\begin{pmatrix}\phi^1 \\ \phi^2 \end{pmatrix}  \right\|^2_{X^3_{(0)}}\\
\leq& \, 0.
\end{align*}

Thus from the Lumer-Philipps theorem (see e.g. \cite[Corollary 4.4]{PazyBook}) we obtain that $A+BK$ generates a $C^0-$semigroup on $X^3_{(0)}$.

\endpf

\subsection{Proof of the rapid stabilization} \label{subsec_ExponentialStability}

This section is dedicated to the proof of the main result, the rapid stabilization stated in Theorem~\ref{main}.

\beginpf
To begin, let us assume that $(\Psi^1_0, \Psi^2_0)\tr \in D(A+BK)$. Then, from Proposition~\ref{Prop:C0semigroup}, we get that the associated solution of~(\ref{slin}) with the feedback law $v(t) = K (\Psi^1(t,.),\Psi^2(t,.)){\tr}$ is given by
\begin{equation}
\begin{pmatrix} \Psi^1(t,.) \\ \Psi^2(t,.) \end{pmatrix} = e^{t (A+BK)} \begin{pmatrix} \Psi^1_0 \\ \Psi^2_0 \end{pmatrix}.
\end{equation}
Let $\begin{pmatrix} \xi^1(t,.) \\ \xi^2(t,.) \end{pmatrix} := T \begin{pmatrix} \Psi^1(t,.) \\ \Psi^2(t,.) \end{pmatrix}$. Using the equality of operators given in Proposition~\ref{Prop:EgaliteOperateur}, it comes that
\begin{align*}
\frac{\mathrm{d}}{\mathrm{d}t} \begin{pmatrix} \xi^1(t,.) \\ \xi^2(t,.) \end{pmatrix}
&= T \frac{\mathrm{d}}{\mathrm{d} t} \begin{pmatrix} \Psi^1(t,.) \\ \Psi^2(t,.) \end{pmatrix}
\\
&= T (A+BK) \begin{pmatrix} \Psi^1(t,.) \\ \Psi^2(t,.) \end{pmatrix}
\\
&= (A T - \lambda T) \begin{pmatrix} \Psi^1(t,.) \\ \Psi^2(t,.) \end{pmatrix}
\\
&= A \begin{pmatrix} \xi^1(t,.) \\ \xi^2(t,.) \end{pmatrix} - \lambda \begin{pmatrix} \xi^1(t,.) \\ \xi^2(t,.) \end{pmatrix}.
\end{align*}
Thus, we obtain the following exponential stability of $( \xi^1(t,.), \xi^2(t,.)){\tr}$
\[
\left\| \begin{pmatrix} \xi^1(t,.) \\ \xi^2(t,.) \end{pmatrix} \right\|_{X^3_{(0)}}
\leq e^{-\lambda t} \left\| \begin{pmatrix} \xi^1_0 \\ \xi^2_0 \end{pmatrix} \right\|_{X^3_{(0)}}
= e^{-\lambda t} \left\| T \begin{pmatrix} \Psi^1_0 \\ \Psi^2_0 \end{pmatrix} \right\|_{X^3_{(0)}} .
\]
From the continuity and invertibility of $T$ in $X^3_{(0)}$ (see Proposition~\ref{Prop:InversibiliteT}) it comes that
\begin{align*}
\left\| \begin{pmatrix} \Psi^1(t,.) \\ \Psi^2(t,.) \end{pmatrix} \right\|_{X^3_{(0)}}
&\leq  \left\| T^{-1} \right\| \left\| \begin{pmatrix} \xi^1(t,.) \\ \xi^2(t,.) \end{pmatrix} \right\|_{X^3_{(0)}} \\
&\leq \left\| T^{-1} \right\| \left\| T \right\| e^{-\lambda t} \left\| \begin{pmatrix} \Psi^1_0 \\ \Psi^2_0 \end{pmatrix} \right\|_{X^3_{(0)}}.
\tag{\theequation} \addtocounter{equation}{1}
\label{stabilite_exponentielle}
\end{align*}
Finally, let $\Psi_0 \in H^3_{(0)}$. Then, $(\Re(\Psi_0) \\ \Im(\Psi_0))\tr \in X^3_{(0)}$. The stability estimate~(\ref{stabilite_exponentielle}) and the density proved in Lemma~\ref{Lemme:domaine_dense} ends the proof of Theorem~\ref{main}.

\endpf

\paragraph{Acknowledgements.}

This work has been partially carried out thanks to the support of ARCHIMEDE LabEx (ANR-11-LABX- 0033), the A*MIDEX project (ANR-11-IDEX-0001-02) funded by the ``Investissements d'Avenir'' French government program managed by the ANR, the ANR grant EMAQS No.ANR-2011-BS01-017-01, the ERC avanced grant 266907 (CPDENL) of the 7th Research Framework Programme (FP7) and the FQRNT.

\appendix
\section{Simplified Saint-Venant Equation Example}
\label{Annexe_SaintVenant}

Let us provide an explicit transformation $(T,K)$ which allows to stabilize exponentially rapidly the simplified Saint-Venant equation
\begin{equation}\label{SVSimple}
\begin{cases}
h_t+v_x=0, & (t,x)\in (0,T)\times (0,1), \\
v_t+h_x=-u(t),& (x,t)\in (0,T)\times (0,1), \\
h(t,0)=v(t,1)=0, &t\in (0,T), \\
h(0,x)=h_0(x), \quad v(0,x)=v_0(x), & x\in (0,1), \\
\end{cases}
\end{equation}
which is controllable in time $T> 2$.

Let $H:=h_x$ and $V:=v_x$. Then, the equation on $(H,V)$ writes
\begin{equation*}
\begin{cases}
H_t+V_x=0, & (t,x)\in (0,T)\times (0,1), \\
V_t+H_x=0,& (t,x)\in (0,T)\times (0,1), \\
H(t,1)=-u(t), &t\in (0,T), \\
V(t,0)=0, &t\in (0,T), \\
H(0,x)=H_0(x), \quad V(0,x)=V_0(x), & x\in (0,1), \\
\end{cases}
\end{equation*}
with $H_0=(h_0)_x$ and $V_0=(v_0)_x$. Consider now $R^1:=H+V$ and $R^2:=H-V$. Then
\begin{equation*}
\begin{cases}
R^1_t+R^1_x=0, & (x,t)\in (0,T)\times (0,1), \\
R^2_t-R^2_x=0,& (x,t)\in (0,T)\times (0,1), \\
(R^1+R^2)(t,1)=-2u(t), &t\in (0,T), \\
(R^1-R^2)(t,0)=0, &t\in (0,T), \\
R^1(0,x)=R^1_0(x), \quad R^2(0,x)=R^2_0(x), & x\in (0,1),
\end{cases}
\end{equation*}
with $R^1_0:=H_0+V_0$ and $R^2_0:=H_0-V_0$. Let us consider a transformation which maps $(R^1,R^2)$ to a solution of a target stable system, that is, $\widetilde{R}^1:=e^{-\lambda x}R^1/\cosh(\lambda)$ and $\widetilde{R}^2:=e^{\lambda x}R^2/\cosh(\lambda)$, for $\lambda>0$. A straightforward computation leads to
\begin{equation}\label{stableSV}
\begin{cases}
\widetilde{R}^1_t+\widetilde{R}^1_x+\lambda\widetilde{R}^1=0, & (x,t)\in (0,T)\times (0,1), \\
\widetilde{R}^2_t-\widetilde{R}^2_x+\lambda\widetilde{R}^2=0,& (x,t)\in (0,T)\times (0,1), \\
(\widetilde{R}^1+\widetilde{R}^2)(t,1)=-2e^{-\lambda}(u(t)/\cosh(\lambda))+2\tanh(\lambda)R^2(t,1), &t\in (0,T) \\
(\widetilde{R}^1-\widetilde{R}^2)(t,0)=0, &t\in (0,T), \\
\widetilde{R}^1(0,x)=\widetilde{R}^1_0(x), \quad \widetilde{R}^2(0,x)=\widetilde{R}^2_0(x), & x\in (0,1).
\end{cases}
\end{equation}
Hence, the exponential stability of \eqref{stableSV} is obtained if $-2e^{-\lambda} u(t)/\cosh(\lambda)+2\tanh(\lambda)R^2(t,1)=0$. In terms of the original variables, it implies that
\begin{align*}
0=&-2e^{-\lambda} u(t)/\cosh(\lambda)+2\tanh(\lambda)R^2(t,1) \\
=&-2e^{-\lambda} u(t)/\cosh(\lambda)+2\tanh(\lambda)(h_x(t,1)-v_x(t,1)) \\
=&-2e^{-\lambda} u(t)/\cosh(\lambda)-2\tanh(\lambda)(v_x(t,1)+u(t)),
\end{align*}
that is, $u(t)=-\tanh(\lambda)v_x(t,1)$.

The target system of \eqref{SVSimple} is given by
\begin{equation}\label{SVSimplestable}
\begin{cases}
\widetilde{h}_t+\widetilde{v}_x+\lambda \widetilde{h}=0, & (x,t)\in (0,T)\times (0,1), \\
\widetilde{v}_t+\widetilde{h}_x + \lambda \widetilde{v}=0,& (x,t)\in (0,T)\times (0,1), \\
\widetilde{h}(t,0)=\widetilde{v}(t,1)=0, &t\in (0,T), \\
\widetilde{h}(0,x)=\widetilde{h}_0(x), \, \widetilde{v}(0,x)=\widetilde{v}_0(x), & x\in (0,1), \\
\end{cases}
\end{equation}
One can recover an explicit expression of the transformation $(T,K)$ leading to this target system using
\begin{equation} \label{Def_T_StVenant}
\begin{aligned}
\widetilde{h}_x+\widetilde{v}_x=&e^{-\lambda x}(h_x+v_x)/\cosh(\lambda), \\
\widetilde{h}_x-\widetilde{v}_x=&e^{\lambda x}(h_x-v_x)/\cosh(\lambda),
\end{aligned}
\end{equation}
which boils down to
\begin{align*}
\widetilde{h}_x=&\left(\cosh(\lambda x)h_x -\sinh(\lambda x)v_x\right)/\cosh(\lambda) \\
\widetilde{v}_x=&\left(-\sinh(\lambda x)h_x+\cosh(\lambda x)v_x\right)/\cosh(\lambda).
\end{align*}
Using the boundary conditions, one obtains the explicit transformation $T$
\begin{align*}
\widetilde{h}(x)=&\dfrac{1}{\cosh(\lambda)}\left[\cosh(\lambda x)h(x)-\lambda \int_0^x \sinh(\lambda y)h(y)\, dy -\sinh(\lambda x)v(x)+\lambda \int_0^x \cosh(\lambda y)v(y)\, dy\right] \\
=&\dfrac{1}{\cosh(\lambda)}\left[ \int_0^1 \left(\delta_{x=y}\cosh(\lambda y)-\lambda \mathds{1}_{(0,x)}(y)\sinh(\lambda y)\right)h(y)\, dy \right. \\
&+ \left. \int_0^1 \left(\lambda \mathds{1}_{(0,x)}(y)\cosh(\lambda y)-\delta_{x=y}\sinh(\lambda y)\right)v(y)\, dy\right], \\
\widetilde{v}(x)=&\dfrac{1}{\cosh(\lambda)}\left[\sinh(\lambda) h(1) - \sinh(\lambda x) h(x)-\int_x^1\lambda \cosh(\lambda y)h(y) \, dy \right. \\
 &+ \left. \cosh(\lambda x) v(x) + \lambda \int_x^1 \sinh(\lambda y)v(y)\, dy\right] \\
=&\dfrac{1}{\cosh(\lambda)}\left[ \int_0^1 \left( \delta_{y=1}\sinh(\lambda) - \delta_{x=y}\sinh(\lambda y)-\lambda \mathds{1}_{(x,1)}(y)\cosh(\lambda y)\right)h(y)\, dy \right. \\
&+\left. \int_0^1 \left(\lambda \mathds{1}_{(x,1)}(y)\sinh(\lambda y)+\delta_{x=y}\cosh(\lambda y)\right)v(y)\, dy\right].
\end{align*}
Moreover, the explicit transformation $K$ writes
\[
u(t)=\tanh(\lambda)\int_0^1 \delta'_{y=1}v(t,y)\, dy.
\]
If one writes, in the same spirit as \eqref{SystemeNoyau}, the kernels equation for \eqref{SVSimple}, then one obtains that the kernels of the transformations $(T,K)$ exhibited here are the solution of this system. One also verifies that, thanks to the factor $1/\cosh(\lambda)$, the $TB=B$ condition is verified by the transformation $T$. Getting back to~\eqref{Def_T_StVenant} one sees that the inverse of $T$ can be computed explicitly performing similar computations.

Moreover, the Fourier coefficients of the kernels system associated to \eqref{SVSimple} have the same expression as \eqref{Formule_gn}, where the eigenvalues/eigenfunctions are replaced by those associated with \eqref{SVSimple} and the Fourier coefficients of the control operator $\mu \varphi_1$ are replaced by the one of \eqref{SVSimple}, that is $1$. 

It is also noticeable that the Fourier coefficients of the kernel $\alpha^2$ are
\[
\alpha_n^2=(-1)^n (\pi n)\tanh(\lambda),
\]
which is adequate with the perturbation argument used in Proposition \ref{prop_regulariteFeedback} to obtain the Fourier coefficients from the $TB=B$ condition. One notice that, for \eqref{SVSimple}, $\alpha^1\equiv 0$. Whether $\alpha^1\equiv 0$ or not in the case of the linearized Schr\"odinger equation cannot be verified with our analysis.
\section{Quadratically close families}
\label{Annexe_Riesz}

This section is devoted to the proof of Lemma~\ref{Lemme_quadratiquement_proche}.

\beginpf
Let $s=2$ or $3$.
To simplify the notations, let,
\[
\tilde{c}^{ij}_{nk}:=c^{ij}_{nk}\beta_n^1\<\mu \varphi_1,\varphi_k \>, \quad \tilde{d}^{ij}_{nk}:=d^{ij}_{nk}\beta_n^2 \<\mu \varphi_1,\varphi_k \>,
\]
where $c_{nk}^{ij}$ and $d_{nk}^{ij}$ are defined by \eqref{DecompositionFourier_gn_hn} and $\beta_n^j$ are defined by~\eqref{Def_beta_jn}.

\medskip
\emph{First step: } let us prove that
\[
\sum_{n\in \NN^*} \left\|
\left(
\begin{array}{c}
\varphi_n/\lambda_n^{s/2}\\
0
\end{array}
\right)
-
\left(
\begin{array}{c}
g^{12}_n/\lambda_n^{(s-2)/2} \\
g^{22}_n/\lambda_n^{(s-2)/2}
\end{array}
\right)
\right\|^2_{{X^{s}_{(0)}}} <+\infty.
\]
First, by denoting $k=m+n$, we have using in particular \eqref{esthyp1}, \eqref{boundmuphi1phik}, \eqref{Formule_gn} and \eqref{Relation_gn_hn}
\begin{align*}
&\sum_{n\in \NN^*} \left\| \frac{\varphi_n}{\lambda_n^{s/2}}-\frac{g^{12}_n}{\lambda_{n}^{(s-2)/2}} \right\|^2_{H^{s}_{(0)}}
= \sum_{n\in \NN^*} \sum_{k\in \NN^*\setminus \{n\} } \left|\dfrac{ \lambda_k^{s/2} \tilde{c}^{12}_{nk}}{\lambda_{n}^{(s-2)/2}} \right|^2
\\
&=\sum_{n\in \NN^*}\bigg(\sum_{\substack{0<|m|<n \\ m\in \ZZ}} + \sum_{\substack{n<m \\ m\in \NN^*}} \bigg) \frac{\lambda_{n+m}^s}{\lambda_n^s} \left| \lambda_n \tilde{c}^{12}_{nn+m} \right|^2
\\
&= \lambda^4 \sum_{n\in \NN^*} \bigg(\sum_{\substack{0<|m|<n \\ m\in \ZZ}} + \sum_{\substack{n<m \\ m\in \NN^*}} \bigg) \frac{\lambda_{n+m}^s}{\lambda_n^s}
\left| \frac{(\lambda^2+4\lambda_n^2) \lambda_{n+m}}{\delta_{nn+m}(\lambda) \lambda_n} \right|^2 \left| \frac{\< \mu \varphi_1, \varphi_{n+m} \>}{\< \mu \varphi_1, \varphi_n \>} \right|^2
\\
&\leq C \lambda^4 \sum_{n\in \NN^*} \bigg(\sum_{\substack{0<|m|<n \\ m\in \ZZ}} + \sum_{\substack{n<m \\ m\in \NN^*}} \bigg) \left(\frac{\lambda_{n+m}}{\lambda_n} \right)^{s-1}
\left| \frac{\lambda^2+4\lambda_n^2}{\delta_{nn+m}(\lambda)} \right|^2.
\tag{\theequation} \addtocounter{equation}{1}
\label{QuadProche1}
\end{align*}
The two sums of~\eqref{QuadProche1} are dealt with separately.

Consider first the case where $0 < |m| <n$ and $m \in \ZZ$.
We have, using in particular \eqref{defdeltank}
\begin{align}
\notag
\frac{(\lambda^2+4\lambda_n^2)}{\delta_{nn+m}(\lambda)}
&= \frac{1}{m^2 n^2} \frac{\frac{\lambda^2}{n^4} + 4\pi^4}{\left( \frac{\lambda^2}{n^4} +
 \left( \pi^2 + \pi^2 \left( 1 +\frac{m}{n}\right) \right)^2  \right) \left( \frac{\lambda^2}{(mn)^2} + \left( \frac{m}{n} \pi^2 + 2\pi^2 \right)^2 \right)}
\\
&\leq \frac{1}{m^2 n^2} \frac{\lambda^4+4\pi^4}{\pi^8},
\label{QuadProche_delta_m/n<1}
\end{align}
and
\begin{equation} \label{QuadProche_lambda_n+m/lambda_n_m/n<1}
\lambda_{n+m} = \left( 1 +\frac{m}{n} \right)^2 \lambda_n \leq 4 \lambda_n.
\end{equation}
Thus, for the first term of the right-hand side of~\eqref{QuadProche1}
\begin{equation} \label{QuadProche1_1}
\sum_{n\in \NN^*} \sum_{\substack{0<|m|<n \\ m\in \ZZ}}
\left(\frac{\lambda_{n+m}}{\lambda_n} \right)^{s-1}
\left| \frac{(\lambda^2+4\lambda_n^2)}{\delta_{nn+m}(\lambda)} \right|^2
\leq C \sum_{n\in \NN^*} \sum_{\substack{0<|m|<n \\ m\in \ZZ}} \frac{1}{m^4 n^4}
< +\infty.
\end{equation}

Consider now the case where $m > n$.
We have
\begin{align}
\notag
\frac{(\lambda^2+4\lambda_n^2)}{\delta_{nn+m}(\lambda)}
&=
\frac{n^4}{m^8} \frac{\frac{\lambda^2}{n^4} + 4\pi^4}{\left( \frac{\lambda^2}{m^4} +  \left(  \left( 1+\frac{n}{m} \right)^2 \pi^2 + \left( \frac{n}{m} \right)^2 \pi^2  \right)^2 \right)  \left( \frac{\lambda^2}{m^4} + \left( \pi^2 + 2 \frac{n}{m} \pi^2 \right)^2 \right)},
\\
&\leq \frac{n^4}{m^8} \frac{\lambda^2+4\pi^4}{\pi^8}
\label{QuadProche_delta_m>n}
\end{align}
and
\begin{equation} \label{QuadProche_lambda_n+m/lambda_n_m>n}
\frac{\lambda_{n+m}}{\lambda_n} = \frac{(n+m)^2}{n^2} = \frac{m^2}{n^2} \left(1 +\frac{n}{m} \right)^2 \leq 4 \frac{m^2}{n^2}.
\end{equation}
Thus, for the second term of the right-hand side of~\eqref{QuadProche1} yields,
\begin{align}
\notag
\sum_{n\in \NN^*} \sum_{\substack{n < m \\ m\in \NN^*}}
\left(\frac{\lambda_{n+m}}{\lambda_n} \right)^{s-1}
\left| \frac{(\lambda^2+4\lambda_n^2)}{\delta_{nn+m}(\lambda)} \right|^2
&\leq C \sum_{n\in \NN^*} \sum_{\substack{n<m \\ m\in \NN^*}} \left(\frac{n}{m}\right)^{8-2(s-1)} \frac{1}{m^8}
\\
&\leq C \sum_{n\in \NN^*} \sum_{\substack{n<m \\ m\in \NN^*}} \frac{1}{m^4} \frac{1}{n^4}
< +\infty.
\label{QuadProche1_2}
\end{align}
Then inequalities~\eqref{QuadProche1}, \eqref{QuadProche1_1} and~\eqref{QuadProche1_2} imply that
\begin{equation} \label{QuadProcheConclu1}
\sum_{n\in \NN^*} \left\| \frac{\varphi_n}{\lambda_n^{s/2}}-\frac{g^{12}_n}{\lambda_{n}^{(s-2)/2}} \right\|^2_{H^{s}_{(0)}} <+\infty.
\end{equation}

\medskip
The other sum is treated in a similar way. Indeed,
\begin{align*}
&\sum_{n\in \NN^*} \left\|\dfrac{g^{22}_n}{\lambda_n^{(s-2)/2}}\right\|^2_{H^s_{(0)}}
= \sum_{n\in \NN^*}\sum_{k\in \NN^*} \left|\dfrac{\lambda_k^{s/2}\tilde{c}^{22}_{nk}}{\lambda_n^{(s-2)/2}}\right|^2 \\
&=\sum_{n \in \NN^*} \bigg( \sum_{\substack{0<|m|<n \\ m\in \ZZ}} + \sum_{\substack{n<m \\ m\in \NN^*}} \bigg) \left( \frac{\lambda_{n+m}}{\lambda_n} \right)^s \left| \lambda_n \tilde{c}_{nn+m}^{22} \right|^2
+ \sum_{n\in \NN^*} \left| \lambda_n \tilde{c}_{nn}^{22} \right|^2
\\
&= \frac{\lambda^2}{4} \sum_{n \in \NN^*} \bigg( \sum_{\substack{0<|m|<n \\ m\in \ZZ}} + \sum_{\substack{n<m \\ m\in \NN^*}} \bigg)
\left( \frac{\lambda_{n+m}}{\lambda_n} \right)^s
\left| \frac{\lambda^2+4 \lambda_n^2}{\delta_{nn+m}(\lambda)} \right|^2
\left| \frac{\lambda^2-\lambda_{n+m}^2+ \lambda_n^2}{\lambda_n} \right|^2
\left| \frac{\< \mu \varphi_1, \varphi_{n+m}\>}{\< \mu \varphi_1, \varphi_n\>} \right|^2
\\
&+ \sum_{n\in \NN^*} \left( \frac{\lambda}{2\lambda_n} \right)^2
\\
&\leq C + C \sum_{n \in \NN^*} \bigg( \sum_{\substack{0<|m|<n \\ m\in \ZZ}} + \sum_{\substack{n<m \\ m\in \NN^*}} \bigg)
\left( \frac{\lambda_{n+m}}{\lambda_n} \right)^{s-3}
\left| \frac{\lambda^2+4 \lambda_n^2}{\delta_{nn+m}(\lambda)} \right|^2
\left| \frac{\lambda^2-\lambda_{n+m}^2+ \lambda_n^2}{\lambda_n} \right|^2.
\tag{\theequation} \addtocounter{equation}{1}
\label{QuadProche2}
\end{align*}

Consider first the case where $0 < |m| <n$.
We have
\begin{equation*}
\left| \frac{\lambda^2-\lambda_{n+m}^2+ \lambda_n^2}{\lambda_n} \right|
\leq C \frac{n^3 |m|}{n^2} \left| \frac{\lambda^2}{n^3 m} - \pi^4 \left( 4 + 6 \frac{m}{n} + 4 \left( \frac{m}{n} \right)^2 + \left(\frac{m}{n} \right)^3 \right) \right|
\leq C n |m| \left( \lambda^2 + 15\pi^4 \right).
\end{equation*}
Using~\eqref{QuadProche_delta_m/n<1}, it comes that
\begin{equation*}
\sum_{n \in \NN^*} \sum_{\substack{0<|m|<n \\ m\in \ZZ}}
\left( \frac{\lambda_{n+m}}{\lambda_n} \right)^{s-3}
\left| \frac{\lambda^2+4 \lambda_n^2}{\delta_{nn+m}(\lambda)} \right|^2
\left| \frac{\lambda^2-\lambda_{n+m}^2+ \lambda_n^2}{\lambda_n} \right|^2
\leq C \sum_{n \in \NN^*} \sum_{\substack{0<|m|<n \\ m\in \ZZ}}
\left( \frac{\lambda_{n+m}}{\lambda_n} \right)^{s-3} \frac{1}{n^2 m^2}.
\end{equation*}
This sum is clearly finite for $s=3$. For $s=2$, we notice that the previous series is the general term of a Cauchy product. Indeed,
\begin{align*}
\sum_{\substack{0<|m|<n \\ m\in \ZZ}}
\left( \frac{\lambda_{n+m}}{\lambda_n} \right)^{s-3} \frac{1}{n^2 m^2}
&= \sum_{\substack{0<|m|<n \\ m\in \ZZ}} \frac{1}{m^2} \frac{1}{(n+m)^2}
\\
&\leq 2\sum_{m=1}^{n-1} \frac{1}{m^2} \frac{1}{(n-m)^2} + \frac{2}{n^2} \sum_{m \in \NN^*} \frac{1}{m^2}.
\end{align*}
Thus, for $s\in \{2,3\}$,
\begin{equation} \label{QuadProche2_1}
\sum_{n \in \NN^*} \sum_{\substack{0<|m|<n \\ m\in \ZZ}}
\left( \frac{\lambda_{n+m}}{\lambda_n} \right)^{s-3}
\left| \frac{\lambda^2+4 \lambda_n^2}{\delta_{nn+m}(\lambda)} \right|^2
\left| \frac{\lambda^2-\lambda_{n+m}^2+ \lambda_n^2}{\lambda_n} \right|^2 < +\infty.
\end{equation}

Consider now the case where $m>n$.
We have
\begin{equation*}
\left| \frac{\lambda^2-\lambda_{n+m}^2+ \lambda_n^2}{\lambda_n} \right|
= \left| \frac{m^4 \left( \frac{\lambda^2}{m^4} - \left(1+\frac{n}{m}\right)^4 \pi^4 + \left(\frac{n}{m}\right)^4  \right)}{n^2 \pi^2} \right|
\leq C \frac{m^4}{n^2},
\end{equation*}
and
\begin{equation*}
\left( \frac{\lambda_{n+m}}{\lambda_n} \right)^{s-3} \leq 1.
\end{equation*}
Using~\eqref{QuadProche_delta_m>n}, it then comes that,
\begin{align}
\notag
\sum_{n \in \NN^*} \sum_{\substack{n<m \\ m\in \NN^*}}
\left( \frac{\lambda_{n+m}}{\lambda_n} \right)^{s-3}
\left| \frac{\lambda^2+4 \lambda_n^2}{\delta_{nn+m}(\lambda)} \right|^2
\left| \frac{\lambda^2-\lambda_{n+m}^2+ \lambda_n^2}{\lambda_n} \right|^2
&\leq C \sum_{n \in \NN^*} \sum_{\substack{n<m \\ m\in \NN^*}} \left( \frac{m^4}{n^2} \frac{n^4}{m^8} \right)^2
\\
&\leq C \sum_{n \in \NN^*} \sum_{\substack{n<m \\ m\in \NN^*}} \frac{1}{m^2} \frac{1}{n^2}
< +\infty.
\label{QuadProche2_2}
\end{align}
Then, inequalities~\eqref{QuadProche2}, \eqref{QuadProche2_1} and~\eqref{QuadProche2_2} imply that
\begin{equation} \label{QuadProcheConclu2}
\sum_{n\in \NN^*} \left\| \frac{g^{22}_n}{\lambda_{n}^{(s-2)/2}} \right\|^2_{H^{s}_{(0)}} <+\infty.
\end{equation}
Together with~\eqref{QuadProcheConclu1} it ends the first step.

\medskip
\emph{Second step: } let us prove that
\[
\sum_{n\in \NN^*} \left\|
\left(
\begin{array}{c}
0 \\
\varphi_n/\lambda_n^{s/2}
\end{array}
\right)
-
\left(
\begin{array}{c}
h^{12}_n/\lambda_n^{(s-2)/2}\\
h^{22}_n/\lambda_n^{(s-2)/2}
\end{array}
\right)
\right\|^2_{{X^s_{(0)}}} <+\infty.
\]

This proof is very similar to the first step. Thus we give the expressions of the different sums but we do not detail every computation:
\begin{align*}
&\sum_{n\in \NN^*} \left\|\dfrac{h^{12}_n}{\lambda_n^{(s-2)/2}}\right\|^2_{H^s_{(0)}}
= \sum_{n\in \NN^*}\sum_{k\in \NN^*} \left|\dfrac{\lambda_k^{s/2}\tilde{d}^{12}_{nk}}{\lambda_n^{(s-2)/2}}\right|^2
\\
&=\sum_{n \in \NN^*} \bigg( \sum_{\substack{0<|m|<n \\ m\in \ZZ}} + \sum_{\substack{n<m \\ m\in \NN^*}} \bigg) \left( \frac{\lambda_{n+m}}{\lambda_n} \right)^s \left| \lambda_n \tilde{d}_{nn+m}^{12} \right|^2
+ \sum_{n\in \NN^*} \left| \lambda_n \tilde{d}_{nn}^{12} \right|^2
\\
&= \lambda^2 \sum_{n \in \NN^*} \bigg( \sum_{\substack{0<|m|<n \\ m\in \ZZ}} + \sum_{\substack{n<m \\ m\in \NN^*}} \bigg)
\frac{\lambda_{n+m}^s}{\lambda_n^s}
\left| \frac{\lambda^2+4\lambda_n^2}{\delta_{nn+m}(\lambda)} \right|^2
\left| \frac{\lambda_{n+m} \left(\lambda^2+\lambda_{n+m}^2-\lambda_n^2\right)}{\lambda^2+2\lambda_n^2} \right|^2
\left| \frac{\< \mu \varphi_1, \varphi_{n+m}\>}{\< \mu \varphi_1, \varphi_n \>} \right|^2
\\
&+ \sum_{n \in \NN^*} \left| \frac{\lambda^2 \lambda_n}{\lambda^2+\lambda_n^2} \right|^2
\\
&< +\infty,
\end{align*}
and
\begin{align*}
&\sum_{n\in \NN^*} \left\| \frac{\varphi_n}{\lambda_n^{s/2}}-\frac{h^{22}_n}{\lambda_{n}^{(s-2)/2}} \right\|^2_{H^{s}_{(0)}}
= \sum_{n\in \NN^*} \sum_{k\in \NN^*\setminus \{n\} } \left|\dfrac{ \lambda_k^{s/2} \tilde{d}^{22}_{nk}}{\lambda_{n}^{(s-2)/2}} \right|^2
\\
&=\sum_{n\in \NN^*}\bigg(\sum_{\substack{0<|m|<n \\ m\in \ZZ}} + \sum_{\substack{n<m \\ m\in \NN^*}} \bigg) \frac{\lambda_{n+m}^s}{\lambda_n^s}
\left| \lambda_n \tilde{d}^{22}_{nn+m} \right|^2
\\
&=\lambda^4 \sum_{n\in \NN^*}\bigg(\sum_{\substack{0<|m|<n \\ m\in \ZZ}} + \sum_{\substack{n<m \\ m\in \NN^*}} \bigg) \frac{\lambda_{n+m}^s}{\lambda_n^s}
\left| \frac{\lambda^2+4\lambda_n^2}{\delta_{nn+m}(\lambda)} \right|^2
\left| \frac{\lambda^2+\lambda_{n+m}^2+\lambda_n^2}{\lambda^2+2\lambda_n^2} \right|^2
\left| \frac{\< \mu \varphi_1, \varphi_{n+m}\>}{\< \mu \varphi_1, \varphi_n \>} \right|^2
\\
&<+\infty.
\end{align*}
This completes the proof of Lemma~\ref{Lemme_quadratiquement_proche}.

\endpf

\section{Rapid stabilization of the linearized system}
\label{Annexe_lambda1}

In this section we detail how Theorem~\ref{main_linearise} can be obtained from the results developed in this article.

Let $\Psi_0 \in \HH_0$ and $\Psi$ the associated solution for a control $u$. Then, if we define,
$\widetilde{\Psi} (t, \cdot) := \Psi(t,\cdot) e^{i \lambda_1 t}$, it comes that $\widetilde{\Psi}$ satisfies
\begin{equation} \label{linearise_lambda1}
\begin{cases}
 i \partial_t \widetilde{\Psi} = - \Delta \widetilde{\Psi} - \lambda_1 \widetilde{\Psi} - u(t) \mu \varphi_1, & (t,x)\in (0,T)\times (0,1)
\\
 \widetilde{\Psi}(t,0) = \widetilde{\Psi}(t,1) = 0,& t\in (0,T)
\\
\widetilde{\Psi}(0,\cdot) = \Psi_0, & x \in (0,1)
\end{cases}
\end{equation}
and $\| \widetilde{\Psi}(t,.) \|_{H^3_{(0)}} = \| \Psi(t,.) \|_{H^3_{(0)}}$.

\paragraph{Rapid stabilization of~(\ref{linearise_lambda1}).} Notice that system~(\ref{linearise_lambda1}) is almost identical to (\ref{SystLin}) except that the spectrum of the underlying operator is shifted by $\lambda_1$. This modifies the state space. Indeed, for every $t \geq 0$,
\[
0 = \Re \< \Psi(t,.) , \Phi_1(t,.) \> = \Re \< \widetilde{\Psi}(t,.) , \varphi_1 \>.
\]
Thus, $\widetilde{\Psi} \in \HH_0$. Notice that due to Theorem~\ref{KC}, one gets that system~(\ref{linearise_lambda1}) is exactly controllable in $\HH_0$. The rest of the analysis is barely modified.
For instance, the Riesz bases have instead the form
\[
\B = \left\{ \begin{pmatrix} g_n^{12} \\ g_n^{22} \end{pmatrix} \: ; \: n \geq 2 \right\} \cup \left\{ \begin{pmatrix} h_n^{12} \\ h_n^{22} \end{pmatrix} \: ; \: n \geq 1 \right\},
\]
where, for example,
\[
g_n^{12} (x) = \sum_{k=1}^{+\infty} \frac{2 \lambda (\lambda_k - \lambda_1) (\lambda_n-\lambda_1)}{(\lambda^2 + (\lambda_k -\lambda_n)^2)(\lambda^2 + (\lambda_k + \lambda_n - 2 \lambda_1)^2)} \beta^1_n \< \mu \varphi_1, \varphi_k \> \varphi_k(x),
\]
with $\beta^1_n$ chosen such that $\< g_n^{12} , \varphi_n \> = 1/\lambda_n$.

This leads to the existence of a feedback law $\tilde{K} \begin{pmatrix} \Re( \widetilde{\Psi}) \\ \Im (\widetilde{\Psi}) \end{pmatrix}$ such that, for the closed-loop system
\[
\left\| \begin{pmatrix} \Re (\widetilde{\Psi}(t,.)) \\ \Im (\widetilde{\Psi}(t,.)) \end{pmatrix} \right\|_{X^3_{(0)}} \leq C e^{-\lambda t} \left\| \begin{pmatrix} \Re (\Psi_0) \\ \Im (\Psi_0) \end{pmatrix} \right\|_{X^3_{(0)}}.
\]

\paragraph{Rapid stabilization of~(\ref{lin}).}
Due to the previous relation, $\widetilde{\Psi}(t,.) = \Psi(t,.) e^{i \lambda_1 t}$ it comes that if $\Psi$ is the solution of~(\ref{lin}) with the feedback law
\[
v(t) = \tilde{K} \begin{pmatrix} \cos(\lambda_1 t) \Re (\Psi(t,.)) - \sin(\lambda_1 t) \Im (\Psi(t,.)) \\\sin(\lambda_1 t) \Re (\Psi(t,.)) + \cos(\lambda_1 t) \Im (\Psi(t,.)) \end{pmatrix},
\]
then
\[
\| \Psi(t,.) \|_{H^3_{(0)}} \leq C e^{-\lambda t} \| \Psi_0 \|_{H^3_{(0)}}.
\]

\bibliographystyle{plain}
\bibliography{biblio}

\end{document}